\documentclass[10pt,preprint,reqno]{amsart}
\usepackage{amsmath,amssymb,amsfonts,amsthm}

\usepackage{amscd, calc, tikz, tikz-cd}
\usepackage[all]{xy}
\usepackage[color=pink,textsize=footnotesize]{todonotes}

\usepackage{enumerate,graphicx,psfrag}
\usepackage{mathrsfs}
\usepackage{hyperref}

\usepackage{color} 
\usepackage{pstricks,pdftricks} 

\usepackage[shortlabels]{enumitem}

\usepackage{cases}  
\usepackage{tikz-feynman}
\usetikzlibrary{decorations.pathmorphing}         

\usepackage[english]{babel}

\usepackage{caption} 
\usepackage{bm}

\newcommand\numberthis{\addtocounter{equation}{1}\tag{\theequation}}

\newtheorem{thm}{Theorem}[section]
\newtheorem{lem}[thm]{Lemma}
\newtheorem{prop}[thm]{Proposition}

\newtheorem{cor}[thm]{Corollary}

\theoremstyle{definition}
\newtheorem{defn}[thm]{Definition}
\newtheorem{nota}[thm]{Notation}
\newtheorem{remark}[thm]{Remark}
\newtheorem{example}[thm]{Example}

\allowdisplaybreaks

\newcommand\ds\displaystyle
\newcommand\ts\textstyle

\renewcommand{\phi}{\varphi}                 
\renewcommand{\epsilon}{\varepsilon}
\newcommand\eset{\varnothing}

\renewcommand\emptyset\eset
\renewcommand\ell{l}

\newcommand\cupdot {\mbox{\hspace{.15em}$\cup$\hspace{-.47em}$\cdot$\hspace{.4em}}}

\newcommand\supp{\operatorname{supp}}

\newcommand\cB{\mathcal{B}}
\newcommand\cC{\mathcal{C}}
\newcommand\cD{\mathcal{D}}
\newcommand\cF{\mathcal{F}}
\newcommand\cH{\mathcal{H}}
\newcommand\cI{\mathcal{I}}

\newcommand\cM{\mathcal{M}}

\newcommand\cP{\mathcal{P}}
\newcommand\cS{\mathcal{S}}

\newcommand\cU{\mathcal{U}}

\newcommand \tE{{\widetilde E}}

\newcommand \tcC{{\widetilde \cC}}
\newcommand \tcM{{\widetilde \cM}}
\newcommand \tcB{{\widetilde \cB}}
\newcommand \tY{{\widetilde Y}}

\newtheoremstyle{case}{}{}{}{}{}{:}{}{}
\theoremstyle{case}

\newcommand \Mod{\mathrm{Mod}}
\newcommand \rank{\mathrm{rank}}
\newcommand\MS{\operatorname{MinSupp}}
\newcommand \cl{\mathrm{cl}}

\usepackage[top=1.25in, bottom=1.25in, left=1in, right=1in]{geometry}                		         		

\usepackage{enumitem}
\usepackage{graphicx}
\usepackage{float}

\usetikzlibrary{decorations.pathmorphing}         

\newcommand{\com}[1]{\ignorespaces}

\begin{document}

\title{Single-element extensions of matroids over skew tracts}

\author{Ting Su}

\email{tsu2@binghamton.edu}

\address{Department of Mathematics, Universität Hamburg, Germany}
\address{Current affiliation: Changjiang Geophysical Exploration and Testing Co., Ltd., China}


\begin{abstract} Matroids over skew tracts provide an algebraic framework simultaneously generalizing the notions of linear subspaces, matroids, oriented matroids, phased matroids, and some other ``matroids with extra structure". A single-element extension of a matroid $\cM$ over a skew tract $T$ is a matroid $\widetilde{\cM}$ over $T$ obtained from $\cM$ by adding one more element. Crapo characterized single-element extensions of ordinary matroids, and Las Vergnas characterized single-element extensions of oriented matroids, in terms of single-element extensions of their rank 2 contractions. The results of Crapo and Las Vergnas do not generalize to matroids over skew tracts, but we will show a necessary and sufficient condition on skew tracts, called Pathetic Cancellation, such that the result can generalize to weak matroids over skew tracts. 

Stringent skew hyperfields are a special case of skew tracts which behave in many ways like skew fields. We find a characterization of single-element extensions of strong matroids over stringent skew hyperfields.
\end{abstract}

\maketitle

\section{Introduction}

A \emph{matroid} (\cite{CR70},\cite{Wel76}) is a combinatorial object that abstracts the notion of linear independence in a vector configuration over an arbitrary field. An \emph{oriented matroid} (\cite{BLV78}, \cite{LV75}) is a combinatorial object that abstracts the notion of linear independence in a vector configuration over an ordered field.
 The theories of matroids and oriented matroids are major branches of combinatorics with applications in many fields of mathematics, including topology, algebra, graph theory, and geometry. 
 
In \cite{Bak17}, Baker and Bowler introduced {\em matroids over tracts}, which is an algebraic framework simultaneously generalizing the notion of linear subspaces, matroids, oriented matroids, phased matroids, and some other ``matroids with extra structure".
Pendavingh (\cite{Pen18}) partially extended the theory of matroids over tracts to skew hyperfields and defined weak matroids over skew hyperfields. In \cite{Su20}, Su presented a theory of \emph{matroids over skew tracts}, which generalize both the theory of matroids over tracts and the theory of weak matroids over skew hyperfields by Pendavingh.

The {\em single-element extension} is an operation on matroids over skew tracts. For the construction and analysis of matroids over skew tracts it is of interest to describe the set of all single-element extensions of a given matroid over a skew tract. For instance, this is useful for proofs by induction on the size of the ground set: every rank $d$ matroid over a skew tract $T$ can be obtained from some rank $d$ matroid over $T$ on $d$ elements by a sequence of single-element extensions. Crapo found a characterization of single-element extensions of ordinary matroids in terms of cocircuits (\cite{Crapo65}), and Las Vergnas found a characterization of single-element extensions of oriented matroids in terms of signed cocircuits (\cite{LV78}). Most importantly, they showed that single-element extensions can be understood in terms of single-element extensions of rank 2 contractions. The results of Crapo and Las Vergnas do not generalize to matroids over skew tracts, but we show a necessary and sufficient condition on a skew tract $T$ called {\em Pathetic Cancellation} such that the results can generalize to weak matroids over $T$.

\emph{Stringent skew hyperfields} are a special case of skew tracts which behave in many ways like skew fields (\cite{BS20}). In particular, the stringency implies the Pathetic Cancellation property. We find a similar characterization of single-element extensions of strong matroids over stringent skew hyperfields. 

In the introduction, we will first talk about the characterization of single-element extensions of oriented matroids. Then, we will talk about single-element extensions of matroids over skew tracts and stringent skew hyperfields.

\subsection{Classical results for matroids and oriented matroids}
Throughout we will assume that readers have a basic knowledge of matroids and oriented matroids. As a general reference, we refer to \cite{Oxl92, BLVS+99}.

\subsubsection{single-element extensions of matroid and oriented matroid}\label{subsect.ext.MOM}

There is a nice characterization of single-element extension of matroids and of oriented matroids in terms of cocircuits. Crapo found it for matroids in \cite{Crapo65}, and Las Vergnas found it for oriented matroids in \cite{LV78}. 

First, let us see single-element extensions of oriented matroids by an example. The Topological Representation Theorem for oriented matroids (\cite{FL78}) says that there is a one-to-one correspondence between (equivalence classes of) arrangements of signed pseudospheres in $S^d$ and simple rank $d+1$ oriented matroids. A pseudosphere arrangement in $S^d$ is a collection of subsets that topologically ``looks like" a collection of equators, and a signed pseudosphere arrangement in $S^d$ is an arrangement of pseudospheres with a choice of positive side to each pseudosphere. 

Figure~\ref{Ext.OM.pseudo.arr} represents an oriented matroid $\cM$ of rank 3 on the ground set [5]; it has five signed pseudolines ($\cS_e: e\in [5]$) (given by solid curves) with small arrows pointing to the positive side. The signed cocircuits correspond to the 0-dimensional points: $\cC^* = \{\pm Y_k \,|\, 1\leq k \leq 8\}$, where $-Y_k(e) = -(Y_k(e))$ for all $e\in [5]$. We extended $\cM$ by $p$, with pseudocircle $S_p$ given by a dashed curve. Then there are six more signed cocircuits for $\mathcal{\widetilde{M}}$: $ \{\pm Z_1, \pm Z_2, \pm Z_3\}$. We can see that there is a unique function $\sigma$ from the signed cocircuits set $\cC^*$ of $\cM$ to $\{+,-,0\}$ such that $\sigma$ signifies, for every signed cocircuit $Y$, whether it is supposed to lie in $S_p^+$ or $S_p^-$, or on $S_p$, defined by
\begin{equation}\label{sigma.OM}
\sigma (Y) = \begin{cases}
  + & \text{ if signed cocircuit $Y$ lies in $S_p^+$,}\\
  - & \text{ if signed cocircuit $Y$ lies in $S_p^-$.}\\
  0 & \text{ if signed cocircuit $Y$ lies on $S_p$,}\\
\end{cases}
\end{equation}
For instance, $\sigma(Y_1) = 0$, $\sigma(Y_2) = +$ and $\sigma(Y_4) = -$. We can determine the extension by looking at the $\sigma$ value for the signed cocircuits on every pseudocircle. For instance, as $\sigma(Y_3)=+$ and $\sigma(Y_5)=-$, the pseudocircle $S_p$ should cross the pseudocircle $S_2$ between $Y_3$ and $Y_5$. This function $\sigma$ is called the \emph{localization} corresponding to the extension $\widetilde{\cM}$.

\begin{figure}[htb]
\centering
\begin{tikzpicture}
\draw (2,0) arc(0:360:2); 
\draw (1,1.732) arc(0:-60:4); \node at (-1.15,-1.85) {$5$};

\draw[dashed][decorate,decoration={snake,amplitude=10,segment length=100}] (1,-1.732) -- (-1,1.732);\node at (1.2,-1.85) {$p$};

\draw (1.732,1) arc(90:150:4); \node at (-1.85,-1.15) {$4$};
\draw (1.732,-1) arc(30:90:4); \node at (-1.85,1.15) {$1$};
\draw (2,0) arc(-28:-(180-28):2.26); \node at (-2.15, 0) {$2$};
\draw (1.931,0.519) arc(-41.29:71.39-180:3.6); \node at (-2.1,-0.519) {$3$};

\draw[fill] (-0.26,-1.18) circle (.05); \node at (-0.38,-0.95) {$Y_5$};

\draw[fill] (0,0.6) circle (.05); \node at (-0.43,0.60) {$Y_1$};

\draw[fill] (1.5,-0.635) circle (.05); \node at (1.5,-1) {$Y_3$};

\draw[fill] (-1.5,-0.635) circle (.05);  \node at (-1.5,-0.3) {$Y_4$};
 
\draw[fill] (0.92,0.92) circle (.05);  \node at (0.7,1.1) {$Y_6$};

\draw[fill] (0.08,-0.82) circle (.05); \node at (0.35,-0.85) {$Z_1$};

\draw[fill] (0.33,-1.17) circle (.05); \node at  (0.32,-1.43) {$Z_2$};

\draw[fill] (0.03,-0.62) circle (.05); \node at (-0.2,-0.45) {$Z_3$};

\draw[fill] (0.29,-0.54) circle (.05); \node at (0.22,-0.3) {$Y_7$};

\draw[fill] (1.1,-0.18) circle (.05); \node at (1.1,-0.45) {$Y_8$};

\draw[fill] (0.685,0.19) circle (.05); \node at   (1,0.23) {$Y_2$};

\draw[->] (-1.63,0.97) -- (-1.7,0.7); 
\draw[->] (-1.95, -0.1) -- (-1.76,0.15); 
\draw[->] (-1.8,-0.57) -- (-1.7, -0.82); 

\draw[->] (-1.65,-0.93) -- (-1.48, -1.13); 

\draw[->] (-0.8,-1.65) -- (-0.55,-1.82); 

\draw[->] (0.95,-1.53) -- (1.2,-1.35);

\node at (0,-2.8) {$(a)$};

\draw (2+7,0) arc(0:360:2); 
\draw (1+7,1.732) arc(0:-60:4); \node at (-1.15+7,-1.85) {$5$};

\draw[dashed][decorate,decoration={snake,amplitude=10,segment length=100}] (1+7,-1.732) -- (-1+7,1.732);\node at (1.2+7,-1.85) {$p$};

\draw (1.732+7,1) arc(90:150:4); \node at (-1.85+7,-1.15) {$4$};
\draw (1.732+7,-1) arc(30:90:4); \node at (-1.85+7,1.15) {$1$};
\draw (2+7,0) arc(-28:-(180-28):2.26); \node at (-2.15+7, 0) {$2$};
\draw (1.931+7,0.519) arc(-41.29:71.39-180:3.6); \node at (-2.1+7,-0.519) {$3$};

\draw[fill] (-0.26+7,-1.18) circle (.05); \node at (-0.38+7,-0.95) {$Y_5$};

\draw[fill] (0+7,0.6) circle (.05); \node at (-0.43+7,0.60) {$Y_1$};

\draw[fill] (1.5+7,-0.635) circle (.05); \node at (1.5+7,-1) {$Y_3$};

\draw[fill] (-1.5+7,-0.635) circle (.05);  \node at (-1.5+7,-0.3) {$Y_4$};
 
\draw[fill] (0.92+7,0.92) circle (.05);  \node at (0.7+7,1.1) {$Y_6$};

\draw[fill] (0.08+7,-0.82) circle (.05); \node at (0.35+7,-0.85) {$Z_1$};

\draw[fill] (0.33+7,-1.17) circle (.05); \node at  (0.32+7,-1.43) {$Z_2$};

\draw[fill] (0.03+7,-0.62) circle (.05); \node at (-0.2+7,-0.45) {$Z_3$};

\draw[fill] (0.29+7,-0.54) circle (.05); \node at (0.22+7,-0.3) {$Y_7$};

\draw[fill] (1.1+7,-0.18) circle (.05); \node at (1.1+7,-0.45) {$Y_8$};

\draw[fill] (0.685+7,0.19) circle (.05); \node at   (1+7,0.23) {$Y_2$};

\node at (7,-2.8) {$(b)$};
\end{tikzpicture}
\caption{The single-element extension of the oriented matroid $\cM$ by $p$ and the single-element extension of the matroid $\underline{\cM}$ by $p$.}
\label{Ext.OM.pseudo.arr}
\end{figure}

As with, let us see single-element extensions of matroids by an example. In Figure~\ref{Ext.OM.pseudo.arr}$(b)$, we extended the underlying matroid $\underline{\cM}$ by $p$, as shown by a dashed curve. Then there are three more cocircuits for $\widetilde{\underline{\cM}}$: $ \{ \underline{Z_1}, \underline{Z_2}, \underline{Z_3}\}$. We can see that there is a unique function $\sigma$ from the cocircuit set $C^*$ to the set $\{1,0\}$ such that $\sigma$ signifies, for every cocircuit $A$, whether it is supposed to lie on the pseudoline $S_p$ or not, defined by
$$\sigma (A) = \begin{cases}
  1 & \text{ if cocircuit $A$ does not lie on the pseudoline $S_p$,}\\
  0 & \text{ if cocircuit $A$ lies on the pseudoline $S_p$.}
\end{cases}$$
For instance, $\sigma(\underline{Y_1}) = 0$ and $\sigma(\underline{Y_2}) = 1$. 
Similarly to the extension for oriented matroid, we can determine this unique extension by looking at the $\sigma$ value for the cocircuits on every pseudoline. This function $\sigma$ is called the \emph{localization} corresponding to the extension $\widetilde{\underline{\cM}}$.

\subsubsection{Reduction of extensions to rank 2} 
Rank 2 oriented matroids are very easy to understand. For example, every rank 2 oriented matroid is realizable. And rank 2 contractions of an oriented matroid correspond to pseudocircles in topological representation. 

In this subsection, we will describe modularity of cocircuits, which helps to reduce questions to rank 2 contractions.




\begin{defn}\label{mod.latt}(\cite{Sta12}) Let $L$ be a lattice. An element $x\in L$ is called an \textbf{atom} if $x \neq \hat{0}$ and there is no $z \in L$ with $\hat{0} < z< x$. Two atoms $x, y \in L$ form a {\bf modular pair} if $l(L_{x\lor y}) = 2$, i.e., $x\neq y$ and there do not exist $z_1,z_2\in L$ with $\hat{0} < z_1< z_2< x \lor y$.
\end{defn}

Let $E$ be a set and let $C$ be a collection of pairwise incomparable nonempty subsets of $E$. The set $U(C):= \{\bigcup S \,|\, S \subseteq C \}$ forms a lattice when equipped with the partial order coming from inclusion of sets, with join corresponding to union and with the meet of $X$ and $Y$ defined to be the union of all sets in $C$ contained in both $X$ and $Y$. So every $X \in C$ is atomic as an element of $U(C)$. We say that $C_1, C_2 \in C$ form \textbf{a modular pair} in $C$ if they are a modular pair in $U(C)$; that is, the height of their join in the lattice $U(C)$ is 2. 

In a matroid $M$ of rank $d$ on $E$ with cocircuit set $C^*$, $C^*$ is a collection of pairwise incomparable nonempty subsets of $E$.

\begin{defn}\label{mod.MOM}
Two cocircuits $A, B \subseteq E$ form \textbf{a modular pair} in the matroid $M$ if $A\neq B$ and $A\cup B$ does not properly contain a union of two distinct elements of $C^*$. 

In an oriented matroid $\cM$ on $E$, two signed cocircuits $X, Y\in \cC^*(\cM)$ form \textbf{a modular pair} in $\cM$ if $\underline{X}$ and $\underline{Y}$ form a modular pair in the underlying matroid $\underline{\cM}$.
\end{defn}

In a pseudosphere arrangement in $S^d$, it is easy to find a modular pair. Two signed cocircuits form a modular pair if and only if the corresponding two points in $S^d$ are on a common pseudocircle. For example, in Figure~\ref{Ext.OM.pseudo.arr}, the signed cocircuits $Y_1$ and $Y_2$ form a modular pair, and the signed cocircuits $Y_1$ and $Y_5$ do not form a modular pair. 

In 1978, Las Vergnas showed that oriented matroid axiomatics can be considerably sharpened (``reduced to a rank 2 situation") if the underlying structure is known to be a matroid: in this case, cocircuit elimination for modular pairs is sufficient to define oriented matroids (\cite{LV78}).

In \cite{Del11}, Delucchi strengthened Las Vergnas' result by showing that the elimination for matroids can be reduced to the Modular Elimination. And this leads to a corresponding strengthening of the cocircuit axioms for oriented matroids.

\begin{thm}(\cite{Del11})
In the standard definition of matroids, the Elimination axiom can be replaced by

(Modular Elimination) Let $C_1, C_2\in C^*$ such that $C_1$ and $C_2$ form a modular pair in $C^*$, and let $e \in C_1\cap C_2$. There is a member $C_3$ of $C^*$ such that $C_3\subseteq (C_1\cup C_2) -e$.
\end{thm}

\begin{thm}(\cite{Del11})
In the standard definition of oriented matroids, the Elimination axiom can be replaced by

(Modular Elimination) Let $X, Y\in \cC^*$ such that $\underline{X}$ and $\underline{Y}$ a modular pair of signed cocircuits and $X\neq -Y$, and let $e\in X^+\cap Y^-$. There exists $Z\in \cC^*$ such that 
$$Z^+\subseteq (X^+ \cup Y^+)\backslash \{e\} \text{ and } Z^-\subseteq (X^- \cup Y^-)\backslash \{e\}.$$
\end{thm}

Now let's turn back to consider single-element extensions of oriented matroids. We will see that an extension can be determined by the function $\sigma$ we defined by Equation~\eqref{sigma.OM} in the last section. Recall that the function $\sigma: \cC^* \rightarrow \{+,-,0\}$ signifies, for every signed cocircuit $Y$, whether it is supposed to lie on the positive side or on the negative side of the extended pseudocircle $S_p$ or on $S_p$ itself. The extended oriented matroid $\widetilde{\cM}$ is uniquely determined by $\sigma$, and we can determine all the signed cocircuits of $\widetilde{\cM}$ by looking at the values of $\sigma$ for the signed cocircuits of $\cM$ on every pseudocircle (rank 2 situation).

Moreover, Las Vergnas (\cite{LV78}) proved that if a function $\sigma$ determines an extension on every pseudocircle, then it determines an extension of the oriented matroid.

\begin{thm}\label{OM.ext} (\cite{BLVS+99}, see also \cite{LV78})
Let $\cM$ be an oriented matroid with signed cocircuit set $\mathcal{C^*(M)}$ and let
$$\sigma: \mathcal{C^*(M)} \rightarrow \{+, -, 0\}$$ 
be a function, satisfying $\sigma(- Y) = -\sigma(Y)$ for all $Y \in \mathcal{C^*(M)}$. 

Then the following statements are equivalent.
\begin{enumerate}[(1)]
\item $\sigma$ defines a single-element extension $\mathcal{\widetilde{M}}$ of $\cM$.

\item $\sigma$ defines a single-element extension of every rank 2 contraction of $\cM$.

\item $\sigma$ defines a single-element extension of every rank 2 minor of $\cM$ on three elements.
\end{enumerate}
\end{thm}

For the single-element extension of a matroid, there is not a nice topological picture like that for oriented matroids, but there is a similar result. 
The function $\sigma$ from the cocircuit set to $\{1,0\}$ also uniquely determines the extended matroid. And if $\sigma$ determines an extension of every rank 2 contraction, it also determines an extension of the matroid. 

\begin{thm}\label{M.ext} (\cite{Crapo65})
Let $M$ be a matroid with cocircuit set $C^*(M)$. Then, for a function $\sigma: C^*(M) \rightarrow \{0,1\}$, which to every $Y \in C^*(M)$ assigns a signature $\sigma(Y) \in \{0,1\}$, the following statements are equivalent:
\begin{enumerate}[(1)]
\item$\sigma$ defines a single-element extension $\widetilde{M}$ of $M$.

\item $\sigma$ defines a single-element extension of every rank 2 contraction of $M$.
\end{enumerate}
\end{thm}

This paper will generalize the results of Crapo and Las Vergnas to matroids over skew tracts.

\subsection{New results for matroids over skew tracts and stringent skew hyperfields}\label{sect.intro.ext}
\subsubsection{Skew hyperfields and tracts} A skew hyperfield is an algebraic structure similar to a skew field except that its addition is multivalued (cf. \cite{Viro10}). 
A \emph{hyperfield} is then a skew hyperfield with commutative multiplication. 
In particular, every skew field is a skew hyperfield. 
There are many interesting examples of skew hyperfields including the following. 

The \textbf{Krasner hyperfield} $\mathbb{K}:=\{0, 1\}$ has the usual multiplication rule and hyperaddition is defined by $0\boxplus x = \{x\}$ for $x \in \mathbb{K}$ and $1\boxplus 1 = \{0,1\}$. 

The \textbf{sign hyperfield} $\mathbb{S}:= \{0, 1, -1\}$ has the usual multiplication rule and hyperaddition is defined by $0\boxplus x = \{x\}, x \boxplus x = \{x\}$ for $x \in \mathbb{S}$ and $1\boxplus -1 = \{0,1, -1\}$. 

The \textbf{phase hyperfield} $\mathbb{P}:= S^{1} \cup \{0\}$ has the usual multiplication rule and hyperaddition is defined by $0\boxplus x = \{x\}$ for $x \in \mathbb{P}$, $x \boxplus -x = \{0, x, -x\}$ and $ x \boxplus y = \{ \frac{a x + b y}{|a x + b y|} \,|\,  a, b \in \mathbb{R}_{>0} \}$ for $x, y \in S^{1}$ and $x\neq -y$, as in Figure~\ref{phase.sum}. 

\begin{figure}[htb]
\centering
\begin{tikzpicture}
\draw (2,0) arc(0:360:2); 

\node at (2.2,-0.3) {$1$};

\node at (-2.3,-0.3) {$-1$};

\node at (0.2,2.3) {$i$};

\node at (0.3,-2.3) {$-i$};

\node at (0.2,0.2) {$0$};

\draw[->] (-3,0) -- (3,0);
\draw[->] (0, -3) -- (0,3);

\draw (1.932,0.518) circle (.08); \node at (2.2,0.6) {$y$};

\draw (-1,1.732) circle (.08); \node at (-1.2,1.9) {$x$};

\draw[ultra thick] (1.91,0.594) arc(15+2.3:120-2.3:2);

\node at (1.5, 1.9){$x\boxplus y$};
\end{tikzpicture}
\caption{}
\label{phase.sum}
\end{figure}

A skew hyperfield $F$ is said to be \textbf{stringent} if, for any $a$, $b$ in $F$, $a \boxplus b$ is a singleton whenever $a\neq -b$ (\cite{BS20}). Skew fields, $\mathbb{K}$ and $\mathbb{S}$ are all stringent. 

\begin{example}\label{dd.exam} The phase hyperfield $\mathbb{P}$ does not satisfy stringency as $1\boxplus i = \{e^{\alpha i} \,|\, \alpha\in (0, \pi/2) \}$.
\end{example}

Hyperfields are a special case of still more general objects called tracts, which appear to be a natural setting for matroid theory (\cite{Bak17}). Tracts generalize fields, hyperfields in the sense of Krasner (\cite{Kra57}), partial fields in the sense of Semple and Whittle (\cite{SW96}), and fuzzy rings in the sense of Dress (\cite{Dre86}), as shown in Section 2 in \cite{Bak17}.

A \textbf{tract} is an abelian group $G$ (written multiplicatively), together with an \textbf{additive relation structure} on $G$, which is a subset $N_G$ of the group semiring $\mathbb{N}[G]$ satisfying:
\begin{enumerate}[(1)]
\item The zero element of $\mathbb{N}[G]$ belongs to $N_G$.
\item The identity element $1$ of $G$ is not in $N_G$.
\item There is a unique element $\epsilon$ of $G$ with $1 + \epsilon \in N_G$.
\item $N_G$ is closed under the natural action of $G$ on $\mathbb{N}[G]$.
\end{enumerate}

\subsubsection{Matroids over skew tracts}

Baker and Bowler (\cite{Bak17}) introduced matroids over a tract $T$ and called them $T$-matroids. A $T$-matroid for a field $T$ corresponds to a linear subspace of some $T^n$. A $\mathbb{K}$-matroid is a usual matroid. An $\mathbb{S}$-matroid is an oriented matroid. A $\mathbb{P}$-matroid is a phased matroid (\cite{AD12}).

Baker and Bowler also provided two natural notions of matroids over a tract $T$, {\em weak $T$-matroids} and {\em strong $T$-matroids}, which diverge for certain tracts. This is a new phenomenon because weak matroids and strong matroids are the same for fields, $\mathbb{K}$ and $\mathbb{S}$. They provided both circuit axioms and Grassman-Pl\"ucker axioms for weak and strong $T$-matroids. The weak circuit axioms looks more like the (signed) circuit axioms for (oriented) matroids and generalize them better than the strong axioms. The strong Grassmann-Pl\"ucker axioms looks more like the chirotope axioms for oriented matroids and generalize them better than the weak axioms. Weak matroids and strong matroids coincide in rank 2.

In \cite{Pen18}, Pendavingh partially extended the theory of matroids over tracts to skew hyperfields. He provided several cryptomorphic axiom systems for weak matroids over skew hyperfields, including circuits axioms and a new axiom system in terms of quasi-Pl\"ucker coordinates.

Su (\cite{Su20}) presented matroids over skew tracts, which generalize both matroids over tracts and weak matroids over skew hyperfields by Pendavingh (\cite{Pen18}). Similarly, Su provided two natural notions of matroids over a skew tract $T$, {\em weak $T$-matroids} and {\em strong $T$-matroids}, corresponding to the two notions of matroids over tracts. Su also presented several cryptomorphic axiom systems for both kinds of $T$-matroids.

The main result of this paper is a characterization of the single-element extensions of a weak matroid over a skew tract $T$ if $T$ satisfies a special property, called {\em Pathetic Cancellation}. This generalizes the characterization of the single-element extensions of matroids (Theorem~\ref{M.ext}) and the characterization of the single-element extensions of oriented matroids (Theorem~\ref{OM.ext}). 

A skew tract $T$ satisfies the \textbf{Pathetic Cancellation Property} if for every $a, b, x, y, z \in T^\times$ with 
$$1 + a- x, -1 + b- y, a + b - z, x + y - z \in N_G,$$
$$ ax = xa, by  = yb, a^{-1}zb^{-1}  = b^{-1}za^{-1}, x^{-1}zy^{-1}  = y^{-1}zx^{-1}, $$
we have $$xb -ay -z \in N_G.$$

\begin{thm}\label{maintheorem}
Let $T$ be a skew tract, let $\cM$ be a weak left (resp. right) $T$-matroid on a set $E$ with $T$-cocircuit set $\mathcal{C^*(M)}$, and let
$$\sigma: \mathcal{C^*(M)} \rightarrow T$$ 
be a right (resp. left) $T^\times$-equivariant function. 

Then $T$ satisfies Pathetic Cancellation if and only if the following statements are equivalent.
\begin{enumerate}[(1)]
\item $\sigma$ defines a weak single-element extension of $\cM$.

\item $\sigma$ defines a weak single-element extension of every rank 2 contraction of $\cM$.

\item $\sigma$ defines a weak single-element extension of every rank 2 minor of $\cM$ on three elements.
\end{enumerate}
\end{thm}

Bowler and Pendavingh showed that the notions of weak and strong $F$-matroid coincide if $F$ is a stringent skew hyperfield in \cite{BP19}. 

We will show that the stringency implies the Pathetic Cancellation property. So a similar characterization of single-element extensions is also true for strong matroids over stringent skew hypefields. 

\begin{thm}
Any stringent skew hyperfield satisfies Pathetic Cancellation.
\end{thm}

There is also an interesting fact showed in Su's thesis (\cite{Su18}) that a hyperfield not satisfying the Pathetic Cancellation property might fail even weaker generalization of Crapo and Las Vergnas' characterization of extensions, as follows.

\begin{remark} \label{strong} (\cite{Su18})
There is a strong $T$-matroid $\cM$ of rank 4 on 6 elements and a $T^\times$-equivariant function $\sigma: \mathcal{C^*(M)} \rightarrow T$ such that $\sigma$ defines a strong single-element extension by a new element $p$ on every contraction of $\cM$ of rank $\leq 3$ but $\sigma$ is not a strong localization.
\end{remark}

\subsection{Structure of the paper}
In Section~\ref{sect.background}, we will give background, including definition of skew tracts, two natural notions of matroids over skew tracts and two cryptomorphic axioms systems for both two kinds of matroids over skew tracts: circuit axioms and quasi-Pl\"ucker axioms. In Section~\ref{sect.ext.local}, we will define extensions and localizations, and present the quasi-Pl\"ucker coordinates and $T$-cocircuit set for the extended $T$-matroid.
In Section~\ref{sect.PCP} and \ref{sect.chara}, we will show the main result of this paper, a characterization of single-element extensions of a weak matroid over a skew tract $T$ if $T$ satisfies Pathetic Cancellation and present the proof. In Section~\ref{sect.DD}, we will show that the stringency implies the Pathetic Cancellation property and present a similar characterization of single-element extensions for strong matroids over stringent skew hypefields.

\subsection{Acknowledgment} Many thanks to my PhD advisor Laura Anderson for all help and guidance on my thesis which is the original version (commutative tracts case) of this paper and for helpful feedback and suggestions. Thanks also to Nathan Bowler for asking whether the results might hold for skew tracts and for helpful discussion and suggestions.

\section{Background}\label{sect.background}
\subsection{Skew hyperfields and skew tracts}
\begin{defn}\label{hyperoperation} A \textbf{hyperoperation} on a set $S$ is a map $\boxplus$ from $S \times S$ to the collection of non-empty subsets of $S$.

If $A$, $B$ are non-empty subsets of $S$, we define
$$A \boxplus B := \bigcup_{a\in A, b\in B} (a\boxplus b) $$
and we say that $\boxplus$ is \textbf{associative} if $a\boxplus(b\boxplus c)=(a\boxplus b) \boxplus c$ for all $a, b, c \in S$.
\end{defn}

All hyperoperations in this paper will be commutative and associative.

\begin{defn}\label{hypergroup} (\cite{Bak16, Pen18}) 
A (commutative) \textbf{hypergroup} is a triple $(G,\boxplus, 0)$ where $\boxplus$ is a commutative and associative hyperoperation on $G$ such that:
\begin{enumerate}[(1)]
\item $0\boxplus x = \{x\}$ for all $x \in G$. 

\item For every $x \in G$ there is a unique element of G (denoted by $-x$ and called the {\bf hyperinverse} of $x$) such that $0 \in x \boxplus -x$.

\item $x \in y \boxplus z$ if and only if $z \in x \boxplus -y$.
\end{enumerate}

A \textbf{skew hyperfield} is a tuple $(F,\odot, \boxplus, 1, 0)$ such that $0 \neq 1$:
\begin{enumerate}[(1)]
\item $(F-\{0\}, \odot, 1)$ is a group with respect to the multiplication $\odot$ of the skew hyperfield $F$.

\item $(F, \boxplus, 0)$ is a hypergroup.

\item (Absorption rule) $x\odot 0 = 0\odot x = 0$ for all $x \in F$.

\item (Distributive Law) $a \odot (x \boxplus y )= (a \odot x) \boxplus (a \odot y)$ and $(x \boxplus y ) \odot a = (x \odot a) \boxplus (y \odot a)$ for all $a, x, y \in F$.
\end{enumerate}

A \textbf{hyperfield} is then a skew hyperfield with commutative multiplication.
\end{defn}

If $F$ is a skew hyperfield, then $F^E$ has a hypergroup structure given by 
$$X \boxplus Y := \{Z \,|\, Z(e) \in X(e) \boxplus Y (e), \forall e\in E \}.$$
$F$ acts on $F^E$ by componentwise multiplication. (Formally, $F^E$ is an $F$-module (cf. Subsection 2.18 in \cite{Bak16}).)

\begin{example}\label{example.hyp} Now we would like to present three more examples of skew hyperfields. 

\begin{enumerate}[(1)]
\item If $F$ is a (skew) field, then $F$ is a (skew) hyperfield with $a\odot b = a \cdot b$ and $a\boxplus b = \{a + b\}$, for $a, b \in F$.



\item Let $D_6$ be the dihedral group of order 6. Let $F: = D_6 \cup \{0\}$ with multiplication given by $x\cdot y = 0$ if $x$ or $y$ is 0 and by the multiplication of $D_6$ otherwise. Hyperaddition is given by
$$ x\boxplus y = \begin{cases}
F\backslash \{0\} & \text{ if } x \neq y. \\
F & \text{ if } x = y.
\end{cases}
$$
for all $x,y \in D_6$ and $0\boxplus x = x\boxplus 0 = \{x\}$ for all $x\in F$.

It is easy to check that $F$ is a skew hyperfield.

\item Let $F: = \mathbb{H}/\mathbb{R}_{>0} = \{[g]=g\cdot \mathbb{R}_{>0} | g\in \mathbb{H}\}$ be the quotient of the skew field of quarternions $\mathbb{H}$ by the positive real numbers. The multiplication of $F$ is given by $[g] \cdot [h] = [gh]$, for $[g], [h] \in \mathbb{H}/\mathbb{R}_{>0} $ and hyperaddition is given by $[g]\boxplus [0] = [g]$ and $[g] \boxplus [h] = \{[f] \subseteq \mathbb{H}/\mathbb{R}_{>0}  \, | \, f \in g\cdot \mathbb{R}_{>0} + h\cdot \mathbb{R}_{>0}\}$, for $[g], [h] \in (\mathbb{H}/\mathbb{R}_{>0} )^{\times}$. 

The quotient hyperfield was introduced by Krasner in \cite{Kra83}. $F$ is a skew hyperfield following from the argument with $\mathbb{R}_{>0}$ a normal subgroup of $\mathbb{H}^\times$.

\end{enumerate}
\end{example}

\begin{defn}\label{tract}(\cite{Su20, Bak17}) A \textbf{skew tract} is a group $G$ (written multiplicatively), together with an \textbf{additive relation structure} on $G$, which is a subset $N_G$ of the group semiring $\mathbb{N}[G]$ satisfying:
\begin{enumerate}[(1)]
\item The zero element of $\mathbb{N}[G]$ belongs to $N_G$.
\item The identity element $1$ of $G$ is not in $N_G$.
\item There is a unique element $\epsilon$ of $G$ with $1 + \epsilon \in N_G$.
\item $N_G$ is closed under the natural left and right actions of $G$ on $\mathbb{N}[G]$.
\end{enumerate}

A \textbf{tract} is a skew tract with the group $G$ abelian.
\end{defn}

One thinks of $N_G$ as those linear combinations of elements of G which ``sum to zero" (the $N$ in $\mathbb{N}[G]$ stands for ``null").

We let $T = G\cup \{0\}$ and $T^{\times} = G$. We often refer to the skew tract $(G, N_G)$ simply as $T$. Because of the following Lemma~\ref{epsilon}, we often write $-1$ instead of $\epsilon$ and $-x$ instead of $\epsilon x$ or $x\epsilon$ for $x \in G$.

\begin{lem}\label{epsilon}(\cite{Su20, Bak17}) Let $T$ be a skew tract.
\begin{enumerate}
\item If $x,y\in G$ satisfy $x+y\in N_G$, then $y = \epsilon x = x\epsilon$.
\item $\epsilon^2 = 1$.
\item $G\cap N_G = \emptyset.$
\end{enumerate}
\end{lem}

For a skew hyperfield $F$ with hyperaddition $\boxplus$, we could define a skew tract $(G, N_G)$ associated to this skew hyperfield by setting $G = F\backslash \{0\}$ and $N_G = \{\sum_{i=1}^k g_i \,|\, \forall i, g_i \in G \text{ and } 0\in \boxplus_{i=1}^k g_i\}$.

If $T$ is a skew tract and $X, Y\in T^E$, then we define,
$$X+Y: = \{Z\,|\, \forall e\in E, X(e)+Y(e) - Z(e) \in N_G\},$$
$T$ acts on $T^E$ by componentwise multiplication.

\begin{defn}\label{tract.homo}(\cite{Su20}) A \textbf{homomorphism} $f: (G, N_G) \rightarrow (H, N_H)$ of skew tracts is a group homomorphism $f : G \rightarrow H$, together with a map $f : \mathbb{N}[G] \rightarrow \mathbb{N}[H]$ satisfying $f(\sum_{i=1}^k g_i) = \sum_{i=1}^k f(g_i)$ for $g_i \in G$, such that if $\sum_{i=1}^k g_i \in N_G$ then $\sum_{i=1}^k f(g_i) \in N_H$.
\end{defn}

\begin{defn}\label{involut}(\cite{Su20})
Let $T$ be a skew tract. An \textbf{involution} of $T$ is a homomorphism $\tau: T \rightarrow T$ such that $\tau^2$ is the identity map. Denote the image of $x \in T$ under $\tau$ by $\overline{x}$.
\end{defn}

\begin{defn}\label{orthg}(\cite{Su20})
Let $T$ be a skew tract endowed with an involution $x \mapsto \overline{x}$. For $X, Y \in T^E$, the \textbf{product} of $X$ and $Y$ is defined with respect to the involution as $$X\cdot Y := \sum_{e\in E} X(e)\cdot \overline{Y(e)}.$$
Note that $X\cdot Y\in \mathbb{N}[G]$. We say that $X, Y$ are \textbf{orthogonal}, denoted by $X \perp Y$, if $X \cdot Y \in N_G$. 

Let $\cC, \cD \subseteq T^E$. We say that $\cC, \cD$ are \textbf{orthogonal}, denoted by $\cC \perp \cD$, if $X \cdot Y \in N_G$ for all $X\in \cC$ and $Y\in \cD$. We say that $\cC, \cD$ are \textbf{$k$-orthogonal}, denoted by $\cC \perp_k \cD$, if $X \perp Y$ for all $X \in \cC$ and $Y \in \cD$ such that $|
\{e\,|\,X(e)\cdot Y(e)\neq 0\}|\leq k$.
\end{defn}

When $T$ is the field $\mathbb{C}$ of complex numbers or the phase hyperfield $\mathbb{P}$, the usual involution on $T$ is complex conjugation. When $T$ is $\mathbb{K}$ or $\mathbb{S}$, the usual involution on $T$ is the identity map.

\subsection{Matroids over skew tracts} 

\begin{nota} Throughout $E$ denotes a non-empty finite set and $T$ denotes a skew tract. For a skew tract $T$, $T^\times$ denotes $T-\{0\}$.

For simplicity, $E \backslash a$ and $E \cup a$ denote $E \backslash \{a\}$ and $E \cup \{a\}$ respectively.

The $\bf{support}$ of $X \in T^E$ is $\underline{X} := \{e \in E \,|\, X(e) \neq 0\}.$ The $\textbf{zero set}$ of $X$ is $X^0 := \{e \in E \,|\, X(e) = 0\}.$

For $S \subseteq T^E$, $\supp(S)$ or $\underline{S}$ denotes the set of supports of elements of $S$, and $\MS(S)$ denotes the set of elements of $S$ of minimal support.
\end{nota}

We will always view a skew tract $T$ as being equipped with an involution $x \mapsto \overline{x}$. 

First, we will talk about modular pairs and modular families in $T$-matroids.

\begin{defn} (\cite{Bak17}) Let $E$ be a set and let $C$ be a collection of pairwise incomparable nonempty subsets of $E$. We say that $C_1, C_2 \in C$ form a \textbf{modular pair} in $C$ if the height of their join in the lattice $U(C)$ is 2. 
We say that $C_1, ... , C_k \in C$ form a \textbf{modular family} in $C$ if the height of their join in the lattice $U(C)$ is the same as the size of the family $k$. 
\end{defn}

\begin{defn}\label{mod} (\cite{Su20})
Let $\mathcal{C}$ be a subset of $T^E$ such that $\underline{\mathcal{C}}$ is a collection of pairwise incomparable nonempty subsets of $E$. We say that $X, Y \in \cC$ form a \textbf{modular pair} in $\cC$ if $\underline{X}$, $\underline{Y}$ form a modular pair in the lattice of unions of supports of elements of $\cC$. We say that $X_1, ... , X_k \in \cC$ form a \textbf{modular family} in $\cC$ if $\underline{X_1},... ,\underline{X_k}$ form a modular family in the lattice of unions of supports of elements of $\cC$.
\end{defn}

The following is a lemma for modular pairs in terms of rank function of a matroid.

\begin{lem}\label{modpair.rank} Let $M$ be a matroid with rank function $r$. If $Y_1$, $Y_2$ form a modular pair of cocircuits of $M$, then $r(E-(Y_1\cup Y_2)) = \rank(M)-2$.
\end{lem}

\begin{proof} Let $d = \rank(M)$. We know that for a cocircuit $Y$, $E-Y$ is a hyperplane of $M$. Thus by definition of hyperplane, $r(E-Y) = d-1$.

As $E-(Y_1\cup Y_2)$ is a flat of $M$, then
$$r(E-(Y_1\cup Y_2)) < r(E-Y_1) = d-1.$$

By way of contradiction, we assume $r(E-(Y_1\cup Y_2)) < d-2$. Then there exists a flat $E-Z$ with $Z\in U(C^*(M))$ such that
\begin{align*}
E-Y_1 & \supsetneq E-Z \supsetneq E-(Y_1\cup Y_2),\\
Y_1 & \subsetneq Z \subsetneq Y_1\cup Y_2.
\end{align*}
Thus $(Y_1, Y_2)$ is not a modular pair, contradicting the assumption.

So $r(E-(Y_1\cup Y_2)) = d-2$.
\end{proof}

\begin{defn}\label{def.weakcir} (\cite{Su20}) Let $E$ be a non-empty finite set and let $T = (G, N_G)$ be a skew tract. A subset $\mathcal{C}$ of $T^E$ is called the \textbf{$T$-circuit set of a weak left $T$-matroid $\cM$ on $E$} if $\mathcal{C}$ satisfies the following axioms:
\begin{enumerate}[(C1)]
\item \label{c1} $0 \notin \mathcal{C}.$

\item \label{c2} (Symmetry) If $X \in \mathcal{C}$ and $\alpha \in T^{\times}$, then $\alpha \cdot X \in \mathcal{C}$.

\item \label{c3} (Incomparability) If $X$, $Y \in \mathcal{C}$ and $\underline{X} \subseteq \underline{Y}$, then there exists $\alpha \in T^{\times}$ such that $Y= \alpha \cdot X$.

\item \label{c4} (Modular Elimination) If $X$, $Y \in \mathcal{C}$ are a modular pair of $T$-circuits and $e\in E$ is such that $X(e)=-Y(e)\neq 0$, there exists a $T$-circuit $Z \in \mathcal{C}$ such that $Z(e)=0$ and $X(f) + Y(f) - Z(f)\in N_G $ for all $f \in E$.
\end{enumerate}
\end{defn}

In the Modular Elimination axiom, we say $Z$ \textbf{eliminates} $e$ between $X$ and $Y$.

A {\bf weak right $T$-matroid} is defined analogously, with $\alpha \cdot X$ replaced by $X \cdot \alpha$ in (C2) and (C3). 

If $\mathcal{C}$ is the set of $T$-circuits of a weak (left or right) $T$-matroid $\cM$ with ground set $E$, then there is an underlying matroid (in the usual sense) $\underline{\cM}$ on $E$ whose circuits are the supports of the $T$-circuits of $\cM$.

\begin{defn}\label{rank.T} (\cite{Su20}) The \textbf{rank} of $\cM$ is defined to be the rank of the underlying matroid $\underline{\cM}$.
\end{defn}

\begin{defn}\label{def.strongcir} (\cite{Su20}) A subset $\mathcal{C}$ of $T^E$ is called the \textbf{$T$-circuit set of a strong left $T$-matroid $\cM$ on $E$} if $\mathcal{C}$ satisfies \ref{c1}, \ref{c2} and \ref{c3} in Definition~\ref{def.weakcir} and the following stronger version of the Modular Elimination axiom \ref{c4}:
\begin{enumerate}[(C4)$'$]
\item \label{c4'} (Strong Modular Elimination) Suppose $X_1,... ,X_k$ and $X$ are $T$-circuits of $\mathcal{M}$ which together form a modular family of size $k+1$ such that $\underline{X} \not\subseteq \bigcup_{1\leq i \leq k}\underline{X_i}$, and for $1\leq i \leq k$ let
$$e_i \in (\underline{X}\cap \underline{X_i})\backslash \bigcup_{\substack{1\leq j \leq k \\ j\neq i} } \underline{X_j}$$
be such that $X(e_i)=-X_i(e_i)\neq 0$. Then there exists a $T$-circuit $Z \in \mathcal{C}$ such that $Z(e_i)=0$ for $1\leq i \leq k $ and $X(f)+ X_1(f)+ \cdots + X_k(f) - Z(f) \in N_G$ for every $f \in E$.
\end{enumerate}
\end{defn}

A {\bf strong right $T$-matroid} is defined analogously, with $\alpha \cdot X$ replaced by $X \cdot \alpha$ in (C2) and (C3). 

From the definition, it is easy to see that any strong left (resp. right) $T$-matroid on $E$ is also a weak left (resp. right) $T$-matroid on $E$. 

A \textbf{projective $T$-circuit of a (weak or strong) left $T$-matroid} $\cM$ is an equivalence class of $T$-circuits of $\cM$ under the equivalence relation $X_1 \sim X_2$ if and only if $X_1 = \alpha \cdot X_2$ for some $\alpha \in T^\times$. Analogously, a \textbf{projective $T$-circuit of a (weak or strong) right $T$-matroid} $\cM$ is an equivalence class of $T$-circuits of $\cM$ under the equivalence relation $X_1 \sim X_2$ if and only if $X_1 = X_2 \cdot \alpha$ for some $\alpha \in T^\times$.

As \cite{Bak17} shows, a matroid over a field $T$ corresponds to a linear subspace of some $T^n$. A $\mathbb{K}$-matroid corresponds to a usual matroid. An $\mathbb{S}$-matroid is an oriented matroid. A $\mathbb{P}$-matroid is a phased matroid, defined in \cite{AD12}. 
Weak matroids and strong matroids coincide over a field, $\mathbb{K}$ and $\mathbb{S}$, but they do not coincide over $\mathbb{P}$.

There is a cryptomorphic characterization of weak and strong matroids over a skew tract $T$ in terms of {\em quasi-Pl\"ucker coordinates}. First we define some terms.

Let $F$ be a subset of $E$ and let $a_1, ... , a_n$ be distinct elements of $E\backslash F$ with $n\in \mathbb{Z}_{>0}$. For simplicity, we write $Fa_1...a_n := F\cup \{a_1, ... , a_n\}$.

Let $N$ be a matroid of rank $d$ with set of bases $\cB$ and let $T$ be a skew tract. For $B, B'\in \cB$, we say $(B, B')$ is an ordered pair of {\bf adjacent bases} if $|B\backslash B'|=1$. We name the set of ordered pairs of adjacent bases $A_N.$

\begin{defn} \label{def.weakQP} (\cite{Su20}) Let $T$ be a skew tract and let $N$ be a matroid on $E$ with set of bases $\cB$. Then $[\cdot] : A_N \rightarrow T$ are \textbf{weak left quasi-Pl\"ucker coordinates} if
\begin{enumerate}[(P1)]
\item \label{p1} $[Fe, Ff]\cdot [Ff, Fe] = 1$\hfill if $|F|=d-1$ and $Fe, Ff \in \cB$.
\item \label{p2} $[Feg, Ffg]\cdot [Fef, Feg]\cdot [Ffg, Fef] = -1$ \hfill if $|F|=d-2$ and $Fef, Feg, Ffg \in \cB$.
\item \label{p3} $[Fe, Ff]\cdot [Ff, Fg]\cdot [Fg, Fe] = 1$ \hfill if $|F|=d-1$ and $Fe, Ff, Fg \in \cB$.
\item \label{p4} $[Feg, Ffg] = [Feh, Ffh]$ \hfill if $|F|=d-2$ and $Feg, Feh, Ffg, Ffh \in \cB$, and $Fef \notin \cB$ or $Fgh \notin \cB$.
\item \label{p5} $ -1 + [Ffh, Fef]\cdot [Feg, Fgh] + [Feh, Fef]\cdot [Ffg, Fgh] \in N_G$ 

\hfill if $|F|=d-2$ and $Feg$, $Feh$, $Ffg$, $Ffh$, $Fef$, $Fgh$ $\in \cB.$
\end{enumerate}
\end{defn}

\begin{defn}\label{def.strongLQP} (\cite{Su20}) Let $d$ be the rank of $N$. $[\cdot] : A_N \rightarrow T$ are \textbf{strong left quasi-Pl\"ucker coordinates} if $[\cdot]$ satisfies \ref{p1}, \ref{p2} and \ref{p3} in Definition~\ref{def.weakQP} and the following stronger versions of \ref{p4} and \ref{p5} :
\begin{enumerate} [(P4)$'$]
\item \label{p4'} For any two subsets $I, J$ of $E$ with $|I| = d+1$, $|J| = d-1$ and $|I\backslash J|\geq 3$, we let $I_1 = \{x\in I \,|\, \text{both } I\backslash x \text{ and } Jx \text{ are bases of } N\}$. If $|I_1| = 2$ and we say $I =\{a, b\}$, then $$[I\backslash a, I\backslash b] = [Jb, Ja].$$
\end{enumerate}
\begin{enumerate} [(P5)$'$]
\item \label{p5'} For any two subsets $I, J$ of $E$ with $|I| = d+1$, $|J| = d-1$ and $|I\backslash J|\geq 3$, we let $I_1 = \{x\in I \,|\, \text{both } I\backslash x \text{ and } Jx \text{ are bases of } N\}$. If $|I_1|\geq 3$, then for any $z\in I_1$,
$$-1 + \sum_{x\in I_1\backslash z} [I\backslash x, I\backslash z]\cdot [Jx, Jz] \in N_G.$$
\end{enumerate}
\end{defn}

The definition of {\bf weak right quasi-Pl\"ucker coordinates} and {\bf strong right quasi-Pl\"ucker coordinates} are obtained by reversing the order of multiplication throughout.

\begin{thm} (\cite{Su20})
Let $E$ be a finite set, let $N$ be a matroid and let $T$ be a skew tract. There is a natural bijection between the collection of strong (resp. weak) circuit sets $\cC\subseteq T^E$ with $\underline{\cC}$ the circuit set of $N$ satisfying all circuit axioms and the collection of maps $[\cdot]: A_N \rightarrow T$ satisfying all strong (resp. weak) quasi-Pl\"ucker axioms.
\end{thm}

Parallel to matroids and oriented matroids, for every strong (resp. weak) left $T$-matroid $\cM$ of rank $d$ on $E$ with $T$-circuit set $\cC(\cM)$, there is a strong (resp. weak) right $T$-matroid $\cM^*$ of rank $|E|-d$ with $T$-circuit set $\cC^*(\cM) := \MS(\{Y\in T^E\,|\, X\cdot Y \in N_G, X\in \cC(\cM)\}-\{\bf{0}\})$, where $T$ is a skew tract endowed with an involution $x\mapsto \overline{x}$. $\cM^*$ is called the {\bf dual $T$-matroid} of $\cM$.
The $T$-circuits of $\cM^{*}$ are called the \textbf{$T$-cocircuits} of $\cM$, and vice versa.

If $\cM$ is a right $T$-matroid, then the dual matroid $\cM^*$ is a left $T$-matroid and the $T$-cocircuit set $\cC^*(\cM)$ is obtained by reversing the order of multiplication throughout.

\begin{lem} \label{PivotProp} (\cite{Su20}) (Dual pivoting property) Let $T$ be a skew tract and let $\cM$ be a left (weak or strong) $T$-matroid with $T$-cocircuit set $\cC^*$ and left quasi-Pl\"ucker coordinates $[\cdot]$. Let $Y\in \cC^*$. Choose a maximal independent set $J$ in $Y^0$. Then for every $y_1, y_2 \in \underline{Y}$, 
$$Y(y_1)Y(y_2)^{-1} = \overline{[ Jy_1, Jy_2]}.$$
\end{lem}

\begin{lem}\label{exist} Let $\cl$ be the closure operator of $\underline{\cM}$.
Let $Y_1, Y_2 \in \mathcal{C^*(M)}$ such that $Y_1, Y_2$ form a modular pair of $T$-cocircuits. Let $I$ be a maximal independent set in $Y_1^0 \cap Y_2^0$. Then $|I|=d-2$, and for every $e \in \underline{Y_1} \cap \underline{Y_2}$, there exists $Y_e \in \mathcal{C^*(M)}$ with
$\underline{Y_e} = E\backslash \cl(e \cup I)$.
\end{lem}

\begin{proof} As $Y_1$, $Y_2$ form a modular pair in $\mathcal{C^*(M)}$, then $\underline{Y_1}$, $\underline{Y_2}$ form a modular pair in $C^*(\underline{\cM})$. Then by Lemma~\ref{modpair.rank}, $\rank(Y_1^0 \cap Y_2^0) = d-2$. So we can choose a maximal independent set $\{z_1,...,z_{d-2}\}$ in $Y_1^0 \cap Y_2^0$.

Let $e \in \underline{Y_1} \cap \underline{Y_2}$. Then $e \not\in Y_1^0 \cup Y_2^0$. Thus $e \not\in \cl(\{z_1,...,z_{d-2}\})$. So $\rank(\{e, z_1,...,z_{d-2}\}) = (d-2) +1 = d-1$, and thus $\{e, z_1,...,z_{d-2}\}$ is a hyperplane of $\underline{\cM}$. Then there exists $Y_e \in \mathcal{C^*(M)}$ with
$\underline{Y_e} = E\backslash \cl(\{ e, z_1, ... ,z_{d-2}\})$.
\end{proof}

\subsection{Minors}
Let $X\in T^E$ and let $A \subseteq E$. Define $X\backslash A\in T^{E\backslash A}$ by $(X\backslash A)(e)=X(e)$ for $e\in E\backslash A$. 

For $\cU\subseteq T^E$, define the {\bf deletion} of $A$ from $\cU$ as 
$$\cU\backslash A = \{ X\backslash A \,|\, X\in \cU, \underline{X}\cap A = \emptyset\}.$$
Define the {\bf contraction} of $A$ in $\cU$ as 
$$\cU/A = \MS(\{X\backslash A \,|\, X\in \cU\}).$$

\begin{thm}\label{minor} (\cite{Su20}) $\mathcal{C}\backslash A$ is the set of $T$-circuits of a strong (resp. weak) left $T$-matroid $\cM\backslash A$ on $E\backslash A$, called the \textbf{deletion} of $\cM$ by $A$, whose underlying matroid is $\underline{\cM}\backslash A$. 

Similarly, $\cC/A$ is the set of $T$-circuits of a strong (resp. weak) left $T$-matroid $\cM/A$ on $E\backslash A$, called the \textbf{contraction} of $\cM$ by $A$, whose underlying matroid is $\underline{\cM}/A$. 

Moreover, $(\cM\backslash A)^{*} = \cM^{*}/A$ and $(\cM/A)^{*} =\cM^{*}\backslash A $.
\end{thm}

\begin{lem}\label{qp.del.cont} (\cite{Su20})
Let $N$ be a matroid on $E$, let $T$ be a skew tract, let $[\cdot]: A_N \rightarrow T$ be weak left (resp. right) quasi-Pl\"ucker coordinates, and let $A$ be a subset of $E$. 
\begin{enumerate}
\item (Contraction) Let $I_A\subseteq A$ be such that $I_A$ is a basis of $N|A$. Define $[\cdot]/A: A_{N/A} \rightarrow T$ by
$$[B, B']/A := [B\cupdot I_A, B'\cupdot I_A]$$
for all $(B, B') \in A_{N/A}.$ 

$[\cdot]/A$ are weak left (resp. right) quasi-Pl\"ucker coordinates of $N/A$. If $[\cdot]$ are strong left (resp. right) quasi-Pl\"ucker coordinates, so is $[\cdot]/A$. The definition of $[\cdot]/A$ is independent of the choice of $I_A$.

\item (Deletion) Let $J_A\subseteq A$ be such that $J_A$ is a basis of $N/(E\backslash A)$. Define $[\cdot]\backslash A : A_{N\backslash A} \rightarrow T$ by
$$[B, B']\backslash A := [B\cup J_A, B'\cup J_A]$$
for all $(B, B') \in A_{N\backslash A}.$

$[\cdot]\backslash A$ are weak left (resp. right) quasi-Pl\"ucker coordinates of $N\backslash A$. If $[\cdot]$ are strong left (resp. right) quasi-Pl\"ucker coordinates, so is $[\cdot]\backslash A$. The definition of $[\cdot]\backslash A$ is independent of the choice of $J_A$.

\item $([\cdot]\backslash A)^* = [\cdot]^*/A$.
\end{enumerate}
\end{lem}

\begin{defn}\label{equivariant} Let $T$ be a skew tract and $E$ be a non-empty finite set. Let $\cU \subseteq T^E$ such that $\cU$ is closed under both left and right nonzero scalar multiplications. 

A function $f: \cU \rightarrow T$ is called {\bf right $T^\times$-equivariant} if it satisfies
$$f(U\cdot \alpha) = f(U)\cdot \alpha,$$
for all $U\in \cU$ and all $\alpha \in T^\times$.

Similarly, a function $g: \cU \rightarrow T$ is called {\bf left $T^\times$-equivariant} if it satisfies
$$g(\beta \cdot V) = \beta \cdot g(V),$$
for all $V\in \cU$ and all $\beta \in T^\times$.

If $T$ is a (commutative) tract, then a function is called {\bf $T^\times$-equivariant} if it is both left and right $T^\times$-equivariant.
\end{defn}

\begin{prop} 
Let $\cM$ be a left (resp. right) $T$-matroid on $E$ and let $f: \cC^*(\cM) \rightarrow T$ be a right (resp. left) $T^\times$-equivariant function. Let $A\subseteq E$. We can get \textbf{induced} maps $f/A: \cC^*(\cM/A) \rightarrow T$ on $\cC^*(\cM/A)$ and $f\backslash A: \cC^*(\cM\backslash A) \rightarrow T$ on $\cC^*(\cM\backslash A)$
by
$$(f/A)(Z) = f(Y) \text{ and } (f\backslash A)(U) = f(X),$$
where $Z = Y\backslash A$, $Y \in \mathcal{C^*(M)}$, $\underline{Y} \cap A = \emptyset$, and $U = X\backslash A$, $X \in \mathcal{C^*(M)}$. We say that $f/A$ is {\bf induced from $\sigma$ by contraction} of $A$ and $f\backslash A$ is {\bf induced from $\sigma$ by deletion} of $A$.
\end{prop}

\subsection{Rescaling}

\begin{defn}(\cite{Su20}) Let $T$ be a skew tract, let $X\in T^E$ and let $\rho:E \rightarrow T^\times$. Then \textbf{right rescaling} $X$ by $\rho$ yields the vector $X\cdot\rho\in T^\times$ with entries $(X\cdot\rho)(e) = X(e)\cdot \rho(e)$ for all $e\in E$. Similarly, \textbf{left rescaling} gives a vector $\rho\cdot X$. We use $\rho^{-1}$ for the function from $E$ to $T^\times$ such that $\rho^{-1}(e) = \rho(e)^{-1}$ for all $e\in E$.
\end{defn}

\begin{lem}(\cite{Su20})\label{rescaling} Let $\cM$ be a weak (resp. strong) left $T$-matroid on $E$ with $T$-circuit set $\cC$ and $T$-cocircuit set $\cD$ and let $\rho:E \rightarrow T^\times$. Then $\cC \cdot\rho^{-1}$ and $\rho \cdot \cD$ are the $T$-circuit set and $T$-cocircuit set of a weak (resp. strong) left $T$-matroid $\cM^\rho$, where
$$\cC\cdot\rho^{-1}:=\{X\cdot\rho^{-1}\,|\, X\in \cC\} \text{ and } \rho\cdot\cD:=\{\rho\cdot Y\,|\, Y\in \cD\}.$$
$\cM^\rho$ is called the $T$-matroid arising from $\cM$ by \textbf{right rescaling}.
\end{lem}

\subsection{Crapo characterization of extension}\label{sect.crapo}
In this subsection, we will present Crapo's result of extension of matroids. 

\begin{defn}(\cite{Crapo65,Oxl92}) Let $M$ be a matroid. A set $\cF$ of flats of $M$ is called a \textbf{modular cut} if it satisfies the following properties:
\begin{enumerate}[(1)]
\item If $F\in \cF$ and $F'$ is a flat of $M$ containing $F$, then $F'\in \cF$.
\item If $F_1, F_2\in \cF$ and $(F_1, F_2)$ is a modular pair, then $F_1\cap F_2 \in \cF$.
\end{enumerate}
\end{defn}

The following lemma shows that every single-element extension of a matroid gives rise to a modular cut.

\begin{lem}\label{ext.modularcut} (\cite{Crapo65, Oxl92}) Let $N$ be an extension of a matroid $M$ by an element $p$ and let $\cF$ be the set of flats $F$ of $M$ such that $F\cup p$ is a flat of $N$ having the same rank as $F$. Then $\cF$ is a modular cut of $M$.
\end{lem}

Moreover, every modular cut gives rise to a unique extension, as the following theorem shows.

\begin{thm}\label{mdcut.ext}(\cite{Crapo65, Oxl92}) Let $\cF$ be a modular cut of a matroid $M$ on a set $E$. Then there is a unique extension $N$ of $M$ on $E\cupdot p$ such that $\cF$ consists of those flats $F$ of $M$ for which $F\cup p$ is a flat of $N$ having the same rank as $F$. Moreover, for all subsets $X$ of $E$,
$$r_N(X) = r_M(X) \text{ and}$$
$$r_N(X\cup p) = \begin{cases}
r_M(X) & \text{ if } \cl_M(X) \in \cF,\\
r_M(X)+1 & \text{ if } \cl_M(X) \notin \cF.
\end{cases}
$$
\end{thm}

On combining Lemma~\ref{ext.modularcut} and Theorem~\ref{mdcut.ext}, we get that there is a one-to-one correspondence between single-element extensions of a matroid and its modular cuts.

There is a more compact way to specify single-element extensions by hyperplanes. Before showing it, we would like to introduce a new term {\em linear subclass}.

\begin{defn} (\cite{Crapo65, Oxl92}) A \textbf{linear subclass} of a matroid $M$ is defined to be a subset $\cH'$ of the set of hyperplanes of $M$ such that if $H_1$ and $H_2$ are members of $\cH'$ for which $r(H_1\cap H_2) = r(M) -2$, and $H_3$ is a hyperplane containing $H_1\cap H_2$, then $H_3\in \cH'$.
\end{defn}

The following lemma shows us that every modular cut gives rise a a linear subclass, and vice versa.

\begin{lem}(\cite{Crapo65, Oxl92}) Let $M$ be a matroid.
\begin{enumerate}[(1)]
\item If $\cF$ is a modular cut of $M$, then the hyperplanes of $M$ in $\cF$ form a linear subclass.
\item If $\cH'$ is a linear subclass of $M$ and $\cF$ consists of all flats $F$ of $M$ for which every hyperplane containing $F$ is in $\cH'$, then $\cF$ is a modular cut of $M$.
\end{enumerate}
\end{lem}

So it is obvious that every single-element extension $N$ of a matroid $M$ on an element $p$ gives rise to a linear subclass of $M$ consisting of all hyperplanes $H$ of $M$ such that $H\cup p$ is a hyperplane of $N$. Moreover, the other direction also holds.

\begin{thm}(\cite{Crapo65, Oxl92})
If $\cH'$ is a linear subclass of a matroid $M$, then there is a unique extension $N$ of $M$ by an element $p$ such that $\cH'$ is the set of hyperplanes $H$ of $M$ for which $H\cup p$ is a hyperplane of $N$.
\end{thm}

Therefore, there is a one-to-one correspondence between single-element extensions of a matroid and its linear subclasses.

We all know that for a matroid $M$ on a set $E$, a set $Y\subseteq E$ is a cocircuit of $M$ if and only if $E-Y$ is a hyperplane. Thus a subset $\cH'$ of the set of hyperplanes could give rise to a collection $\cD'$ of cocircuits such that $\cD' := \{E-H\,|\, H\in \cH'\}$, and vice versa.

\begin{cor} Let $\cD'$ be a collection of cocircuits of a matroid $M$ and let $\cH'$ be a subset of the set of hyperplanes $H$ of $M$ such that $E-H\in \cD'$. If $\cH'$ is a linear subclass of $M$, then there is a unique extension $N$ of $M$ by an element $p$ such that $\cD'$ is the collection of cocircuits $Y$ of $M$ for which $Y$ is a cocircuit of $N$.
\end{cor}

For a cocircuit $Y$ of $M$, either $Y$ or $Y\cup p$ is a cocircuit of $N$. The following holds obviously.

\begin{cor}\label{matroid.ext.sigma} Let $M$ be a matroid with cocircuit set $C^*(M)$. Then, for a function $\sigma: C^*(M) \rightarrow \{0,1\}$, which to every $Y \in C^*(M)$ assigns a signature $\sigma(Y) \in \{0,1\}$, the following statements are true:
\begin{enumerate}[(1)]
\item Let $N$ be an extension of $M$ by an element $p$. Then $\sigma$ can be uniquely defined by
$$\sigma(Y) = \begin{cases}
0 & \text{ if } Y \text{ is a cocircuit of } N, \\
1 & \text{ if } Y\cup p \text{ is a cocircuit of } N.
\end{cases}
$$
\item Let $\cH'$ be a subset of the set of hyperplanes $H$ of $M$ such that $\sigma(E-H) = 0$. If $\cH'$ is a linear subclass of $M$, then there is a unique extension $N$ of $M$ by an element $p$ such that for $Y\in C^*(M)$, $\sigma(Y) = 0$ for which $Y$ is a cocircuit of $N$.
\end{enumerate}
\end{cor}

We say $\sigma$ {\bf defines a single-element extension} of $M$ if it satisfies (2) in Corollary~\ref{matroid.ext.sigma}. 

Following is also a useful lemma due to Crapo.

\begin{lem} \label{ext.hypl} (\cite{Crapo65, BLVS+99})
The hyperplanes of a single-element extension $N$ of a matroid $M$ by $p$ are given by one of the following two types:
\begin{enumerate}[(1)]
\item either $H$ or $H \cup p$, for every hyperplane $H$ of $M$, and

\item the sets $G \cup p$, where $G = H_1 \cap H_2$ is an intersection of two hyperplanes of $M$ for which there is no hyperplane $H \supset G$ of $M$ such that $H \cup p$ is a hyperplane of $N$.
\end{enumerate}
\end{lem}






\section{Extension and Localization}\label{sect.ext.local}
In this section, we will introduce single-element extensions of matroids over skew tracts and localizations, and show that localizations characterize extensions. 

\begin{defn}\label{extens} Let $T$ be a skew tract and let $\cM$ be a weak (resp. strong) left (resp. right) $T$-matroid on $E$. A weak (resp. strong) \textbf{extension} of $\cM$ is a weak (resp. strong) left (resp. right) $T$-matroid $\mathcal{\widetilde{M}}$ on a ground set $\widetilde{E}$ that contains $E$, such that the deletion $\mathcal{\widetilde{M}}\backslash (\widetilde{E}-E)$ is equal to $\cM$.

$\mathcal{\widetilde{M}}$ is a weak (resp. strong) \textbf{single-element extension of $\cM$} if $|\widetilde{E}\backslash E| = 1$, that is, $\widetilde{E} = E \cup p$ for some $p \notin E$. So, $\mathcal{\widetilde{M}}\backslash p = \cM$.

We will exclude the \emph{trivial} case where $p$ is a coloop of $\mathcal{\widetilde{M}}$ (equivalently, if $\rank(\mathcal{\widetilde{M}}) = \rank(\cM)+1$) as $\mathcal{C}^*(\widetilde{\cM}) = \{Y\in T^{E\cup p} \,|\, Y(p)=0 \text{ and } Y\backslash p \in \mathcal{C}^*(\cM)\} \cup \{Y\in T^{E\cup p} \,|\, Y(p)\neq 0\text{ and } Y(e)=0 \text{ for all } e\in E\}$. Thus all single-element extensions considered in this paper are $\textbf{non-trivial}$.
\end{defn}

Note that any strong $T$-matroid on $E$ is a weak $T$-matroid on $E$, and any weak $T$-matroid of rank $\leq 3$ on $E$ is also a strong $T$-matroid on $E$.

If $\cM$ is a weak (resp. strong) $T$-matroid, then of course extensions of $\cM$ are weak (resp. strong) $T$-matroid and localizations are weak (resp. strong). We will omit ``weak (resp. strong)" if theorems, propositions and lemmas are true for both weak and strong cases.

If $Y\in T^E$, $p\notin E$ and $\alpha \in T$, then we denote by $(Y,\alpha)$ to denote the extension of $Y$ to $E\cup p$ with $Y(p) = \alpha$ when there is no confusion.

\begin{prop}\label{sigma.intro}  
Let $T$ be a skew tract and let $\cM$ be a weak (resp. strong) left $T$-matroid on $E$ with $T$-cocircuit set $\mathcal{C^*(M)} \subseteq T^E$. Let $\mathcal{\widetilde{M}}$ be a weak (resp. strong) single-element extension of $\cM$ by an element $p$ with $T$-cocircuit set $\mathcal{C^*(\widetilde{M})} \subseteq T^{E\cup p}$. Then for every $T$-cocircuit $Y \in \mathcal{C^*(M)}$, there is a unique way to extend $Y$ to a $T$-cocircuit of $\mathcal{\widetilde{M}}$: there is a unique function
$$\sigma: \mathcal{C^*(M)} \rightarrow T$$
such that
$$\{(Y, \sigma(Y))\,|\, Y \in \mathcal{C^*(M)}\} \subseteq \mathcal{C^*(\widetilde{M})}.$$

That is, $(Y, \sigma(Y))$ is a $T$-cocircuit of $\mathcal{\widetilde{M}}$ for every $T$-cocircuit $Y$ of $\cM$. Furthermore, $\sigma$ is right $T^\times$-equivariant.
\end{prop}

The statement also holds when $\cM$ is a (weak or strong) right $T$-matroid and $\sigma$ is left $T^\times$-equivariant.

\begin{defn}\label{localiz}
The functions $\sigma$ that correspond to weak (resp. strong) single-element extensions (via Proposition~\ref{sigma.intro}) are called weak (resp. strong) $\textbf{localizations}$.
\end{defn}

The definition of {\bf localizations for right $T$-matroids} is obtained by reversing the order of multiplication throughout.

\begin{proof}[Proof of Proposition~\ref{sigma.intro}] Let $Y \in \mathcal{C^*(M)}$. By Theorem~\ref{minor}, $Y = \widetilde{Y}\backslash p$ for some $ \widetilde{Y} \in \mathcal{C^*(\widetilde{M})}$. By Incomparability in $\mathcal{C^*(\widetilde{M})}$, $ \widetilde{Y}$ is uniquely defined. Define $\sigma(Y) = \widetilde{Y}(p)$. Thus $(Y, \sigma(Y)) \in \mathcal{C^*(\widetilde{M})}$.

Also for all $\alpha \in T^\times$, $(Y, \sigma(Y))\cdot \alpha \in \mathcal{C^*(\widetilde{M})}$, by Symmetry. That is, $(Y\cdot \alpha, \sigma(Y)\cdot \alpha) \in \mathcal{C^*(\widetilde{M})}$.

So $$\sigma(Y\cdot \alpha) = \sigma(Y)\cdot \alpha.$$
and so $\sigma$ is right $T^\times$-equivariant.
\end{proof}

For the remainder of this section, we take $\cM$ to be a left $T$-matroid of rank $d$ on $E$ with $T$-cocircuit set $\cC^*(\cM)$ and left quasi-Pl\"ucker coordinates $[\cdot]$, $\widetilde{\cM}$ to be an extension of $\cM$ by $p$, and $\sigma: \mathcal{C^*(M)} \rightarrow T$ to be the corresponding localization of $\cM$. Let $\cl$ be the closure operator of $\underline{\cM}$.

\begin{lem} \label{ext.basis} Let $\cB$ be the set of bases of $\underline{\cM}$. The set of bases of $\underline{\tcM}$ is
$$\tcB = \cB \cup \{Fp \subseteq E\cup p \,|\, |F| = d-1, E \backslash \cl(F) = \underline{Y} \text{ and } \sigma(Y) \neq 0 \text{ for some } Y \in \cC^*(\cM) \}.$$
\end{lem}

\begin{proof} As $\underline{\cM} = \underline{\mathcal{\widetilde{M}}\backslash p}=\underline{\tcM}\backslash p$ and $\rank(\underline{\cM})=\rank(\underline{\tcM})$, we know that the independent sets $\cI(\underline{\cM}) = \cI(\underline{\tcM})\cap \cP(E)$. Thus $\cB\subseteq \tcB$ and for $B\in \tcB\backslash \cB$, $B = Fp$ for some $F\in \cI(\underline{\cM})$. So we only need to prove that $$\tcB\backslash \cB = \{Fp \subseteq E\cup p \,|\, |F| = d-1, E \backslash \cl(F) = \underline{Y} \text{ and } \sigma(Y) \neq 0 \text{ for some } Y \in \cC^*(\cM) \}.$$

For $\subseteq$, let $B\in \tcB\backslash \cB$ with $B = Fp$ and $F\in \cI(\underline{\cM})$. As $|F| = d-1$, then there exists $Y \in \cC^*(\cM)$ such that $E \backslash \cl(F) = \underline{Y}$. We claim that $\sigma(Y)\neq 0$. Otherwise, $(Y,0) \in \cC^*(\tcM)$. Then $p\notin \tE\backslash \cl_{\underline{\tcM}}(F)$ and this contradicts that $Fp$ is a basis of $\underline{\tcM}$.

For $\supseteq$, let $F$ be an independent set in $\underline{\cM}$ with $|F| = d-1$. Then there exists $Y \in \cC^*(\cM)$ such that $E \backslash \cl(F) = \underline{Y}$. If $\sigma(Y) \neq 0$, then $p\in \tE\backslash \cl_{\underline{\tcM}}(F)$. So $Fp\in \tcB$. If $\sigma(Y) = 0$, then $p\notin \tE\backslash \cl_{\underline{\tcM}}(F)$. So $Fp\notin \tcB$. 
\end{proof}

\begin{lem}\label{pivot}
Let $[\cdot]_{\tcM}$ be left quasi-Pl\"ucker coordinates for $\tcM$. Let $F \subseteq E$ be independent in $\underline{\cM}$ with $|F| = d-1$ and let $Y \in \mathcal{C^*(M)}$ with $\underline{Y} = E\backslash \cl(F)$. Let $y \in \underline{Y}$. If $\sigma(Y) \neq 0$, then
$$\sigma(Y)\cdot Y(y)^{-1} = \overline{[Fp, Fy]_\tcM}.$$ 
\end{lem}
\begin{proof}This is proved by Lemma~\ref{PivotProp}.
\end{proof}

\begin{prop}\label{ext.uniq} 
The left $T$-coordinates $[\cdot]_p: A_{\underline{\tcM}} \rightarrow T$ defined by
\begin{equation}
[Fs, Ft]_p=\begin{cases}
[Fs, Ft]	& \text{ if } p\notin Fst,\\
\overline{\sigma(Y)\cdot Y(t)^{-1}}  & \text{ if } p = s\text{ and } Y\in \mathcal{C^*(M)} \text{ with } \underline{Y}= E\backslash\cl(F),\\
\overline{ Y(s)\cdot \sigma(Y)^{-1}}  & \text{ if } p = t \text{ and } Y\in \mathcal{C^*(M)} \text{ with } \underline{Y}= E\backslash\cl(F),\\
[Gsg, Gtg] & \text{ if } p\in F, G=F\backslash \{p\}, Gst\notin \cB \text{ and } Gsg, Gtg\in \cB \\
&\,\,\, \text{ for some } g\in E, \\
-\overline{Y_t(s)\cdot \sigma(Y_t)^{-1}\cdot\sigma(Y_s)\cdot Y_s(t)^{-1}} & \text{ if } p\in F, G=F\backslash \{p\}, Gst\in \cB, Y_i\in \mathcal{C^*(M)} \text{ with } \\
& \,\,\, \, \, \underline{Y_i}= E\backslash\cl(Gi) \text{ for } i\in \{s, t\}.
\end{cases}
\label{ext.phi.form}
\end{equation}
are left quasi-Pl\"ucker coordinates of $\widetilde{\cM}$. In particular, $\mathcal{\widetilde{M}}$ is uniquely determined by $\sigma$.
\end{prop}

\begin{proof} 
Let $Fs$, $Ft$ be adjacent bases of $\underline{\tcM}$ and let $[\cdot]_{\tcM}$ be left quasi-Pl\"ucker coordinates of $\widetilde{\cM}$. 

If $p\notin Fst$, then $Fs$ and $Ft$ are also adjacent bases of $\underline{\cM}$. Then there exists $Y\in \cC^*(\cM)$ such that $\underline{Y} = E\backslash \cl(F)$ and $\tY :=( Y, \sigma(Y))\in \cC^*(\tcM)$.
So $[Fs, Ft]_{\tcM} = \overline{\tY(s) \cdot \tY(t)^{-1}} = \overline{Y(s) \cdot Y(t)^{-1}} = [Fs, Ft]$ by Lemma~\ref{PivotProp}. 

If $p \in \{s, t\}$, then there exists $Y\in \cC^*(\cM)$ such that $\underline{Y} = E\backslash \cl(F)$. If $p=s$, then by Lemma~\ref{pivot}, we have $[Fp, Ft]_{\tcM} = \overline{\sigma(Y)\cdot Y(t)^{-1}}$. If $p=t$, then by \ref{p1} and Lemma~\ref{pivot}, we have $[Fs, Fp]_{\tcM} = ([Fp, Fs]_{\tcM})^{-1} = (\overline{\sigma(Y)\cdot Y(s)^{-1}})^{-1} = \overline{Y(s)\cdot \sigma(Y)^{-1}}$.

If $p\in F$, then let $G = F\backslash \{p\}$. As $Fs, Ft\in \tcB$, then both $Gs$ and $Gt$ are independent sets in $\underline{\cM}$. Now we consider $Gst$. 

If $Gst\notin \tcB$, then we let $g\in E$ such that $Gsg, Gtg \in \tcB$. We also have $Gsg, Gtg \in \cB$. So by \ref{p4} and the above argument, $[Fs, Ft]_{\tcM} = [Gsg, Gtg]_{\tcM} = [Gsg, Gtg].$

If $Gst \in \tcB$, then $Gst\in \cB$. So there exist $Y_i\in \cC^*(\cM)$ such that $\underline{Y_i} = E\backslash \cl(Gi)$ for $i\in \{s,t\}$. As $Gst\in \cB$, then $t\in \underline{Y_s}$ and $s\in \underline{Y_t}$. By \ref{p2}, we have $[Gsp, Gtp]_{\tcM} \cdot [Gst, Gsp]_{\tcM} \cdot [Gtp, Gst]_{\tcM} = -1$. Then by \ref{p1} and the above argument, we have
$$[Gsp, Gtp]_{\tcM} = - [Gst, Gtp]_{\tcM} \cdot [Gsp, Gst]_{\tcM} = -\overline{Y_t(s)\cdot \sigma(Y_t)^{-1}\cdot\sigma(Y_s)\cdot Y_s(t)^{-1}}.$$

So, $[\cdot]_{\widetilde{\cM}}$ is defined as Formula~\eqref{ext.phi.form} and thus uniquely determined by $\sigma$. Therefore, $\mathcal{\widetilde{M}}$ is uniquely determined by $\sigma$.
\end{proof}

If $\cM$ is a right $T$-matroid, then the right quasi-Pl\"ucker coordinates of the extension $\tcM$ are found analogously by reversing the order of multiplication throughout.

In Subsection~\ref{subsect.ext.MOM}, we talked about how to find the signed cocircuit set of the extended oriented matroid, determined by a localization. Now we will generalize that to matroids over skew tracts.

First, for matroid over skew tracts, we would like to find all $X \in \mathcal{C^*(\widetilde{M})}$ that arise via modular elimination by $p$. The following lemma characterizes those $T$-cocircuits.

\begin{lem} \label{mod.elim}
Let $(Y_1, Y_2)$ be a modular pair of $T$-cocircuits of $\cM$ with $\sigma(Y_1) = -\sigma(Y_2) \neq 0$. Then there exists a unique $X \in \mathcal{C^*(\widetilde{M})}$ such that $X(p) = 0$ and $Y_1(f)+Y_2(f) - X(f) \in N_G$ for all $f \in E$. Further, if $e \in \underline{Y_1} \cap \underline{Y_2}$, then
$$X(e) = - Y_1(e) \cdot \sigma(Y_1)^{-1}\cdot \sigma(Y_e)\cdot Y_e(e_1)^{-1} \cdot Y_2(e_1)$$	
where $Y_e \in \mathcal{C^*(M)}$ with $\underline{Y_e}=E\backslash \cl((Y_1^0\cap Y_2^0)\cup e)$ and $e_1\in \underline{Y_2} \backslash \underline{Y_1}$.
\end{lem}

Recall that such a $Y_e$ was guaranteed by Lemma~\ref{exist}.

\begin{proof} The Modular Elimination axiom in $\widetilde{\cM}$ shows that there exists $X \in \mathcal{C^*(\widetilde{M})}$ such that $X$ eliminates $p$ between $(Y_1, \sigma(Y_1))$ and $(Y_2, \sigma(Y_2))$, that is $X(p) = 0$ and $Y_1(f)+Y_2(f) - X(f) \in N_G$ for all $f \in E$.

Now, let $e \in \underline{Y_1} \cap \underline{Y_2}$. Let us find the value of $X(e)$.

Let $F$ be a maximal independent set in $Y_1^0 \cap Y_2^0$. Then $|F| = d-2$. Let $e_1 \in \underline{Y_2} \backslash \underline{Y_1}$. Thus $X(e_1) = Y_2(e_1)$ and $\underline{Y_1} = E\backslash \cl(Fe_1)$. So $Fe_1p, Fe_1e, Fep\in \tcB.$ Let $Y_e \in \mathcal{C^*(M)}$ with $\underline{Y_e} = E\backslash \cl(Fe)$.

By Dual pivoting property, we have that $\overline{X(e) \cdot X(e_1)^{-1}} = [Fep, Fe_1 p]_p$. By Proposition~\ref{ext.uniq}, we also have that 
$$[Fep, Fe_1 p]_p = -\overline{Y_1(e)\cdot \sigma(Y_1)^{-1}\cdot\sigma(Y_e)\cdot Y_e(e_1)^{-1}}. $$
Then $$\overline{X(e) \cdot X(e_1)^{-1}} = -\overline{Y_1(e)\cdot \sigma(Y_1)^{-1}\cdot\sigma(Y_e)\cdot Y_e(e_1)^{-1}}.$$
So 
$$X(e) = - Y_1(e) \cdot \sigma(Y_1)^{-1}\cdot \sigma(Y_e)\cdot Y_e(e_1)^{-1} \cdot Y_2(e_1).$$
\end{proof}

If $\cM$ is a right $T$-matroid, then the value of $X(e)$ is obtained by reversing the order of multiplication throughout.

\begin{defn}\label{cocir.mod.elim}
 We denote $X$ as in Lemma~\ref{mod.elim} by $\Mod (Y_1, Y_2, p)$.
\end{defn}

\begin{prop} \label{ext.cocir}
The set of $T$-cocircuits of $\mathcal{\widetilde{M}}$ is
\begin{align*}
\mathcal{C^*(\widetilde{M})} = & \{(Y, \sigma(Y))\,|\, Y \in \mathcal{C^*(M)}\} \cup \\
& \{ \Mod(Y_1, Y_2, p)\,|\, (Y_1, Y_2) \text{ is a modular pair of $T$-cocircuits in $\mathcal{C^*(M)}$, and }\sigma(Y_1) = -\sigma(Y_2) \neq 0 \}.
\end{align*}
\end{prop}

\begin{proof}
($\supseteq$): Proposition~\ref{sigma.intro} and Lemma~\ref{mod.elim} show that the right side of the equation is contained in $\mathcal{C^*(\widetilde{M})}$.

($\subseteq$): As $\underline{\mathcal{\widetilde{M}}\backslash p} = \underline{\mathcal{\widetilde{M}}}\backslash p$, then $\underline{\mathcal{\widetilde{M}}}$ is also a single-element extension of $\underline{\cM}$ by $p$.

Let $X \in \mathcal{C^*(\widetilde{M})}$. Thus $X^0$ is a hyperplane of $\underline{\tcM}$. Lemma~\ref{ext.hypl} describes 2 possible types for $X^0$.

If $X^0$ is of type (1), then either $X^0$ or $X^0 \backslash p$ is a hyperplane of $\underline{\cM}$. So there exists $Y' \in \mathcal{C^*(M)}$ such that $\underline{Y'} = \underline{X} \backslash p$ or $\underline{Y'} = \underline{X}$. As $(Y', \sigma(Y')) \in \mathcal{C^*(\widetilde{M})}$, by Incomparability, $X = (Y', \sigma(Y'))\cdot \alpha$, for some $\alpha \in T^\times$. So $X \in \{(Y, \sigma(Y))\,|\, Y \in \mathcal{C^*(M)}\}$.

If $X^0 $ is of type (2), then $p \notin \underline{X}$. That is, $X = (Z, 0)$, for some $Z \in T^E$. And $Z^0 = H_1 \cap H_2$, where $H_1$ and $H_2$ are two hyperplanes of $\underline{\cM}$ and there is no hyperplane $H \supset Z^0$ of $\underline{\cM}$ such that $H \cup p$ is a hyperplane of $\underline{\tcM}$.

Let $Y_1, Y_2 \in \mathcal{C^*(M)}$ such that $Y_1^0 = H_1$ and $Y_2^0 = H_2$. Then $\underline{Z} = \underline{Y_1} \cup \underline{Y_2}$. As $\rank(H_1) = \rank(H_2) = d-1$, then $\rank(H_1 \cap H_2) \leq d-2$. Since $\rank(X^0) = d-1$ and $H_1 \cap H_2 = Z^0 = X^0 \backslash p$, then $\rank(H_1 \cap H_2) \geq d-2$. So $\rank(H_1 \cap H_2) = d-2$. Thus $(Y_1, Y_2)$ is a modular pair of $T$-cocircuits. As $H_i \supset Z^0$, $H_i \cup p$ is not a hyperplane of $\underline{\tcM}$, for $i=1,2$. So $\sigma(Y_i) \neq 0$, for $i=1,2$. Without loss of generality, choose $Y_1, Y_2$ such that $\sigma(Y_1) = -\sigma(Y_2) \neq 0$ and $Y_1(e) = Z(e)$, for $e \in \underline{Y_1} \backslash \underline{Y_2}$.

By Modular Elimination in $\widetilde{\cM}$, there exists $X_1 \in \mathcal{C^*(\widetilde{M})}$ such that $X_1(p) = 0$ and for every $e \in E$, $Y_1(e) + Y_2(e) - X_1 (e) \in N_G$. So $\underline{X_1} \subseteq \underline{Y_1} \cup\underline{Y_2}$.

As $\underline{X} = \underline{Z} = \underline{Y_1} \cup \underline{Y_2}$, $\underline{X_1} \subseteq \underline{X}$. By Incomparability, $X = X_1\cdot \alpha$, for some $\alpha \in T^\times$. As $X_1(e) = Y_1(e)$, for all $e \in\underline{Y_1} \backslash \underline{Y_2}$, then $ \alpha = 1$. So $X = X_1$.

So $X = \Mod(Y_1, Y_2, p)$. 
\end{proof}

\begin{cor}\label{sup.mod} The support of $\Mod(Y_1, Y_2, p)$ is equal to $\underline{Y_1} \cup \underline{Y_2}$.
\end{cor}

It follows from the proof of Proposition~\ref{ext.uniq} and the proof of Proposition~\ref{ext.cocir} that the characterization of quasi-Pl\"ucker coordinates $[\cdot]_p$ and the characterization of $T$-cocircuit set $\cC^*(\tcM)$ also hold if $\cM$ is a strong $T$matroid.



\section{Pathetic Cancellation Property}\label{sect.PCP}
Recall (Theorem~\ref{OM.ext}) that if a function $\sigma: \mathcal{C^*(M)} \rightarrow \mathbb{S}$ defines a single-element extension of every rank 2 contraction of an oriented matroid $\cM$, then it is a localization. But it is not always true for other skew tracts. 

Now we will introduce a property for a skew tract $T$, which is necessary and sufficient for Theorem~\ref{OM.ext} to generalize to weak matroids over skew tracts (Theorem~\ref{maintheorem}). Theorem~\ref{counterex} proves the necessity and Theorem~\ref{weakext.equi} proves the sufficiency.

\begin{defn}\label{pathetic}
Let $T$ be a skew tract. We say $T$ satisfies the \textbf{Pathetic Cancellation Property} if for every $a, b, x, y, z \in T^\times$ with 
\begin{align*}
 1 + a- x &\in N_G, \\
 ax & = xa, \\
-1 + b- y &\in N_G,\\
by & = yb,\\
a + b - z &\in N_G,\\
a^{-1}zb^{-1} & = b^{-1}za^{-1},\\ 
x + y - z &\in N_G,\\
x^{-1}zy^{-1} & = y^{-1}zx^{-1},
\end{align*}
we have 
$$xb -ay -z \in N_G.$$
\end{defn}

The name \emph{Pathetic Cancellation} comes from how the property behaves in a field. When $T$ is a field, the hypothesis can be stated as 
$$x=1+a, y=-1 + b,z=a + b, z= x + y.$$
and the conclusion can be stated as $z = xb-ay.$

It is a \emph{cancellation} in a field as
$$xb -ay = (1+ a)b -a(-1+ b) = (a+ b) + (ab -ab) = a+ b = z.$$ And it is {\em Pathetic} because it needs all hypotheses, but when $T$ is a field, the hypotheses become less and any three of them imply the fourth one and give the conclusion via cancellation.

The following example shows that $\mathbb{P}$ does not satisfy Pathetic Cancellation. 



\begin{example}\label{exam2}
Let $a, b, x, y, z$ be elements of $\mathbb{P}$ such that $a = \frac{-1+i}{\sqrt{2}}$, $b = \frac{1+i}{\sqrt{2}}$, $x = \frac{\sqrt{3}+i}{2}$, $y= \frac{-1+\sqrt{3}i}{2}$, and $z = i$. It is easy to check that $a, b, x, y, z$ satisfies all hypotheses in Definition~\ref{pathetic}. However, $z\notin xb\boxplus -ay = \{\frac{-1+\sqrt{3}+i(1+\sqrt{3})}{2\sqrt{2}}\}$, as shown in Figure~\ref{ph.reln}.

We depict $a$, $b$, $x$, $y$ and $z$ by labelled points on a picture of $S^1$. The left circle in Figure~\ref{ph.reln} depicts $a$, $b$, $x$, $y$ and $z$. The dashed lines indicate the hyperadditions $a\boxplus 1$, $-1\boxplus b$, $a\boxplus b$ and $x\boxplus y$. And the right circle depicts $-ay$ and $xb$.


\end{example}

\begin{figure}[htb]
\centering
\begin{tikzpicture}
\draw (3,0) arc(0:180:3); 
\draw[fill] (3,0) circle (.08); \node at (3,-0.4) {$1$};
\draw[fill] (-3,0) circle (.08); \node at (-3,-0.4) {$-1$};
\draw[fill] (-2.121,2.121) circle (.08); \node at (-2,1.8) {$a$};
\draw[fill] (2.121,2.121) circle (.08); \node at (2,1.8) {$b$};
\draw[->] (-3.5,0) -- (3.5,0);
\draw[->] (0, -0.5) -- (0,4);
\draw[dashed] (3.2,0) arc(0:135:3.2);
\draw[fill] (2.77,1.6) circle (.08);\node at (3.07,1.7) {$x$};
\draw[dashed] (-3.4,0) arc(180:45:3.4);
\draw[fill] (-1.7, 2.94) circle (.08);\node at (-2, 3.2) {$y$};
\draw[dashed] (1.98,1.98) arc(45:135:2.8); 

\draw[dashed] (3.12,1.8) arc(30:120:3.6);

\draw[fill] (0, 2.8) circle (.08); \node at (0.2,2.6) {$z$};

\draw (9,1) arc(0:180:2);
\draw[fill] (7, 3) circle (.08); \node at (6.8, 3.2) {$z$};
\draw[fill] (7.517, 2.932) circle (.08); \node at (8.2,3.2) {$-ay = xb$};
\draw[->] (4.5,1) -- (9.5,1);
\draw[->] (7, 0.5) -- (7,3.5);
\end{tikzpicture}
\caption{Illustration for Example~\ref{exam2}}
\label{ph.reln}
\end{figure}



\begin{defn}
Let $\cM$ be a $T$-matroid. We say that a function $\sigma: \cC^*(\cM)\rightarrow T$ \textbf{defines} a single-element extension of a minor of $\cM$ if and only if the map induced from $\sigma$ is a localization of this minor.
\end{defn}

\begin{thm}\label{counterex} 
Let $T$ be a skew tract not satisfying Pathetic Cancellation. There is a both weak and strong left $T$-matroid $\cM$ on $E$ and a right $T^\times$-equivariant function $\sigma: \mathcal{C^*(M)} \rightarrow T$ such that $\sigma$ defines a single-element extension by a new element $p$ on every rank 2 contraction of $\cM$ but $\sigma$ is not a localization.
\end{thm}

Before showing the proof, we will first introduce a useful lemma for determining localizations of rank 2 $T$-matroids on three elements.

\begin{lem}\label{rank2.elim.loc} Let $T$ be a skew tract, let $\cM$ be a uniform rank 2 left $T$-matroid on three elements with $T$-cocircuit set $\cC^*(\cM)$, and let $\sigma: \mathcal{C^*(M)} \rightarrow T$ be a right $T^\times$-equivariant function. Then $\sigma$ is a localization if and only if there exist $Y_1, Y_2, Y_3 \in \mathcal{C^*(M)}$ such that $Y_3$ eliminates some element $e$ between $Y_1$ and $Y_2$, and $\sigma(Y_1) + \sigma(Y_2) - \sigma(Y_3) \in N_G$.
\end{lem}
\begin{proof} Let $E= \{e_1, e_2, e_3\}$ be the ground set of $\cM$.

($\Rightarrow$): We suppose that $\sigma$ is a localization. As $\cM$ is a rank 2 $T$-matroid, then every pair of $T$-cocircuits of $\cM$ form a modular pair. Let $e\in E$ and $(Y_1, Y_2)$ be a modular pair of $T$-cocircuits of $\cM$ with $Y_1(e) = -Y_2(e) \neq 0$. Without loss of generality, we may assume $e = e_3$, $Y_1(e_1) = 0$ and $Y_2(e_2) = 0$. There exists $Y_3\in \mathcal{C^*(M)}$ such that $Y_3$ eliminates $e_3$ between $Y_1$ and $Y_2$. 
Thus $Y_3(e_3) = 0$ and $Y_1(f) + Y_2(f) - Y_3(f) \in N_G$ for all $f\in E$. So $Y_3(e_1) = Y_2(e_1)$, $Y_3(e_2) = Y_1(e_2)$ and $\underline{Y_3} = \{e_1, e_2\}$. 

Let $\widetilde{\cM}$ be the extension of $\cM$ on $E\cup p$ determined by $\sigma$. As $(Y_1, \sigma(Y_1))$,$(Y_2, \sigma(Y_2))$ form a modular pair in $\cC^*(\widetilde{\cM})$ with $Y_1(e) = -Y_2(e)$, then there exists $\widetilde{Y}\in \cC^*(\widetilde{\cM})$ such that $\widetilde{Y}$ eliminate $e$ between $(Y_1, \sigma(Y_1))$ and $(Y_2, \sigma(Y_2))$. Thus $\widetilde{Y}(e)=0$, $Y_1(f)+Y_2(f)-\widetilde{Y}(f)\in N_G$ for all $f\in E$, and $\sigma(Y_1) + \sigma(Y_2)-\widetilde{Y}(p)\in N_G$. So $\widetilde{Y}(e_1) = Y_2(e_1)$, $\widetilde{Y}(e_2) = Y_1(e_2)$ and $\{e_1, e_2\} \subseteq \underline{\widetilde{Y}} \subseteq \{e_1, e_2, p\}$. 

Now let us look at the $T$-cocircuits $(Y_3, \sigma(Y_3))$ and $\widetilde{Y}$. As $\{e_1, e_2\} \subseteq \underline{(Y_3, \sigma(Y_3))} \subseteq \{e_1, e_2, p\}$ and $\{e_1, e_2\} \subseteq \underline{\widetilde{Y}} \subseteq \{e_1, e_2, p\}$, then either $\underline{(Y_3, \sigma(Y_3))} \subseteq \underline{\widetilde{Y}}$ or $\underline{\widetilde{Y}} \subseteq \underline{(Y_3, \sigma(Y_3))}$. By Symmetry, $\widetilde{Y} = (Y_3, \sigma(Y_3))\cdot \alpha$ for some $\alpha\in T^\times$. As $Y_3(e_1) = Y_2(e_1) = \widetilde{Y}(e_1)$, then $\alpha = 1$. Thus $\widetilde{Y} = (Y_3, \sigma(Y_3))$. So $\sigma(Y_1) + \sigma(Y_2)- \sigma(Y_3) \in N_G$.

($\Leftarrow$): We suppose that there exist $Y_1, Y_2, Y_3 \in \mathcal{C^*(M)}$ such that $Y_3$ eliminates some element $e$ between $Y_1$ and $Y_2$, and $\sigma(Y_1) + \sigma(Y_2)- \sigma(Y_3) \in N_G$. Then without loss of generality, we may assume $e=e_3$, $Y_1(e_1) = 0$ and $Y_2(e_2) = 0$. We list a representative from each projective $T$-cocircuit of $\cM$ in Table~\ref{t1}.

\begin{table}[htb]
\centering
\begin{tabular}{ |l | c |c|c| }\hline
& $e_1$ & $e_2$ & $e_3$  \\ \hline 
  $Y_2$ &  $Y_2(e_1)$ & 0 & $-Y_1(e_3)$ \\ \hline
 $Y_1$ &  0 & $Y_1(e_2)$ & $Y_1(e_3) $  \\ \hline
 $Y_3$ & $Y_2(e_1)$ & $Y_1(e_2)$ & 0  \\ \hline
\end{tabular}
\smallskip
\caption{}
\label{t1}
\end{table}

Let $Z \in T^{E\cup p}$ be such that $$Z = (Y_2(e_1)\cdot \sigma(Y_2)^{-1}, - Y_1(e_2) \cdot \sigma(Y_1)^{-1}, -Y_1(e_3)\cdot (\sigma(Y_1)^{-1} + \sigma(Y_2)^{-1}), 0)$$ 
and let $\alpha\in T^\times$ be such that $\alpha = \sigma(Y_1)^{-1} + \sigma(Y_2)^{-1}$. Then we define a set $\tcC \subseteq T^{E\cup p}$ by
$$\widetilde{\cC} := \{(Y_1, \sigma(Y_1))\cdot \beta\,|\, \beta\in T^\times \}\cup \{(Y_2, \sigma(Y_2))\cdot \beta\,|\, \beta\in T^\times \}\cup \{(Y_3, \sigma(Y_3))\cdot \beta \,|\, \beta\in T^\times \}\cup \{Z\cdot \beta \,|\, \beta\in T^\times \}.$$

We would like to show that $\widetilde{\cC}$ is the $T$-cocircuit set of a rank 2 left $T$-matroid $\widetilde{\cM}$ on the set $E\cup p$ and $\widetilde{\cM}$ is the extension of $\cM$ determined by $\sigma$.

As $\cM$ is of rank 2, then every pair of $T$-cocircuits of $\cM$ is a modular pair.

As $\sigma(Y_1) + \sigma(Y_2)- \sigma(Y_3) \in N_G$, it is easy to verify that 
\begin{align*}
Z & = \Mod(-Y_1 \cdot \sigma(Y_1)^{-1}, Y_2 \cdot \sigma(Y_2)^{-1}, p),\\
Z & = \Mod(Y_2 \cdot \alpha, -Y_3\cdot \sigma(Y_1)^{-1}, p), \\
Z & = \Mod(-Y_1 \cdot \alpha, Y_3\cdot \sigma(Y_2)^{-1}, p).
\end{align*}
So $\sigma$ is a localization.
\end{proof}

The statement also holds for a uniform rank 2 right $T$-matroid on three elements and a left $T^\times$-equivariant function.

Now we will give the proof of Theorem~\ref{counterex}.
\begin{proof}[Proof of Theorem~\ref{counterex}] As $T$ does not satisfy Pathetic Cancellation, then there exist $a, b, x, y, z \in T^\times$ with $ 1 + a- x \in N_G$, $ax = xa$, $-1 + b- y \in N_G$, $by = yb$, $a + b - z \in N_G$, $a^{-1}zb^{-1} = b^{-1}za^{-1}$, $x + y - z \in N_G$ and $x^{-1}zy^{-1} = y^{-1}zx^{-1}$, but $xb -ay -z \notin N_G$.

Now we define a weak left $T$-matroid $\cM$ on the set $E = \{y_1, y_2, y_3, y_4\}$ of rank 3 with $T$-cocircuit set $\mathcal{C^*(M)}$. We list a representative from each projective $T$-cocircuit of $\cM$ in Table~\ref{t2}.

\begin{table}[htb]
\centering
\begin{tabular}{|l|c|c|c|c|} \hline
& $y_1$ & $y_2$ & $y_3$ & $y_4$ \\ \hline
$ Y_{23}$ & 1 & 0 & 0 & 1 \\ \hline
$Y_{13}$ & 0 & 1 & 0 & $b^{-1}$ \\ \hline
$Y_{12}$ & 0 & 0 & 1 & $a^{-1}$\\ \hline
$Y_{24}$ & 1 & 0 & $-a$ & 0 \\ \hline
$Y_{34}$ & -1 & $b$ & 0 &0 \\ \hline
 $Y_{14}$ &0 & $b$ & $-a$ &0 \\ \hline
\end{tabular}
\smallskip
\caption{}
\label{t2}
\end{table}

Figure~\ref{collinear} shows the collinearities for cocircuits of the underlying matroid on elements of $E$.

\begin{figure}[htb]
\centering
\begin{tikzpicture}
\draw (-6,0)--(7,0); \node at (-6.3,0) {$y_2$};
\draw (-5.6,-0.6)--(0.8,5.8); \node at (-5.8, -0.8) {$y_3$};
\draw (0,-0.9)--(0,6.2); \node at (0, -1.2) {$y_1$};
\draw[dashed] (1.7,-0.7)--(-4.9,5.9); \node at (1.9, -0.9) {$p$};
\draw (7,-1/2)--(-5,11/2); \node at (7.3, -0.6) {$y_4$};
\draw [fill] (-5,0) circle (.1cm);\node at (-5.3,0.4) {$Y_{23}$};
\draw [fill] (0,5) circle (.1cm);\node at (0.4,4.8) {$Y_{13}$};
\draw [fill] (0,0) circle (.1cm);\node at (-0.3,-0.3) {$Y_{12}$};
\draw [fill] (-4/3,11/3) circle (.1cm);\node at (-1.4,4.1) {$Y_{34}$};
\draw [fill] (6,0) circle (.1cm);\node at (6.3,0.3) {$Y_{24}$};
\draw [fill] (0,3) circle (.1cm);\node at (0.4,3.3) {$Y_{14}$};
\draw [fill] (-2,3) circle (.1cm);\node at (-2.4,3) {$Y_2$};
\draw [fill] (1,0) circle (.1cm);\node at (0.8,-0.3) {$Y_1$};
\draw [fill] (0,1) circle (.1cm);\node at (0.3,1.3) {$Z$};
\draw [fill] (-4,5) circle (.1cm);\node at (-4.3,4.7) {$Y_3$};
\end{tikzpicture}
\caption{Collinearities for the counterexample}
\label{collinear}
\end{figure}

It is easy to check $\mathcal{C^*(M)}$ satisfies Modular Elimination. As $\cM$ is of rank 3 with fewer than 6 elements, $\cM$ is both a weak and a strong left $T$-matroid. 

Let $\sigma: \mathcal{C^*(M)} \rightarrow T$ be the right $T$-equivariant function such that
$$\sigma(Y_{23}) = \sigma(Y_{13}) = -\sigma(Y_{12}) = 1,\sigma(Y_{24}) = x, \sigma(Y_{34}) = y \text{ and } \sigma(Y_{14}) = z.$$

First, we will verify that $\sigma$ defines a single-element extension of every rank 2 contraction of $\cM$ using Lemma~\ref{rank2.elim.loc}.

We can see that $Y_{13}\backslash y_1, Y_{12}\backslash y_1$ and $Y_{14}\backslash y_1$ are representatives for the $T$-cocircuits of the rank 2 contraction $\cM/y_1$. And $Y_{14}\backslash y_1$ eliminates $y_4$ between $(Y_{13}\backslash y_1) \cdot Y_{13}(y_4)^{-1}$ and $(Y_{12}\backslash y_1) \cdot (-Y_{12}(y_4)^{-1})$. We also have
$$\sigma(Y_{13}) \cdot Y_{13}(y_4)^{-1} + \sigma(Y_{12}) \cdot (-Y_{12}(y_4)^{-1}) -\sigma(Y_{14}) = b + a - z \in N_G.$$

So by Lemma~\ref{rank2.elim.loc}, $\sigma$ induces a function $\sigma/y_1$ which is a localization on $\cM/y_1$.

Similarly, $Y_{23}\backslash y_2, Y_{12}\backslash y_2$ and $Y_{24}\backslash y_2$ are representatives for the $T$-cocircuits of the rank 2 contraction $\cM/y_2$. And $Y_{24}\backslash y_2$ eliminates $y_4$ between $Y_{23}\backslash y_2$ and $(Y_{12}\backslash y_2)\cdot (-Y_{12}(y_4)^{-1})$. We also have
$$\sigma(Y_{23}) + \sigma(Y_{12})\cdot (-Y_{12}(y_4)^{-1}) - \sigma(Y_{24})= 1+ a -x \in N_G.$$

So by Lemma~\ref{rank2.elim.loc}, $\sigma$ induces a function $\sigma/y_2$ which is a localization on $\cM/y_2$.

Similarly, $Y_{23}\backslash y_3, Y_{13}\backslash y_3$ and $Y_{34}\backslash y_3$ are representatives for the $T$-cocircuits of the rank 2 contraction $\cM/y_3$. And $Y_{34}\backslash y_3$ eliminates $y_4$ between $-Y_{23}\backslash y_3$ and $(Y_{13}\backslash y_3) \cdot Y_{13}(y_4)^{-1}$. We also have
$$ -\sigma(Y_{23}) + \sigma(Y_{13})\cdot Y_{13}(y_4)^{-1} - \sigma(Y_{34}) = -1+ b - y \in N_G.$$

So by Lemma~\ref{rank2.elim.loc}, $\sigma$ induces a function $\sigma/y_3$ which is a localization on $\cM/y_3$.

Similarly, $Y_{24}\backslash y_4, Y_{34}\backslash y_4$ and $Y_{14}\backslash y_4$ are representatives for the $T$-cocircuits of the rank 2 contraction $\cM/y_4$. And $Y_{14}\backslash y_4$ eliminates $y_1$ between $Y_{24}\backslash y_4$ and $Y_{34}\backslash y_4$. We also have 
$$\sigma(Y_{24}) + \sigma(Y_{34}) -\sigma(Y_{14}) = x + y -z \in N_G.$$

So by Lemma~\ref{rank2.elim.loc}, $\sigma$ induces a function $\sigma/y_4$ which is a localization on $\cM/y_4$.

So $\sigma$ defines a single-element extension of every rank 2 contraction of $\cM$.

Next we would like to show that $\sigma$ is not a localization of $\cM$. Suppose by way of contradiction that $\sigma$ is a localization of $\cM$. Let $\widetilde{\cM}$ be the single-element extension of $\cM$ determined by $\sigma$. Now we list a representative from each projective $T$-cocircuit of $\widetilde{\cM}$ in the form of $(Y,\sigma(Y))$, $Y\in \cC^*(\cM)$ in Table~\ref{t3}.
\begin{table}[htb]
\centering
\begin{tabular}{|l|c|c|c|c|c|} \hline
& $y_1$ & $y_2$ & $y_3$ & $y_4$ & $p$ \\ \hline
$(Y_{23}, \sigma(Y_{23}))$ & 1 & 0 & 0 & 1 &1\\ \hline
$(Y_{13}, \sigma(Y_{13}))$ & 0 & 1 & 0 & $b^{-1}$ &1\\ \hline
$(Y_{12}, \sigma(Y_{12}))$ & 0 & 0 & 1 & $a^{-1}$ & -1\\ \hline
$(Y_{24}, \sigma(Y_{24}))$ & 1 & 0 & $-a$ & 0 & $x$ \\ \hline
$(Y_{34}, \sigma(Y_{34}))$ &-1 & $b$ & 0 & 0 & $y$ \\ \hline
$(Y_{14}, \sigma(Y_{14}))$ & 0 & $b$ & $-a$ & 0 & $z$ \\ \hline
\end{tabular}
\smallskip
\caption{}
\label{t3}
\end{table}

By Proposition~\ref{ext.cocir}, the extension $\mathcal{\widetilde{M}}$ will also have $T$-cocircuits
\begin{align*}
Y_1 & = \Mod(Y_{23}, Y_{12}, p),\\
Y_2 & = \Mod(-Y_{23}, Y_{13}, p),\\
Z & = \Mod(Y_{13}, Y_{12}, p),\\
Y_3 & = \Mod(Y_{24} \cdot \sigma(Y_{24})^{-1}, Y_{34} \cdot (-\sigma(Y_{34})^{-1}), p).
\end{align*}

So by Lemma~\ref{mod.elim}, we have
\begin{align*}
Y_1(y_1) & = Y_{23}(y_1) =1 & \text{ since } & Y_{12}(y_1) = 0, \\
Y_1(y_2) & = 0  & \text{ since } & Y_{23}(y_2) = Y_{12}(y_2) = 0, \\
Y_1(y_3) & = Y_{12}(y_3) =1  & \text{ since } & Y_{23}(y_3) = 0, \\
Y_1(y_4) & = - Y_{23}(y_4) \sigma(Y_{23})^{-1} \sigma(Y_{24}) Y_{24}(y_3)^{-1} Y_{12}(y_3) = xa^{-1},
\end{align*}
and so
$$Y_1= ( 1,  0,  1, xa^{-1}, 0 ).$$

Similarly by Lemma~\ref{mod.elim}, we have
\begin{align*}
Y_2(y_1) & = - Y_{23}(y_1) =-1 & \text{ since } & Y_{13}(y_1) = 0,\\
Y_2(y_2) & = Y_{13}(y_2) =1  & \text{ since } & Y_{23}(y_2) = 0,\\
Y_2(y_3) & = 0  & \text{ since }&  Y_{23}(y_3) = Y_{13}(y_3) = 0,\\
Y_2(y_4) & = - (-Y_{23}(y_4)) \sigma(-Y_{23})^{-1} \sigma(Y_{34}) Y_{34}(y_2)^{-1} Y_{13}(y_2) = -yb^{-1},
\end{align*}
and so
$$
Y_2= ( -1, 1, 0, -yb^{-1} , 0 ).$$

Similarly by Lemma~\ref{mod.elim}, we have
\begin{align*}
Z(y_1) & = 0 & \text{ since } & Y_{13}(y_1) = Y_{12}(y_1) = 0,\\
Z(y_2) & = Y_{13}(y_2)=1  & \text{ since } & Y_{12}(y_2) = 0, \\
Z(y_3) & = Y_{12}(y_3) =1  & \text{ since } & Y_{13}(y_3) = 0, \\
Z(y_4) & = - Y_{13}(y_4) \sigma(Y_{13})^{-1} \sigma(Y_{14}) Y_{14}(y_3)^{-1} Y_{12}(y_3) = b^{-1} z a^{-1},
\end{align*}
and so
$$ Z= ( 0, 1, 1, b^{-1} z a^{-1}, 0 ).$$

We can see that $Y_1^0 = \{y_2, p\}$ and $Y_2^0 = \{y_3, p\}$, thus $Y_1$ and $Y_2$ form a modular pair. As 
$Y_1(y_1) = 1$, $Y_2(y_1) = -1$ and $Z^0 = \{y_1, p\}$, thus $Z$ should eliminate $y_1$ between $Y_1$ and $Y_2$. Then
$$ Y_1(y_2) + Y_2(y_2) - Z(y_2) = 0 + 1 - 1 = 1 - 1 \in N_G,$$
$$ Y_1(y_3) + Y_2(y_3) - Z(y_3) = 1 + 0 - 1 = 1 - 1 \in N_G,$$
but $$ Y_1(y_4) + Y_2(y_4) - Z(y_4) = xa^{-1} -yb^{-1} - b^{-1} z a^{-1} = a^{-1} x - yb^{-1} - a^{-1} z b^{-1}= a^{-1} (xb-ay - z) b^{-1} \notin N_G,$$ by Pathetic Cancellation. We got a contradiction.

So $\sigma$ is not a localization of $\cM$.
\end{proof}

There is also an example of a right $T$-matroid and a left $T^\times$-equivariant function defined analogously.

This proves the necessity of Pathetic Cancellation for extending Crapo and Las Vergnas' results. We will prove the sufficiency in the next section. 

\section{Characterization of localization for weak $T$-matroid}\label{sect.chara}

Our main theorem (Theorem~\ref{weakext.equi}) in this section proves the converse of Theorem~\ref{counterex}. 

\begin{thm}\label{weakext.equi}
Let $T$ be a skew tract satisfying Pathetic Cancellation, let $\cM$ be a weak left $T$-matroid on $E$, and let
$$\sigma: \mathcal{C^*(M)} \rightarrow T$$ 
be a right $T^\times$-equivariant function. 

Then the following statements are equivalent.
\begin{enumerate}[(1)]
\item $\sigma$ is a weak localization of $\cM$.

\item $\sigma$ defines a weak single-element extension of every rank 2 contraction of $\cM$.

\item $\sigma$ defines a weak single-element extension of every rank 2 minor of $\cM$ on three elements.
\end{enumerate}
\end{thm}

The statement also holds if $\cM$ is a weak right $T$-matroid and $\sigma$ is left $T^\times$-equivariant.

Before showing the main theorem, we will first introduce some useful lemmas.

\begin{lem}\label{ext.phi.indep} Let $\cM$ be a (weak or strong) left $T$-matroid of rank $d$ on $E$, let $[\cdot]$ be left quasi-Pl\"ucker coordinates of $\cM$, and let $\sigma: \cC^*(\cM) \rightarrow T$ be a right $T^\times$-equivariant function. Let $\cl$ be the closure operator of $\underline{\cM}$.
Then for any $F\subseteq E$ independent in $\underline{\cM}$ with $|F|=d-1$, $y\in E\backslash \cl(F)$, $Y\in \cC^*(\cM)$ such that $\underline{Y} = E\backslash \cl(F)$, we have 
$$Y(y) \cdot \sigma(Y)^{-1}$$
is independent of the choice of $Y$.
\end{lem}
\begin{proof} Let $F$ be an independent subset of $E$ with $|F|=d-1$ and let $y\in E\backslash \cl(F)$. Then there exists a unique projective $T$-cocircuit $D$ such that $\underline{Y} = E\backslash \cl(F)$ for any $T$-cocircuit $Y\in D$.

Let $Y_1, Y_2\in D$. Then by definition, there exists $\alpha \in T^\times$ such that $Y_1 = Y_2\cdot \alpha$. Moreover, $y\in \underline{Y_1} = \underline{Y_2}$. We have
\begin{align*}
Y_2(y) \cdot \sigma(Y_2)^{-1} & = Y_2(y) \cdot \alpha \cdot \alpha^{-1} \cdot \sigma(Y_2)^{-1} \\
& = (Y_2(y) \cdot \alpha) \cdot (\sigma(Y_2)\cdot \alpha)^{-1} \\
& = (Y_2(y) \cdot \alpha) \cdot \sigma(Y_2\cdot \alpha)^{-1} \\
& = Y_1(y) \cdot \sigma(Y_1)^{-1} .
\end{align*}

So $Y(y) \cdot \sigma(Y)^{-1}$ is independent of the choice of $Y\in D$.
\end{proof}

So for a weak left $T$-matroid $\cM$ of rank $d$ with left quasi-Pl\"ucker coordinates $[\cdot]$ and a right $T^\times$-equivariant function $\sigma: \cC^*(\cM) \rightarrow T$, we can define a collection $\tcB$ of subsets of $E\cup p$ by 
$$\tcB = \cB \cup \{Fp \subseteq E\cup p \,|\, |F| = d-1, E \backslash \cl(F) = \underline{Y} \text{ and } \sigma(Y) \neq 0 \text{ for some } Y \in \cC^*(\cM) \},$$
where $\cl$ is the closure operator and $\cB$ is the set of bases of $\underline{\cM}$.

It is trivial to see that every set in $\tcB$ has the same size. Then we define the collection $A_{\tcB}$ of adjacent sets by $A_{\tcB} : = \{(B, B')\in \tcB \times \tcB\, |\, |B\backslash B'| =1\} \subseteq T^{E\cup p}$ and then define left $T$-coordinates $[\cdot]_p: A_{\tcB} \rightarrow T$ by
$$ [Fs, Ft]_p=\begin{cases}
[Fs, Ft]	& \text{ if } p\notin Fst,\\
\overline{\sigma(Y)\cdot Y(t)^{-1}}  & \text{ if } p = s\text{ and } Y\in \mathcal{C^*(M)} \text{ with } \underline{Y}= E\backslash\cl(F),\\
\overline{ Y(s)\cdot \sigma(Y)^{-1}}  & \text{ if } p = t \text{ and } Y\in \mathcal{C^*(M)} \text{ with } \underline{Y}= E\backslash\cl(F),\\
[Gsg, Gtg] & \text{ if } p\in F, G=F\backslash \{p\}, Gst\notin \cB \text{ and } Gsg, Gtg\in \cB \\
&\,\,\, \text{ for some } g\in E, \\
-\overline{Y_t(s)\cdot \sigma(Y_t)^{-1}\cdot \sigma(Y_s)\cdot Y_s(t)^{-1}} & \text{ if } p\in F, G=F\backslash \{p\}, Gst\in \cB, Y_i\in \mathcal{C^*(M)} \text{ with } \\
& \,\,\, \, \, \underline{Y_i}= E\backslash\cl(Gi) \text{ for } i\in \{s,t\}.
\end{cases}
$$

\begin{lem}\label{phi.loc.equi}
$[\cdot]_p$ satisfies \ref{p3} and \ref{p5} in weak left quasi-Pl\"ucker coordinates axioms if and only if $\sigma$ is a localization of $\cM$. 
\end{lem}
\begin{proof} ($\Leftarrow$): This is proved by Proposition~\ref{ext.uniq}.

($\Rightarrow$): By definition, $[\cdot]_p$ satisfies \ref{p1}, \ref{p2} and \ref{p4} in weak left quasi-Pl\"ucker coordinates axioms. Then if $[\cdot]_p$ satisfies \ref{p3} and \ref{p5} in weak left quasi-Pl\"ucker coordinates axioms, $[\cdot]_p$ are weak left quasi-Pl\"ucker coordinates of $\underline{\tcM}$. It is easy to check that $[\cdot]_p$ satisfies the Dual Pivoting Property with $\cC^*(\cM)$.

So $\sigma$ is a localization.
\end{proof}

Now we would like to introduce a lemma extending one direction of Lemma~\ref{rank2.elim.loc} to weak left $T$-matroids. 

\begin{lem}\label{loc.add}
Let $\cM$ be a (weak or strong) left $T$-matroid on $E$ and let $\sigma: \mathcal{C^*(M)} \rightarrow T$ be a localization of $\cM$. Let $(Y_1, Y_2)$ be a modular pair of $T$-cocircuits of $\cM$ with $Y_1(e) = -Y_2(e) \neq 0$ for some $e\in E$. Let $Y_3 \in \mathcal{C^*(M)}$ eliminate $e$ between $Y_1$ and $Y_2$. Then $\sigma(Y_1) + \sigma(Y_2) - \sigma(Y_3) \in N_G$.
\end{lem}
\begin{proof} We denote $(Y_i,\sigma(Y_i))$ by $\widetilde{Y_i}$ for $i\in \{1,2\}$. As $(Y_1, Y_2)$ is a modular pair of $T$-cocircuits of $\cM$ with $Y_1(e) = -Y_2(e) \neq 0$, then $\widetilde{Y_1}$ and $\widetilde{Y_2}$ form a modular pair of $T$-cocircuits of $\widetilde{\cM}$ with $\widetilde{Y_1}(e) = -\widetilde{Y_2}(e) \neq 0$. Let $Z\in \cC^*(\widetilde{\cM})$ eliminate $e$ between $\widetilde{Y_1}$ and $\widetilde{Y_2}$. Thus $Z(e)=0$ and $\widetilde{Y_1}(f) + \widetilde{Y_2}(f)- Z(f) \in N_G$ for all $f\in E\cup p$. Then $\sigma(Y_1) + \sigma(Y_2)- Z(p) \in N_G$.

We would like to show that $Z = (Y_3, \sigma(Y_3))$, from which it follows that $\sigma(Y_1) + \sigma(Y_2) - \sigma(Y_3) \in N_G$.

By Proposition~\ref{ext.cocir}, we know that
\begin{align*}
\mathcal{C^*(\widetilde{M})} = & \{(W, \sigma(W))\,|\, W \in \mathcal{C^*(M)}\} \cup \\
& \{ \Mod(W_1, W_2, p)\,|\, (W_1, W_2) \text{ is a modular pair of $T$-cocircuits in $\mathcal{C^*(M)}$, and } \sigma(W_1) = -\sigma(W_2) \neq 0 \}.
\end{align*} 

We first claim that $Z\notin \{ \Mod(Y_1, Y_2, p)\,|\, Y_1, Y_2 \in \mathcal{C^*(M)}$, $(Y_1, Y_2)$ is a modular pair of $T$-cocircuits, and $ \sigma(Y_1) = -\sigma(Y_2) \neq 0 \}$. Otherwise, there exist $(Z_1, Z_2)$ a modular pair of $T$-cocircuits of $\cM$ such that $\sigma(Z_1) = -\sigma(Z_2) \neq 0$ and $Z=\Mod(Z_1, Z_2,p)$. By Corollary~\ref{sup.mod}, $\underline{Z} = \underline{Z_1}\cup \underline{Z_2}$. As $Z$ eliminates $e$ between $\widetilde{Y_1}$ and $\widetilde{Y_2}$, then $\underline{Z} \subseteq (\underline{\widetilde{Y_1}}\cup \underline{\widetilde{Y_2}})\backslash \{e\}$. As $Z(p)=0$, then $\underline{Z} \subseteq (\underline{Y_1}\cup \underline{Y_2})\backslash \{e\}\subset \underline{Y_1}\cup \underline{Y_2}$. So $\underline{Z_1}, \underline{Z_2}\subset \underline{Z_1}\cup \underline{Z_2} \subset \underline{Y_1}\cup \underline{Y_2}.$ As $(Y_1, Y_2)$ is a modular pair of $T$-cocircuits of $\cM$, by Definition~\ref{mod}, the height of $\underline{Y_1}\cup \underline{Y_2}$ is 2 in the lattice $U(\cC^*(\cM))$. We got a contradiction. So our claim holds.

Then from our claim, we have that $Z\in \{(Y, \sigma(Y))\,|\, Y \in \mathcal{C^*(M)}\}$. Let $x\in \underline{Y_1}\backslash \underline{Y_2}$. Then there exists $Z'\in \mathcal{C^*(M)}$ such that $Z = (Z',\sigma(Z'))$. As $Z$ eliminates $e$ between $\widetilde{Y_1}$ and $\widetilde{Y_2}$, then $\underline{Z} \subseteq (\underline{\widetilde{Y_1}}\cup \underline{\widetilde{Y_2}})\backslash \{e\}$ and $Z(x)= Y_1(x)$. Thus $\underline{Z'} \subseteq (\underline{Y_1}\cup \underline{Y_2})\backslash \{e\}$ and $Z'(x)= Y_1(x)$.

As $Y_3$ eliminates $e$ between $Y_1$ and $Y_2$, then $\underline{Y_3} \subseteq (\underline{Y_1}\cup \underline{Y_2})\backslash \{e\}$ and $Y_3(x)= Y_1(x)$.

Suppose by way of contradiction that $Z' \neq Y_3$. As $\underline{Z'} \subseteq (\underline{Y_1}\cup \underline{Y_2})\backslash \{e\}$ and $\underline{Y_3} \subseteq (\underline{Y_1}\cup \underline{Y_2})\backslash \{e\}$, then $\underline{Z'}\cup \underline{Y_3}\subseteq (\underline{Y_1}\cup \underline{Y_2})\backslash \{e\} \subset \underline{Y_1}\cup \underline{Y_2}$. As $(Y_1, Y_2)$ is a modular pair of $T$-cocircuits of $\cM$, by Definition~\ref{mod}, the height of $\underline{Y_1}\cup \underline{Y_2}$ is 2 in the lattice $U(\cC^*(\cM))$, contradicting that
$$\underline{Z'}, \underline{Y_3}\subset \underline{Z'}\cup \underline{Y_3} \subset \underline{Y_1}\cup \underline{Y_2}.$$
So $Z'=Y_3$, and then $Z = (Y_3,\sigma(Y_3))$.

So $\sigma(Y_1) + \sigma(Y_2) - \sigma(Y_3) \in N_G$.
\end{proof}

The statement also hold if $\cM$ is a (weak or strong) right $T$-matroid.

Then following is a useful lemma showing that $\sigma \cdot \alpha$ is a localization if $\sigma$ is, for all $\alpha\in T^\times$.

\begin{lem}\label{loc.eq}
Let $\cM$ be a (weak or strong) left $T$-matroid on $E$, let $\sigma: \cC^*(\cM) \rightarrow T$ be a localization, and let $\widetilde{\cM}$ be the extension of $\cM$ by an element $p$ determined by $\sigma$. Let $\alpha \in T^\times$.
Then $\alpha \cdot \sigma$ is also a localization of $\cM$ and the
corresponding extension is the left $T$-matroid $\widetilde{\cM}^\rho$ arsing from $\tcM$ by right rescaling, where $\rho: E\cup p \rightarrow T^\times$ is a function defined by 
$$\rho(e) = \begin{cases}
1 & \text{if } e\in E, \\
\alpha & \text{if } e = p.
\end{cases}
$$
\end{lem}

\begin{proof} First we consider the left $T$-matroid $\widetilde{\cM}^\rho$ arsing from $\tcM$ by right rescaling. By Lemma~\ref{rescaling}, we know that $\cC^*(\widetilde{\cM}^\rho) = \rho \cdot \cC^*(\tcM)$. By definition of $\rho$, we have
$$\cC^*(\widetilde{\cM}^\rho\backslash p)  = (\rho \cdot \cC^*(\tcM))/p = \MS (\{(\rho \cdot Y)\backslash p \,|\, Y \in \cC^*(\tcM)\}) = \MS (\{Y\backslash p \,|\, Y \in \cC^*(\tcM)\}) =  \cC^*(\cM).
$$
So by definition, $\widetilde{\cM}^\rho$ is also a weak extension of $\cM$ by the same element $p$.

Now we would like to show that $\alpha \cdot \sigma$ is the localization corresponding to $\widetilde{\cM}^\rho$. Let $Y\in \cC^*(\cM)$. By definition, we know that $(Y, \sigma(Y)) \in \cC^*(\tcM)$. Then 
$$\cC^*(\widetilde{\cM}^\rho) \ni \rho \cdot (Y, \sigma(Y)) = (Y, \alpha \cdot \sigma(Y)).$$
So by Proposition~\ref{sigma.intro}, $\alpha \cdot \sigma$ is the localization corresponding to $\widetilde{\cM}^\rho$.
\end{proof}

Now we would like to prove Theorem~\ref{weakext.equi}. But the proof will be a very long calculation. So we will first introduce a lemma which shows a general version of the direction $(1) \Rightarrow (2)$.

\begin{lem} \label{ext.contr}
Let $\cM$ be a (weak or strong) left $T$-matroid on $E$ and let $\sigma: \mathcal{C^*(M)} \rightarrow T$ be a localization of $\cM$. Then $\sigma$ defines a single-element extension of every contraction of $\cM$.
\end{lem}

\begin{proof} Let $\widetilde{\cM}$ be the extension of $\cM$ by an element $p$ determined by $\sigma$ and let $E_0 \subseteq E$. We would like to consider the contraction $\cM/E_0$ of $\cM$. By Theorem~\ref{minor}, we know that
$$\cC^*((\tcM/E_0)\backslash p) = \cC^*((\tcM\backslash p)/E_0) = \cC^*(\cM/E_0).$$
So $\tcM/E_0$ is an extension of $\cM/E_0$ by the element $p$.

Let $\sigma': \cC^*(\cM/E_0) \rightarrow T$ be the function induced from $\sigma$ by contraction of $E_0$. We would like to show that $\sigma'$ is the localization of $\cM/E_0$ corresponding to $\tcM/E_0$. It suffices to show that
\begin{align*}
\cC^*(\widetilde{\cM}/E_0) = \{&(Z, \sigma'(Z))\,|\, Z \in \cC^*(\cM/E_0)\} \cup \\
 \{ &\Mod(Z_1, Z_2, p)\,|\, (Z_1, Z_2) \text{ is a modular pair of $T$-cocircuits of } \cM/E_0 \text{ and } \sigma(Z_1) = -\sigma(Z_2) \neq 0 \}.
\end{align*}

By Theorem~\ref{minor} and Proposition~\ref{ext.cocir},
\begin{align*}
\cC^*(\widetilde{\cM}/E_0) = & \cC^*(\tcM)\backslash E_0 \\
= &\{ Y\backslash E_0 \,|\, Y \in \cC^*(\widetilde{\cM}), \underline{Y} \cap E_0 = \emptyset \}\\
= & \{(Y, \sigma(Y))\backslash E_0\,|\, Y \in \mathcal{C^*(M)}, \underline{Y} \cap E_0 = \emptyset\} \cup \{ \Mod(Y_1, Y_2, p) \backslash E_0\,|\, (Y_1, Y_2) \text{ is a modular pair}\\
& \text{ of $T$-cocircuits of } \cM, \sigma(Y_1) = -\sigma(Y_2) \neq 0 \text{ and } \underline{\Mod(Y_1, Y_2, p)} \cap E_0 = \emptyset \}.\\
=& \{(Y\backslash E_0, \sigma(Y))\,|\, Y \in \mathcal{C^*(M)}, \underline{Y} \cap E_0 = \emptyset\} \cup \{ \Mod(Y_1, Y_2, p) \backslash E_0\,|\, (Y_1, Y_2) \text{ is a modular pair}\\
& \text{ of $T$-cocircuits of } \cM, \sigma(Y_1) = -\sigma(Y_2) \neq 0 \text{ and } (\underline{Y_1}\cup \underline{Y_2}) \cap E_0 = \emptyset \}.\\
=&\{(Z, \sigma'(Z))\,|\, Z \in \cC^*(\cM/E_0)\} \cup 
 \{ \Mod(Z_1, Z_2, p)\,|\, (Z_1, Z_2) \text{ is a modular pair of $T$-cocircuits}\\
 & \text{ of } \cM/E_0 \text{ and } \sigma(Z_1) = -\sigma(Z_2) \neq 0 \}.
\end{align*}
\end{proof}

The statement also holds if $\cM$ is a (weak or strong) right $T$-matroid.

Now we will introduce a lemma which shows a general version of the direction $(2) \Rightarrow (3)$ in Theorem~\ref{weakext.equi}.

\begin{lem} \label{ext.del}
Let $\cM$ be a (weak or strong) left $T$-matroid on $E$ and let $\sigma: \mathcal{C^*(M)} \rightarrow T$ be a localization of $\cM$. Then $\sigma$ defines a single-element extension of every deletion of $\cM$.
\end{lem}

\begin{proof} Let $\widetilde{\cM}$ be the extension of $\cM$ by an element $p$ determined by $\sigma$ and let $E_0 \subseteq E$. We will prove that $\sigma$ defines a single-element extension of $\cM\backslash E_0$ when $|E_0| = 1$ and the claim of the lemma follows by induction on the size of $E_0$. Let $E_0 = \{e_0\}$.

By Theorem~\ref{minor}, we know that
$$\cC^*((\tcM\backslash e_0)\backslash p) = \cC^*((\tcM\backslash p)\backslash e_0) = \cC^*(\cM\backslash e_0).$$
So $\tcM\backslash e_0$ is an extension of $\cM\backslash e_0$ by the element $p$.

Let $\sigma': \cC^*(\cM\backslash e_0) \rightarrow T$ be the function induced from $\sigma$ by deletion of $e_0$. We would like to show that $\sigma'$ is the localization of $\cM\backslash e_0$ corresponding to $\tcM\backslash e_0$.

Let $Z \in \cC^*(\cM\backslash e_0)$. Then by definition, there exists $Y\in \cC^*(\cM)$ such that $Z = Y\backslash e_0$. By Proposition~\ref{sigma.intro}, $\tY: = (Y, \sigma(Y)) \in \cC^*(\widetilde{\cM}).$ We claim that $\tY\backslash e_0 \in \cC^*(\tcM\backslash e_0)$. Otherwise, there exists $W\in \cC^*(\tcM)$ such that $\underline{W\backslash e_0} \subset \underline{\tY\backslash e_0}$. Then either $\underline{W\backslash \{e_0, p\}} \subset \underline{\tY\backslash \{e_0, p\}} = \underline{Z}$ or $\underline{W\backslash \{e_0, p\}} = \underline{\tY\backslash \{e_0, p\}}$ and $W(p) = 0$. The first case contradicts the definition of $Z$. 

Now we consider the second case. As $\underline{W\backslash \{e_0, p\}} = \underline{\tY\backslash \{e_0, p\}}$, then by comparability there exists $\alpha \in T^\times$ such that $W\backslash \{e_0, p\} = (\tY\backslash \{e_0, p\}) \cdot \alpha$. Without loss of generality, we assume that $\alpha = -1$. As $W(p) = 0$, then $W(e_0) \neq 0$ and $\tY(e_0) = 0$. So $\rank_\tcM (W^0\cap \tY^0) = \rank(\tcM) - 2$ and so $W$ and $\tY$ form a modular pair of $T$-cocircuits in $\tcM$. Let $e_1\in \underline{W}\cap \underline{\tY} = \underline{W\backslash \{e_0,p\}}$. By Modular Elimination, there exists $V\in \cC^*(\tcM)$ such that $V(e_1) = 0$ and $W(f) + \tY(f) - V(f) \in N_G$. Then $\underline{V} \subseteq (\underline{W}\cup \underline{\tY})\backslash e_1 = (\underline{Z}\cup \{e_0, p\})\backslash e_1$. So $\underline{V\backslash \{e_0, p\}} \subseteq \underline{Z}\backslash e_1 \subset \underline{Z}$, contradicting the definition of $Z$.

So our claim is true, that is $\tY\backslash e_0 \in \cC^*(\tcM\backslash e_0)$. So 
$$(Z, \sigma'(Z)) = (Y\backslash e_0, \sigma(Y)) = \tY\backslash e_0 \in \cC^*(\tcM\backslash e_0).$$

So $\sigma'$ is the localization of $\cM\backslash e_0$ corresponding to $\tcM\backslash e_0$, and so $\sigma$ defines a single-element extension of $\cM\backslash e_0$.
\end{proof}

The statement also holds if $\cM$ is a (weak or strong) right $T$-matroid.

Then we will introduce a lemma which shows a general version of the direction $(3) \Rightarrow (2)$ in Theorem~\ref{weakext.equi}.

\begin{lem}\label{ext.minor}
Let $\cM$ be a $T$-matroid on $E$ of rank 2 and let $\sigma: \mathcal{C^*(M)} \rightarrow T$ be a function. If $\sigma$ defines a single-element extension of every minor of $\cM$ on three elements, then $\sigma$ is a localization of $\cM$.
\end{lem}

\begin{proof} Let $\cl$ be the closure operator and $\cB$ be set of bases of $\underline{\cM}$. Let $[\cdot]: A_\cB \rightarrow T$ be left quasi-Pl\"ucker coordinates of $\cM$.
We define a collection $\tcB$ of subsets of $E\cup p$ and the left $T$-coordinates $[\cdot]_p: A_\tcB \rightarrow T$ as in Lemma~\ref{phi.loc.equi}. Then
$$\tcB: = \cB \cup \{\{e, p\} \,|\, e\in E, E\backslash \cl(\{e\}) = \underline{Y} \text{ and }\sigma(Y) \neq 0 \text{ for some } Y\in \cC^*(\cM)\}$$
and
\begin{equation}
[\{e,s\}, \{e,t\}]_p=\begin{cases}
[\{e,s\}, \{e,t\}]	& \text{ if } p\notin \{e,s,t\},\\
\overline{\sigma(Y)\cdot Y(t)^{-1}}  & \text{ if } p = s\text{ and } Y\in \mathcal{C^*(M)} \text{ with } \underline{Y}= E\backslash\cl(\{e\}),\\
\overline{ Y(s)\cdot \sigma(Y)^{-1}}  & \text{ if } p = t \text{ and } Y\in \mathcal{C^*(M)} \text{ with } \underline{Y}= E\backslash\cl(\{e\}),\\
[\{s,g\}, \{t,g\}] & \text{ if } p = e, \{s,t\}\notin \cB \text{ and } \{s,g\}, \{t,g\}\in \cB\\
&\,\,\text{ for some } g\in E, \\
-\overline{Y_t(s)\cdot \sigma(Y_t)^{-1}\cdot \sigma(Y_s)\cdot Y_s(t)^{-1}} & \text{ if } p = e, \{s,t\}\in \cB, Y_i\in \mathcal{C^*(M)} \text{ with } \\
&\,\,\, \underline{Y_i}= E\backslash\cl(\{i\}) \text{ for } i\in \{s,t\}.
\end{cases}
\label{ext.phi.form.rank2}
\end{equation}

Now we just need to show that $[\cdot]_p$ satisfies \ref{p3} and \ref{p5} in weak left quasi-Pl\"ucker coordinates axioms. So we need to show that
\begin{enumerate}[(a)]
\item for any $e, x_1, x_2, x_3 \in E\cup p$ with $\{e, x_i\}\in \tcB$ for all $i\in \{1,2,3\}$, we have
\begin{equation} \label{eq.rank2.1}
[\{e, x_1\}, \{e, x_2\}]_p \cdot [\{e, x_2\}, \{e, x_3\}]_p \cdot [\{e, x_3\}, \{e, x_1\}]_p = 1,
\end{equation}
\item for any $y_1, y_2, y_3, y_4 \in E\cup p$ with $\{y_i, y_j\} \in \tcB$ for all $i, j\in \{1,2,3,4\}$ and $i\neq j$, we have
\begin{equation}\label{eq.rank2.2}
 -1 + [\{y_2, y_4\}, \{y_1, y_2\}]_p\cdot [\{y_1, y_3\}, \{y_3, y_4\}]_p+ [\{y_1, y_4\}, \{y_1, y_2\}]_p \cdot [\{y_2, y_3\}, \{y_3, y_4\}]_p \in N_G.
\end{equation}
\end{enumerate}

First, we prove Equation~\eqref{eq.rank2.1}. If $p\notin \{e, x_1, x_2, x_3\}$, then from Formula~\eqref{ext.phi.form.rank2} we know that $[\cdot]$ and $[\cdot]_p$ coincide on each pair of the bases in \eqref{eq.rank2.1}. So Equation~\eqref{eq.rank2.1} holds because $[\cdot]$ satisfies \ref{p3}. If $p\in \{e, x_1, x_2, x_3\}$, then we consider a minor $\cM_3$ of $\cM$ on $\{e, x_1, x_2, x_3\}\backslash p$. By assumption, $\sigma$ defines a single-element extension of $\cM_3$. So Equation~\eqref{eq.rank2.1} holds.

Now, we prove the inclusion~\eqref{eq.rank2.2}. If $p\notin \{y_1, y_2, y_3, y_4\}$, then from Formula~\eqref{ext.phi.form.rank2} we know that $[\cdot]$ and $[\cdot]_p$ coincide on each pair of the bases in \eqref{eq.rank2.2}. So the inclusion~\eqref{eq.rank2.2} holds because $[\cdot]$ satisfies \ref{p5}. If $p \in \{y_1, y_2, y_3, y_4\}$. Now we consider a minor $\cM_3'$ of $\cM$ on $\{y_1, y_2, y_3, y_4\}\backslash p$. By assumption, $\sigma$ defines a single-element extension of $\cM_3'$. So the inclusion~\eqref{eq.rank2.2} holds.
\end{proof}


Finally, we will prove the main theorem (Theorem~\ref{weakext.equi}) in this section.

\begin{proof}[Proof of Theorem~\ref{weakext.equi}] (1) $\Rightarrow$ (2): This is proved by Lemma~\ref{ext.contr}.

(2) $\Rightarrow$ (3): This is proved by Lemma~\ref{ext.del}.

(3) $\Rightarrow$ (2): This is proved by Lemma~\ref{ext.minor}.

(2) $\Rightarrow$ (1): Let $d = \rank(\cM)$. Let $r$ be the rank function, $\cl$ be the closure operator and $\cB$ be the set of bases of of $\underline{\cM}$. Let $[\cdot]: A_\cB \rightarrow T$ be weak left quasi-Pl\"ucker coordinates of $\cM$.

We define a collection $\tcB$ of subsets of $E\cup p$ by 
$$\tcB = \cB \cup \{Fp \subseteq E\cup p \,|\, |F| = d-1, E \backslash \cl(F) = \underline{Y} \text{ and } \sigma(Y) \neq 0 \text{ for some } Y \in \cC^*(\cM) \},$$
as in Lemma~\ref{phi.loc.equi}. Now we consider the left $T$-coordinates $[\cdot]_p: A_{\tcB} \rightarrow T$ defined in Lemma~\ref{phi.loc.equi}:
\begin{equation}
[Fs, Ft]_p=\begin{cases}
[Fs, Ft] & \text{ if } p\notin Fst,\\
\overline{\sigma(Y)\cdot Y(t)^{-1}}  & \text{ if } p = s\text{ and } Y\in \mathcal{C^*(M)} \text{ with } \underline{Y}= E\backslash\cl(F),\\
\overline{ Y(s)\cdot \sigma(Y)^{-1}}  & \text{ if } p = t \text{ and } Y\in \mathcal{C^*(M)} \text{ with } \underline{Y}= E\backslash\cl(F),\\
[Gsg, Gtg] & \text{ if } p\in F, G=F\backslash p, Gst\notin \cB \text{ and } Gsg, Gtg\in \cB \\
&\,\,\, \text{ for some } g\in E, \\
-\overline{Y_t(s)\cdot \sigma(Y_t)^{-1}\cdot\sigma(Y_s)\cdot Y_s(t)^{-1}} & \text{ if } p\in F, G=F\backslash p, Gst\in \cB, Y_i\in \mathcal{C^*(M)} \text{ with } \\
& \,\,\, \, \, \underline{Y_i}= E\backslash\cl(Gi) \text{ for } i\in \{s,t\}.
\end{cases}
\label{ext.phi.form.thm}
\end{equation}

By Lemma~\ref{phi.loc.equi}, we just need to show that $[\cdot]_p$ satisfies \ref{p3} and \ref{p5} in weak left quasi-Pl\"ucker coordinates axioms. So we need to show 
\begin{enumerate} [(A)]
\item for any $x_1, x_2, x_3 \in E\cup p$ and $H\subseteq E\cup p$ with $Hx_1$, $Hx_2$, $Hx_3\in \tcB$, we have
\begin{equation} \label{mainfirsteq}
[Hx_1, Hx_2]_p \cdot [Hx_2, Hx_3]_p \cdot [Hx_3, Hx_1]_p = 1,
\end{equation}

\item for any $y_1, y_2, y_3, y_4 \in E\cup p$ and $F\subseteq E\cup p$ with $Fy_iy_j\in \tcB$ for all $i, j \in \{1, 2, 3, 4\}$ and $i\neq j$, we have
\begin{equation}\label{maineq}
 -1 + [Fy_2y_4, Fy_1y_2]_p\cdot [Fy_1y_3, Fy_3y_4]_p+ [Fy_1y_4, Fy_1y_2]_p \cdot [Fy_2y_3, Fy_3y_4]_p \in N_G.
\end{equation}
\end{enumerate}

First, we will prove Equation~\eqref{mainfirsteq}. 

If $p\notin H\cup \{x_1, x_2, x_3\}$, then we know that $[\cdot]$ and $[\cdot]_p$ coincide on each pair of the bases in \eqref{mainfirsteq}. So Equation~\eqref{mainfirsteq} holds because $[\cdot]$ satisfies \ref{p3}. 

If $p\in \{x_1, x_2, x_3\}$, without loss of generality we assume $p= x_1$. Then there exists $Y\in \cC^*(\cM)$ with $\underline{Y} = E\backslash \cl(H)$ and $\sigma(Y)\neq 0$. Then $ x_2, x_3 \in \underline{Y}$. By Formula~\eqref{ext.phi.form.thm}, we get that
\begin{align*}
[Hx_1, Hx_2]_p \cdot [Hx_2, Hx_3]_p \cdot [Hx_3, Hx_1]_p & = [Hp, Hx_2]_p \cdot [Hx_2, Hx_3]_p \cdot [Hx_3, Hp]_p \\
& = \overline{\sigma(Y)Y(x_2)^{-1}} \cdot \overline{Y(x_2)Y(x_3)^{-1}} \cdot \overline{Y(x_3)\sigma(Y)^{-1}} \\
& = 1.
\end{align*}

If $p\in H$, then let $G = H\backslash \{p\}$. As $r(G) = d-2$, then let $\cM' = \cM/G$ and $[\cdot]'=[\cdot]/G$. Then $\cM'$ is a rank 2 contraction of $\cM$ and $[\cdot]'$ are weak left quasi-Pl\"ucker coordinates of $\cM'$. Let $\sigma': \cC^*(\cM') \rightarrow T$ be the function induced from $\sigma$. As $\sigma$ defines a single-element extension of every rank 2 contraction of $\cM$, then $\sigma'$ is a localization of $\cM'$ and we call the corresponding extension $\tcM'$. So by Lemma~\ref{ext.basis}, the set of bases of $\underline{\tcM'}$ is $$\tcB' = \cB(\underline{\cM'}) \cup \{\{e, p\} \subseteq (E\backslash G)\cup p \,|\, (E \backslash G) \backslash \cl_{\underline{\cM'}}(\{e\}) = \underline{X} \text{ and } \sigma'(X) \neq 0 \text{ for some } X \in \cC^*(\cM') \}.$$
By Proposition~\ref{ext.uniq}, the corresponding weak left quasi-Pl\"ucker coordinates $[\cdot]'_p: A_{\tcB'} \rightarrow T$ are given by
\begin{align*}
& \,\, [\{e,s\}, \{e,t\}]'_p \\
= & \begin{cases}
[\{e,s\}, \{e,t\}]	' & \text{ if } p\notin \{e,s,t\},\\
\overline{\sigma'(X)\cdot X(t)^{-1}}  & \text{ if } p = s\text{ and } X\in \mathcal{C^*(M')} \text{ with } \underline{X}= (E\backslash G)\backslash \cl_{\underline{\cM'}}(\{e\}),\\
\overline{ X(s)\cdot \sigma'(X)^{-1}}  & \text{ if } p = t \text{ and } X\in \mathcal{C^*(M')} \text{ with } \underline{X}= (E\backslash G)\backslash \cl_{\underline{\cM'}}(\{e\}),\\
[\{s,g\}, \{t,g\}]' & \text{ if } p = e, \{s,t\}\notin \cB(\underline{\cM'}) \text{ and } \{s,g\}, \{t,g\}\in \cB(\underline{\cM'})\\
&\,\,\text{ for some } g\in E\backslash G, \\
-\overline{X_t(s)\cdot \sigma'(X_t)^{-1}\cdot\sigma'(X_s)\cdot X_s(t)^{-1}} & \text{ if } p = e, \{s,t\}\in \cB(\underline{\cM'}), X_i\in \mathcal{C^*(M')} \text{ with } \\
&\,\,\, \underline{X_i}= (E\backslash G)\backslash \cl_{\underline{\cM'}}(\{i\}) \text{ for } i\in \{s,t\}.
\end{cases}
\end{align*}

Let $k, j \in \{x_1, x_2, x_3\}$ with $k\neq j$. If $Gkj \in \cB$, then there exists $Y_k, Y_j\in \cC^*(\cM)$ with $\underline{Y_k} = E\backslash \cl(Gk)$, $\sigma(Y_k)\neq 0$, $\underline{Y_j} = E\backslash \cl(Gj)$ and $\sigma(Y_j) \neq 0$. Let $X_k = Y_k\backslash G$ and $X_j = Y_j\backslash G$. Then $X_k, X_j\in \cM'$. So 
\begin{align*}
[Hk, Hj]_p & = [Gpk, Gpj]_p \\
& = - \overline{Y_j(k)\cdot \sigma(Y_j)^{-1}\cdot\sigma(Y_k)\cdot Y_k(j)^{-1}} \\
& = - \overline{X_j(k)\cdot \sigma'(X_j)^{-1}\cdot\sigma'(X_k)\cdot X_k(j)^{-1}} \\
& = [\{p,k\}, \{p,j\}]'_p.
\end{align*}
If $Gkj \notin \cB$, then let $g\in E$ such that $Gkg, Gjg\in \cB$. So
$$[Hk, Hj]_p = [Gpk, Gpj]_p = [Gkg, Gjg] = [\{k,g\}, \{j, g\}]' = [\{k,p\}, \{j, p\}]'_p.$$

So Equation~\eqref{mainfirsteq} holds because $[\cdot]'_p$ satisfies \ref{p3}. 

Now we would like to prove the inclusion~\eqref{maineq}:
$$-1 + [Fy_2y_4, Fy_1y_2]_p\cdot [Fy_1y_3, Fy_3y_4]_p+ [Fy_1y_4, Fy_1y_2]_p \cdot [Fy_2y_3, Fy_3y_4]_p \in N_G.$$
By Formula~\eqref{ext.phi.form.thm} and our argument above, we know that $[\cdot]_p$ satisfies \ref{p1}, \ref{p2}, \ref{p3} and \ref{p4}. So we can multiply both sides of \eqref{maineq} on the left by $[Fy_1y_2, Fy_1y_4]_p$ and on the right by $[Fy_3y_4, Fy_2y_3]_p$, and get
$$-[Fy_1y_2, Fy_1y_4]_p \cdot [Fy_3y_4, Fy_2y_3]_p - [Fy_2y_4, Fy_1y_4]_p \cdot [Fy_1y_3, Fy_2y_3]_p + 1 \in N_G.$$
Then from this inclusion and \eqref{maineq}, we know that switching $y_1$ and $y_2$ or switching $y_1$ and $y_4$ in \eqref{maineq} will result in an inclusion equivalent to the inclusion~\eqref{maineq}.

Now we begin our proof and will divide it into three cases.

\textbf{Case 1:} $p\notin F\cup \{y_1, y_2, y_3, y_4\}$.

From Formula~\eqref{ext.phi.form.thm}, we know that $[\cdot]$ and $[\cdot]_p$ coincide on each pair of the bases in \eqref{maineq}. So the inclusion~\eqref{maineq} holds because $[\cdot]$ satisfies \ref{p5}.

\textbf{Case 2:} $p \in \{y_1, y_2, y_3, y_4\}$. Then there exists $Y_i \in \cC^*(\cM)$ such that $\underline{Y_i} = E\backslash \cl(Fi)$ and $\sigma(Y_i) \neq 0$ for $i\in \{y_1, y_2, y_3, y_4\} \backslash \{p\}$. Moreover, $\{y_1, y_2, y_3, y_4\} \backslash \{p, i\} \subseteq \underline{Y_i}$ for $i\in \{y_1, y_2, y_3, y_4\} \backslash \{p\}$. We denote $(Y_i, \sigma(Y_i))$ by $\tY_i$.
\begin{enumerate}[\emph{Case~{2.\arabic*}:}]
\item If $p\in \{ y_1, y_2, y_4\}$, then without loss of generality, we may assume $p = y_1$. Now we consider $Y_2$ and $Y_3$. As $r(F) = d-2$, then $Y_2$ and $Y_3$ form a modular pair in $\cC^*(\cM)$. So $(\tY_2, \tY_3)$ is a modular pair in $\cC^*(\tcM)$. By symmetry in weak circuit axioms, we may assume that $\sigma(Y_2) = - \sigma(Y_3)$. Let $Z = \Mod(Y_2, Y_3, p)$. Then $Z(y_2) = Y_3(y_2)$ and $Y_{2}(y_4) +Y_{3} (y_4)-Z(y_4) \in N_G$. So by Formula~\eqref{ext.phi.form.thm} and Dual Pivoting property,
\begin{align*}
&  -1 + [Fy_2y_4, Fy_1y_2]_p\cdot [Fy_1y_3, Fy_3y_4]_p+ [Fy_1y_4, Fy_1y_2]_p \cdot [Fy_2y_3, Fy_3y_4]_p \\
= &  -1 + [Fy_2y_4, Fy_2p]_p\cdot [Fy_3p, Fy_3y_4]_p+ [Fy_4p, Fy_2p]_p \cdot [Fy_2y_3, Fy_3y_4]_p \\
= &  -1 + \overline{Y_{2} (y_4)\sigma(Y_{2})^{-1}} \cdot \overline{\sigma(Y_{3}) Y_{3}(y_4)^{-1}} + \overline{Z(y_4)Z(y_2)^{-1}} \cdot \overline{Y_{3}(y_2)Y_{3}(y_4)^{-1}}\\
= &  -1 - \overline{Y_{2} (y_4) Y_{3}(y_4)^{-1}} + \overline{Z(y_4)Y_{3}(y_4)^{-1}}\\
= &  - (\overline{Y_{3}(y_4) +Y_{2} (y_4)-Z(y_4)}) \cdot \overline{Y_{3}(y_4)^{-1}}\\
\in & N_G.
\end{align*}

\item If $p = y_3$, then we consider $Y_1$ and $Y_2$. Similarly, $(\tY_1, \tY_2)$ is a modular pair in $\cC^*(\tcM)$. By symmetry in weak circuit axioms, we may assume that $\sigma(Y_1) = - \sigma(Y_2)$. Let $Z = \Mod(Y_1, Y_2, p)$. Then $Z(y_1) = Y_2(y_1)$, $Z(y_2) = Y_1(y_2)$ and $Y_{1}(y_4) +Y_{2} (y_4) - Z(y_4) \in N_G$. So by Formula~\eqref{ext.phi.form.thm} and Dual Pivoting property,
\begin{align*}
&  -1 + [Fy_2y_4, Fy_1y_2]_p\cdot [Fy_1y_3, Fy_3y_4]_p+ [Fy_1y_4, Fy_1y_2]_p \cdot [Fy_2y_3, Fy_3y_4]_p \\
= & -1 + [Fy_2y_4, Fy_1y_2]_p\cdot [Fy_1p, Fy_4p]_p+ [Fy_1y_4, Fy_1y_2]_p \cdot [Fy_2p, Fy_4p]_p \\
= & -1 + \overline{Y_{2} (y_4)Y_{2} (y_1)^{-1}} \cdot \overline{Z(y_1) Z(y_4)^{-1}} + \overline{Y_{1}(y_4)Y_{1}(y_2)^{-1}} \cdot \overline{Z(y_2)Z(y_4)^{-1}}\\
= &  -1 + \overline{Y_{2} (y_4) Z(y_4)^{-1}} + \overline{Y_{1}(y_4) Z(y_4)^{-1}}\\
= & (\overline{-Z(y_4) + Y_{2}(y_4) +Y_{1} (y_4)}) \cdot \overline{Z(y_4)^{-1}} \\
\in & N_G.
\end{align*}

\end{enumerate}

\textbf{Case 3:} $p \in F$. Let $G = F \backslash \{p\}$. So $\rank_{\tcM}(G) =\rank_{\cM}(G) = d-3$.

As $r(Gy_iy_j) = d-1$ for any two distinct elements $i,j\in \{1,2,3,4\}$, then $r(Gy_iy_jy_k)\geq d-1$ for any three distinct elements $i,j,k\in \{1,2,3,4\}$. There is an easy case and two harder cases to check:

\begin{enumerate}[\emph{Case~{3.\arabic*}:}]
\item $r(Gy_1y_2y_3y_4) = d-1$. Then $r(Gy_iy_jy_k) = d-1$ for any three distinct elements $i,j,k\in \{1,2,3,4\}$. So $\cl(Gy_1y_2) = \cl(Gy_1y_3) = \cl(Gy_1y_4) = \cl(Gy_2y_3) = \cl(Gy_2y_4) = \cl(Gy_3y_4)$. Let $e\in E$ such that $Gy_1y_2e \in \cB$. Then $Gy_1y_3e, Gy_1y_4e, Gy_2y_3e, Gy_2y_4e, Gy_3y_4 e\in \cB$. As $[\cdot]$ satisfies \ref{p5}, we have
\begin{align*}
& -1 + [Fy_2y_4, Fy_1y_2]_p\cdot [Fy_1y_3, Fy_3y_4]_p+ [Fy_1y_4, Fy_1y_2]_p \cdot [Fy_2y_3, Fy_3y_4]_p\\
= & -1 + [Gy_2y_4p, Gy_1y_2p]_p\cdot [Gy_1y_3p, Gy_3y_4p]_p+ [Gy_1y_4p, Gy_1y_2p]_p \cdot [Gy_2y_3p, Gy_3y_4p]_p\\
= & -1 + [Gy_2y_4e, Gy_1y_2e]\cdot [Gy_1y_3e, Gy_3y_4e]+ [Gy_1y_4e, Gy_1y_2e]\cdot [Gy_2y_3e, Gy_3y_4e]\\
\in & N_G.
\end{align*}

\item $r(Gy_1y_2y_3y_4) = d$, and there exists $f\in \{y_1, y_2, y_3, y_4\}$ such that $r(Gy_1y_2y_3y_4 \backslash \{f\}) = d-1$.

\begin{enumerate}[\emph{Case~{3.2.\arabic*}:}]
\item If $f \in \{y_1, y_2, y_4\}$, then without loss of generality we assume $f = y_2$ and so $r(Gy_1y_3y_4) = d-1$. So $\cl(Gy_1y_3) = \cl(Gy_1y_4) = \cl(Gy_3y_4) = \cl(Gy_1y_3y_4)$ and $Gy_1y_2y_4, Gy_1y_2y_3, Gy_2y_3y_4 \in \cB$. So there exists $Y_{134}\in \cC^*(\cM)$ with $\underline{Y_{134}} = E\backslash \cl(Gy_1y_3y_4) \ni y_2$ and $\sigma(Y_{134}) \neq 0$.
Then we have 
$$[Gy_1y_4p, Gy_1y_2y_4]_p = \overline{\sigma(Y_{134}) \cdot Y_{134}(y_2)^{-1}} = [Gy_3y_4p, Gy_2y_3y_4]_p.$$
As $[\cdot]_p$ satisfies \ref{p3}, then
\begin{align*}
 [Fy_1y_4, Fy_1y_2]_p = [Gpy_1y_4, Gpy_1y_2]_p = - [Gy_1y_2y_4, Gy_1y_2 p]_p\cdot  [Gy_1y_4p, Gy_1y_2y_4]_p, \\
 [Fy_2y_3, Fy_3y_4]_p = [Gpy_2y_3, Gpy_3y_4]_p = - [Gy_2y_3y_4, Gy_3y_4 p]_p\cdot [Gy_2y_3p, Gy_2y_3y_4]_p.
\end{align*}
So \begin{align*}
& [Fy_1y_4, Fy_1y_2]_p \cdot [Fy_2y_3, Fy_3y_4]_p  \\
 = & [Gy_1y_2y_4, Gy_1y_2 p]_p\cdot [Gy_1y_4p, Gy_1y_2y_4]_p \cdot [Gy_2y_3y_4, Gy_3y_4 p]_p\cdot [Gy_2y_3p, Gy_2y_3y_4]_p \\
 = & [Gy_1y_2y_4, Gy_1y_2 p]_p\cdot [Gy_2y_3p, Gy_2y_3y_4]_p.
\end{align*}
As $[\cdot]_p$ satisfies \ref{p4}, then
$$ [Fy_1y_3, Fy_3y_4]_p = [Gy_1y_3y_2, Gy_3y_4y_2].$$
So to prove the inclusion~\eqref{maineq}, we only need to verify,
$$-1+ [Gy_2y_4p, Gy_1y_2p]_p \cdot [Gy_1y_3y_2, Gy_3y_4y_2] + [Gy_1y_2y_4, Gy_1y_2 p]_p\cdot [Gy_2y_3p, Gy_2y_3y_4]_p \in N_G.$$

Let $J = Gy_2$. Let $\cM_2 = \cM/J$ and $[\cdot]_2 = [\cdot]/J$. $\cM_2$ is a rank 2 contraction of $\cM$ and $[\cdot]_2$ are weak left quasi-P\"ucker coordinates of $\cM_2$. Let $\sigma_2 : \cC^*(\cM_2) \rightarrow T$ be the function induced from $\sigma$. As $\sigma$ defines a single-element extension of every rank 2 contraction of $\cM$, then $\sigma_2$ is a localization of $\cM_2$ and we call the corresponding extension $\tcM_2$. So by Lemma~\ref{ext.basis}, the set of bases of $\underline{\tcM_2}$ is
$$\tcB_2 = \cB(\underline{\tcM_2}) \cup \{\{e, p\} \subseteq (E\backslash J)\cup p \,|\, (E \backslash J) \backslash \cl_{\underline{\cM_2}}(\{e\}) = \underline{X} \text{ and } \sigma_2(X) \neq 0 \text{ for some } X \in \cC^*(\cM_2) \}.$$
By Proposition~\ref{ext.uniq}, the corresponding weak left quasi-Pl\"ucker coordinates $[\cdot]_p^2: A_{\tcB_2} \rightarrow T$ are given by
\begin{align*}
& \,\, [\{e,s\}, \{e,t\}]_p^2\\
= & \begin{cases}
[\{e,s\}, \{e,t\}]_2 & \text{ if } p\notin \{e,s,t\},\\
\overline{\sigma_2(X)\cdot X(t)^{-1}}  & \text{ if } p = s\text{ and } X\in \cC^*(\cM_2) \text{ with } \underline{X}= (E\backslash J)\backslash \cl_{\underline{\cM_2}}(\{e\}),\\
\overline{ X(s)\cdot \sigma_2(X)^{-1}}  & \text{ if } p = t \text{ and } X\in \cC^*(\cM_2) \text{ with } \underline{X}= (E\backslash J)\backslash \cl_{\underline{\cM_2}}(\{e\}),\\
[\{s,g\}, \{t,g\}]_2 & \text{ if } p = e, \{s,t\}\notin \cB(\underline{\tcM_2}) \text{ and } \{s,g\}, \{t,g\}\in \cB(\underline{\tcM_2})\\
&\,\,\text{ for some } g\in E\backslash J, \\
-\overline{X_t(s)\cdot \sigma_2(X_t)^{-1}\cdot\sigma_2(X_s)\cdot X_s(t)^{-1}} & \text{ if } p = e, \{s,t\}\in \cB(\underline{\tcM_2}), X_i\in \cC^*(\cM_2) \text{ with } \\
&\,\,\, \underline{X_i}= (E\backslash J)\backslash \cl_{\underline{\cM_2}}(\{i\}) \text{ for } i\in \{s,t\}.
\end{cases}
\end{align*}

As $Gy_1y_2y_4, Gy_2y_3y_4 \in \cB$, then there exists $Y_{2i}\in \cC^*(\cM)$ such that 
$\underline{Y_{2i}} = E \backslash \cl(Giy_2) = E \backslash \cl(Ji)$, $\sigma(Y_{2i})\neq 0$ for $i\in \{1,3,4\}$. Let $X_{2i} = Y_{2i}\backslash J$ for $i\in \{1,3,4\}$. Then $X_{2i} \in \cC^*(\cM_2)$ for $i\in \{1,3,4\}$. So 
\begin{align*}
[Gy_2y_4p, Gy_1y_2p]_p & = -\overline{Y_{12}(y_4) \sigma(Y_{12})^{-1}\sigma(Y_{24}) Y_{24}(y_1)^{-1}} \\
& = -\overline{X_{12}(y_4) \sigma_2(X_{12})^{-1}\sigma_2(X_{24}) X_{24}(y_1)^{-1}} \\
& = [\{p,y_4\}, \{p,y_1\}]_p^2,\\
[Gy_1y_3y_2, Gy_3y_4y_2] & = [\{y_1,y_3\}, \{y_3,y_4\}]_p^2,\\
[Gy_1y_2y_4, Gy_1y_2 p]_p\cdot [Gy_2y_3p, Gy_2y_3y_4]_p &  = \overline{Y_{12}(y_4)\sigma(Y_{12})^{-1}} \cdot \overline{\sigma(Y_{23})Y_{23}(y_4)^{-1}}\\
&  = \overline{X_{12}(y_4)\sigma_2(X_{12})^{-1}} \cdot \overline{\sigma_2(X_{23})X_{23}(y_4)^{-1}}\\
& = [\{y_1, y_4\}, \{y_1, p\}]_p^2 \cdot [\{y_3, p\}, \{y_3,y_4\}]_p^2.
\end{align*}
So to prove the inclusion~\eqref{maineq}, we only need to verify,
$$-1 + [\{p,y_4\}, \{p,y_1\}]_p^2\cdot [\{y_1,y_3\}, \{y_3,y_4\}]_p^2 + [\{y_1, y_4\}, \{y_1, p\}]_p^2 \cdot [\{y_3, p\}, \{y_3,y_4\}]_p^2\in N_G.$$
As $\sigma_2$ is a localization of $\cM_2$, then $[\cdot]_p^2$ satisfies \ref{p5}.

So the inclusion~\eqref{maineq} holds.

\item If $f = y_3$, then $r(Gy_1y_2y_4) = d-1$. So 
$\cl(Gy_1y_2) = \cl(Gy_1y_4) = \cl(Gy_2y_4) = \cl(Gy_1y_2y_4)$ and $Gy_1y_3y_4, Gy_1y_2y_3, Gy_2y_3y_4 \in \cB$. So there exists $Y_{124}\in \cC^*(\cM)$ with $\underline{Y_{124}} = E\backslash \cl(Gy_1y_2y_4) \ni y_3$ and $\sigma(Y_{124}) \neq 0$. Then we have 
\begin{align*}
[Fy_2y_4, Fy_1y_2]_p & = [Gpy_2y_4, Gpy_1y_2]_p = [Gy_2y_4y_3, Gy_1y_2 y_3\}], \\
[Fy_1y_4, Fy_1y_2]_p & = [Gpy_1y_4, Gpy_1y_2]_p = [Gy_1y_4y_3, Gy_1y_2 y_3].
\end{align*}
So to prove the inclusion~\eqref{maineq}, we only need to verify,
$$-1 + [Gy_2y_4y_3, Gy_1y_2 y_3\}] \cdot [Gpy_1y_3, Gpy_3y_4]_p + [Gy_1y_4y_3, Gy_1y_2 y_3]\cdot [Gpy_2y_3, Gpy_3y_4]_p.$$

Let $K = Gy_3$. Let $\cM_2' = \cM/K$ and $[\cdot]'_2 = [\cdot]/K$. $\cM_2'$ is a rank 2 contraction of $\cM$ and $[\cdot]_2'$ are weak left quasi-P\"ucker coordinates of $\cM_2'$. Let $\sigma_2' : \cC^*(\cM_2') \rightarrow T$ be the function induced from $\sigma$. As $\sigma$ defines a single-element extension of every rank 2 contraction of $\cM$, then $\sigma_2'$ is a localization of $\cM_2'$ and we call the corresponding extension $\tcM_2'$. So by Lemma~\ref{ext.basis}, the set of bases of $\underline{\tcM_2'}$ is
$$\tcB_2' = \cB(\underline{\tcM_2'}) \cup \{\{e, p\} \subseteq (E\backslash K)\cup p \,|\, (E \backslash K) \backslash \cl_{\underline{\cM_2'}}(\{e\}) = \underline{X} \text{ and } \sigma_2'(X) \neq 0 \text{ for some } X \in \cC^*(\cM_2') \}.$$
By Proposition~\ref{ext.uniq}, the corresponding weak left quasi-Pl\"ucker coordinates $[\cdot]_p^2: A_{\tcB_2'} \rightarrow T$ are given by
\begin{align*}
& \,\, [\{e,s\}, \{e,t\}]_p^{2'} \\
= & \begin{cases}
[\{e,s\}, \{e,t\}]_{2}' & \text{ if } p\notin \{e,s,t\},\\
\overline{\sigma_2'(X)\cdot X(t)^{-1}}  & \text{ if } p = s\text{ and } X\in \cC^*(\cM_2') \text{ with } \underline{X}= (E\backslash K)\backslash \cl_{\underline{\cM_2'}}(\{e\}),\\
\overline{ X(s)\cdot \sigma_2'(X)^{-1}}  & \text{ if } p = t \text{ and } X\in \cC^*(\cM_2') \text{ with } \underline{X}= (E\backslash K)\backslash \cl_{\underline{\cM_2'}}(\{e\}),\\
[\{s,g\}, \{t,g\}]_{2}' & \text{ if } p = e, \{s,t\}\notin \cB(\underline{\tcM_2'}) \text{ and } \{s,g\}, \{t,g\}\in \cB(\underline{\tcM_2'})\\
&\,\,\text{ for some } g\in E\backslash K, \\
-\overline{X_t(s)\cdot \sigma_2'(X_t)^{-1}\cdot\sigma_2'(X_s)\cdot X_s(t)^{-1}} & \text{ if } p = e, \{s,t\}\in \cB(\underline{\tcM_2'}), X_i\in \cC^*(\cM_2') \text{ with } \\
&\,\,\, \underline{X_i}= (E\backslash K)\backslash \cl_{\underline{\cM_2'}}(\{i\}) \text{ for } i\in \{s,t\}.
\end{cases}
\end{align*}

As $Gy_1y_3y_4, Gy_2y_3y_4 \in \cB$, then there exists $Y_{3i}\in \cC^*(\cM)$ such that 
$\underline{Y_{3i}} = E \backslash \cl(Giy_3) = E \backslash \cl(Ki)$, $\sigma(Y_{3i})\neq 0$ for $i\in \{1,2,4\}$. Let $X_{3i} = Y_{3i}\backslash K$ for $i\in \{1,2,4\}$. Then $X_{3i} \in \cC^*(\cM_2')$ for $i\in \{1,2,4\}$.
So
\begin{align*}
[Gy_2 y_4y_3, Gy_1y_2y_3] & = [\{y_2, y_4\}, \{y_1,y_2\}]_p^{2'}, \\
[Gpy_1y_3, Gpy_3y_4]_p & = -\overline{Y_{34}(y_1) \sigma(Y_{34})^{-1}\sigma(Y_{13}) Y_{13}(y_4)^{-1}} \\
& = -\overline{X_{34}(y_1) \sigma_2'(X_{34})^{-1}\sigma_2'(X_{13}) X_{13}(y_4)^{-1}} \\
& = [\{y_1, y_3\}, \{y_3,y_4\}]_p^{2'}, \\
[Gy_1y_4y_3, Gy_1y_2y_3] & = [\{y_1, y_4\}, \{y_1,y_2\}]_p^{2'}, \\
[Gpy_2y_3, Gpy_3y_4]_p & = -\overline{Y_{34}(y_2) \sigma(Y_{34})^{-1}\sigma(Y_{23}) Y_{23}(y_4)^{-1}} \\
& = -\overline{X_{34}(y_2) \sigma_2'(X_{34})^{-1}\sigma_2'(X_{23}) X_{23}(y_4)^{-1}} \\
 & = [\{y_2, y_3\}, \{y_3,y_4\}]_p^{2'}.
\end{align*}
So to prove the inclusion~\eqref{maineq}, we only need to verify,
$$-1 + [\{y_2, y_4\}, \{y_1,y_2\}]_p^{2'} \cdot [\{y_1, y_3\}, \{y_3,y_4\}]_p^{2'} + [\{y_1, y_4\}, \{y_1,y_2\}]_p^{2'}\cdot [\{y_2, y_3\}, \{y_3,y_4\}]_p^{2'}.$$

As $\sigma_2'$ is a localization of $\cM_2'$, then $[\cdot]_p^{2'}$ satisfies \ref{p5}.

So the inclusion~\eqref{maineq} holds.
\end{enumerate}

\item For any three distinct elements $i,j,k \in \{1, 2, 3, 4\}$, $r(Gy_iy_jy_k) = d$. Then for each partition $\{i, j\}\cupdot \{k, l\}$ of $\{1, 2, 3, 4\}$, there exists $Y_{ij} \in \cC^*(\cM)$ such that $\underline{Y_{ij}} = E\backslash \cl(Gy_iy_j)$, $\sigma(Y_{ij})\neq 0$ and $y_k, y_l \in \underline{Y_{ij}}$. By symmetry in weak circuit axioms, we may assume $Y_{23}(y_1) = Y_{13}(y_2) = Y_{12}(y_3) = Y_{24}(y_1) = -Y_{34}(y_1) = 1$ and $Y_{14}(y_2) = Y_{13}(y_4)^{-1} Y_{23}(y_4)$. By Lemma~\ref{loc.eq}, we can choose $\sigma(Y_{23}) = 1$. We list values of $Y_{23}$, $Y_{13}$, $Y_{12}$, $Y_{34}$, $Y_{24}$ and $Y_{14}$ on elements $\{y_1, y_2, y_3, y_4\}$ and the corresponding $\sigma$ values in Table~\ref{table1}.

\begin{table}[htb]
\centering
\begin{tabular}{ |c | c |c|c|c|c| }\hline
& $y_1$ & $y_2$ & $y_3$ & $y_4$ & $\sigma$ \\ \hline 
 $Y_{23}$ & 1 & 0 & 0  & $Y_{23}(y_4)$ & 1\\ \hline
 $Y_{13}$ & 0 & 1 & 0 & $Y_{13}(y_4)$ & $\sigma(Y_{13})$ \\ \hline
 $Y_{12}$ & 0 & 0 & 1& $Y_{12}(y_4)$ & $\sigma(Y_{12})$ \\ \hline
 $Y_{24}$ & 1 & 0 & $Y_{24}(y_3)$ & 0 & $\sigma(Y_{24})$  \\ \hline
 $Y_{34}$ & -1 & $Y_{34}(y_2)$ & 0 & 0 & $\sigma(Y_{34})$ \\ \hline
 $Y_{14}$ & 0 & $Y_{13}(y_4)^{-1} Y_{23}(y_4)$ & $Y_{14}(y_3)$ &0& $\sigma(Y_{14})$ \\ \hline
\end{tabular}
\smallskip
\caption{}
\label{table1}
\end{table}

We now find relationships between some numbers in this table.

Let $\alpha\in T^\times$ with $\alpha = Y_{12}(y_4)^{-1} Y_{23}(y_4).$ We consider the rank 2 contraction $\cM/Gy_2$. Let $\sigma_{y_2}: \cC^*(\cM/Gy_2) \rightarrow T$ be the function induced from $\sigma$. As $\sigma$ is a localization on the rank 2 contraction $\cM/Gy_2$, then $\sigma_{y_2}$ is a localization of $\cM/Gy_2$. Now we consider the $T$-cocircuits $Y_{23}\backslash Gy_2, Y_{12}\backslash Gy_2 \cdot (-\alpha)$ and $Y_{24}\backslash Gy_2$ of $\cM/Gy_2$. We list their values on elements $\{y_1, y_3, y_4\}$ and the corresponding $\sigma_{y_2}$ values in Table~\ref{table2},

\begin{table}[htb]
\centering
\begin{tabular}{ |c| c |c|c|c|c| }\hline
& $y_1$  & $y_3$ & $y_4$ & $\sigma_{y_2}$ \\ \hline 
 $Y_{23}\backslash G y_2$ & 1 & 0  & $Y_{23}(y_4)$ & 1 \\ \hline
 $Y_{12}\backslash G y_2 \cdot (-\alpha)$ & 0  & $-\alpha$ & $-Y_{23}(y_4)$ & $\sigma(Y_{12}) \cdot (-\alpha)$  \\ \hline
 $Y_{24}\backslash Gy_2$ & 1  & $Y_{24}(y_3)$ & 0 &  $\sigma(Y_{24})$ \\ \hline
\end{tabular}
\smallskip
\caption{}
\label{table2}
\end{table}

Note that $Y_{24}\backslash Gy_2$ eliminates $y_4$ between $Y_{23}\backslash G y_2$ and $ Y_{12}\backslash Gy_2 \cdot (-\alpha)$. So 
$$Y_{24}(y_3) = -\alpha.$$ 
Also as $\sigma_{y_2}$ is a localization of $\cM/Gy_2$, then by Lemma~\ref{loc.add}
\begin{equation}\label{eq4}
 1 -\sigma(Y_{12}) \alpha -\sigma(Y_{24}) \in N_G. 
\end{equation}

We still consider $\cM/Gy_2$. $(Y_{23}\backslash Gy_2, Y_{12}\backslash Gy_2 \cdot (-\sigma(Y_{12})^{-1}))$ forms a modular pair in $\cC^*(\cM/Gy_2)$ and $\sigma_{y_2}(Y_{23}\backslash Gy_2) = 1 = - \sigma_{y_2}(Y_{12}\backslash Gy_2 \cdot (-\sigma(Y_{12})^{-1}))$. 
We apply Lemma~\ref{mod.elim} with choosing $Y_1 = Y_{23}\backslash Gy_2$, $Y_2 = Y_{12}\backslash Gy_2 \cdot (-\sigma(Y_{12})^{-1})$, $e_1 = y_3$, $e= y_4$ and $Y_e = Y_{24}\backslash Gy_2$, and get
\begin{align*}
\Mod(Y_1, Y_2, p)(y_4) & = - Y_{23}(y_4) \cdot \sigma(Y_{23})^{-1}\cdot \sigma(Y_{24})\cdot Y_{24}(y_3)^{-1} \cdot (Y_{12}(y_3)\cdot (- \sigma(Y_{12})^{-1}))\\
& = -Y_{23}(y_4) \cdot \sigma(Y_{24})\cdot \alpha^{-1} \cdot  \sigma(Y_{12})^{-1}.
\end{align*}
Similarly, we apply Lemma~\ref{mod.elim} with choosing $Y_1' = Y_{12}\backslash Gy_2 \cdot (-\sigma(Y_{12})^{-1})$, $Y_2' = Y_{23}\backslash Gy_2$, $e_1 = y_1$, $e= y_4$ and $Y_e = Y_{24}\backslash Gy_2$, and get
\begin{align*}
\Mod(Y_1', Y_2', p)(y_4) & = -(-Y_{12}(y_4) \cdot \sigma(Y_{12})^{-1}) \cdot (-1) \cdot \sigma(Y_{24})\cdot Y_{24}(y_1)^{-1} \cdot Y_{23}(y_1) \\
& = -Y_{12}(y_4) \cdot \sigma(Y_{12})^{-1}\cdot \sigma(Y_{24}).
\end{align*}
As $\Mod(Y_1, Y_2, p)(y_4) = \Mod(Y_1', Y_2', p)(y_4)$, then
$$-Y_{23}(y_4) \cdot \sigma(Y_{24})\cdot \alpha^{-1} \cdot  \sigma(Y_{12})^{-1} = -Y_{12}(y_4) \cdot \sigma(Y_{12})^{-1}\cdot \sigma(Y_{24}).$$
Multiplying on the left by $-Y_{23}(y_4)^{-1}$ and we get
\begin{align*}
\sigma(Y_{24})\cdot \alpha^{-1} \cdot \sigma(Y_{12})^{-1} & = \alpha^{-1} \cdot \sigma(Y_{12})^{-1} \cdot \sigma(Y_{24}). \\
\sigma(Y_{24})\cdot (\sigma(Y_{12}) \cdot \alpha)^{-1} & = (\sigma(Y_{12})\cdot \alpha)^{-1} \cdot \sigma(Y_{24}). \numberthis \label{eq8}
\end{align*}

Next, let $\beta\in T^\times$ with $\beta = Y_{13}(y_4)^{-1} Y_{23}(y_4)$. We consider the rank 2 contraction $\cM/Gy_3$. Let $\sigma_{y_3}: \cC^*(\cM/Gy_3) \rightarrow T$ be the function induced from $\sigma$. As $\sigma$ is a localization on the rank 2 contraction $\cM/Gy_3$, then $\sigma_{y_3}$ is a localization of $\cM/Gy_3$. Now we consider the $T$-cocircuits $-Y_{23}\backslash Gy_3, Y_{13}\backslash Gy_3 \cdot \beta$ and $Y_{34}\backslash Gy_3$. We list their values on elements $\{y_1, y_2, y_4\}$ and the corresponding $\sigma_{y_3}$ values in Table~\ref{table3}.

\begin{table}[htb]
\centering
\begin{tabular}{ |c | c |c|c|c| }\hline
& $y_1$ & $y_2$ & $y_4$ & $\sigma_{y_3}$ \\ \hline 
 $-Y_{23}\backslash Gy_3$ & -1 & 0 & $-Y_{23}(y_4)$ & -1 \\ \hline
 $Y_{13}\backslash Gy_3 \cdot \beta$ & 0 & $\beta$ & $Y_{23}(y_4)$ & $ \sigma(Y_{13}) \cdot \beta$ \\ \hline
$Y_{34}\backslash Gy_3$ & -1 & $Y_{34}(y_2)$ &0 & $\sigma(Y_{34})$ \\ \hline
\end{tabular}
\smallskip
\caption{}
\label{table3}
\end{table}

Note that $Y_{34}\backslash Gy_3$ eliminates $y_4$ between $-Y_{23}\backslash Gy_3$ and $Y_{13}\backslash Gy_3 \cdot \beta$. So 
$$Y_{34}(y_2) =\beta.$$ 
Also as $\sigma_{y_3}$ is a localization of $\cM/Gy_3$, then by Lemma~\ref{loc.add}
\begin{equation}\label{eq5}
-1 + \sigma(Y_{13}) \beta - \sigma(Y_{34}) \in N_G.
\end{equation}

We still consider $\cM/Gy_3$. $(-Y_{23}\backslash Gy_3, Y_{13}\backslash Gy_3 \cdot \sigma(Y_{13})^{-1})$ forms a modular pair in $\cC^*(\cM/Gy_3)$ and $\sigma_{y_3}(-Y_{23}\backslash Gy_3) = -1 = - \sigma_{y_3}(Y_{13}\backslash Gy_3 \cdot \sigma(Y_{13})^{-1})$. 
We apply Lemma~\ref{mod.elim} with choosing $Y_1 = -Y_{23}\backslash Gy_3$, $Y_2 = Y_{13}\backslash Gy_3 \cdot \sigma(Y_{13})^{-1}$, $e_1 = y_2$, $e= y_4$ and $Y_e = Y_{34}\backslash Gy_3$, and get
\begin{align*}
\Mod(Y_1, Y_2, p)(y_4) & = -(-Y_{23}(y_4))\cdot (-\sigma(Y_{23})^{-1})\cdot \sigma(Y_{34})\cdot Y_{34}(y_2)^{-1} \cdot (Y_{13}(y_2) \cdot \sigma(Y_{13})^{-1})\\
& = - Y_{23}(y_4) \cdot \sigma(Y_{34})\cdot \beta^{-1} \cdot \sigma(Y_{13})^{-1}.
\end{align*}
Similarly, we apply Lemma~\ref{mod.elim} with choosing $Y_1' =Y_{13}\backslash Gy_3 \cdot \sigma(Y_{13})^{-1}$, $Y_2' = -Y_{23}\backslash Gy_3$, $e_1 = y_1$, $e= y_4$ and $Y_e = Y_{34}\backslash Gy_3$, and get
\begin{align*}
\Mod(Y_1', Y_2', p)(y_4) & = - (Y_{13}(y_4) \cdot \sigma(Y_{13})^{-1}) \cdot 1 \cdot \sigma(Y_{34})\cdot Y_{34}(y_1)^{-1} \cdot (-Y_{23}(y_1)) \\
& = - Y_{13}(y_4) \cdot \sigma(Y_{13})^{-1}\cdot \sigma(Y_{34}).
\end{align*}
As $\Mod(Y_1, Y_2, p)(y_4) = \Mod(Y_1', Y_2', p)(y_4)$, then
$$- Y_{23}(y_4) \cdot \sigma(Y_{34})\cdot \beta^{-1} \cdot \sigma(Y_{13})^{-1} = - Y_{13}(y_4) \cdot \sigma(Y_{13})^{-1}\cdot \sigma(Y_{34}).$$
Multiplying on the left by $-Y_{23}(y_4)^{-1}$ and we get
\begin{align*}
\sigma(Y_{34}) \cdot \beta^{-1} \cdot \sigma(Y_{13})^{-1} & = \beta^{-1} \cdot \sigma(Y_{13})^{-1} \cdot \sigma(Y_{34}). \\
\sigma(Y_{34}) \cdot (\sigma(Y_{13})\cdot \beta)^{-1} &= (\sigma(Y_{13})\cdot \beta)^{-1} \cdot \sigma(Y_{34}). \numberthis\label{eq9}
\end{align*}

Next, we consider the rank 2 contraction $\cM/Gy_1$. Let $\sigma_{y_1}: \cC^*(\cM/Gy_1) \rightarrow T$ be the function induced from $\sigma$. As $\sigma$ is a localization on the rank 2 contraction $\cM/Gy_1$, then $\sigma_{y_1}$ is a localization of $\cM/Gy_1$. Now we consider the $T$-cocircuits $Y_{13}\backslash Gy_1\cdot \beta $, $Y_{12}\backslash Gy_1 \cdot (-\alpha)$ and $Y_{14}\backslash Gy_1$. We list their values on elements $\{y_2, y_3, y_4\}$ and the corresponding $\sigma_{y_1}$ values in Table~\ref{table4}.

\begin{table}[htb]
\centering
\begin{tabular}{ |c | c |c|c|c|c| }\hline
& $y_2$ & $y_3$ & $y_4$ & $\sigma_{y_1}$ \\ \hline 
$Y_{13}\backslash Gy_1\cdot \beta$  & $\beta$ & 0 & $Y_{23}(y_4)$ & $\sigma(Y_{13}) \cdot \beta$ \\\hline
 $Y_{12}\backslash Gy_1 \cdot (-\alpha)$ & 0 & $-\alpha$ & $-Y_{23}(y_4)$ & $\sigma(Y_{12}) \cdot (-\alpha)$  \\ \hline
$Y_{14}\backslash Gy_1$  & $\beta$ & $Y_{14}(y_3)$ &0 & $\sigma(Y_{14})$  \\ \hline
\end{tabular}
\smallskip
\caption{}
\label{table4}
\end{table}

Note that $Y_{14}\backslash Gy_1$ eliminates $y_4$ between $Y_{13}\backslash Gy_1\cdot \beta $ and $Y_{12}\backslash Gy_1 \cdot (-\alpha)$. So
$$Y_{14}(y_3) = -\alpha.$$ 
Also as $\sigma_{y_1}$ is a localization of $\cM/Gy_1$, then by Lemma~\ref{loc.add}
\begin{equation}\label{eq6}
\sigma(Y_{13})\beta -\sigma(Y_{12})\alpha -\sigma(Y_{14}) \in N_G.
\end{equation}

We still consider $\cM/Gy_1$. $(Y_{13}\backslash Gy_1 \cdot \sigma(Y_{13})^{-1}, Y_{12}\backslash Gy_1 \cdot (-\sigma(Y_{12})^{-1}))$ forms a modular pair in $\cC^*(\cM/Gy_1)$ and $\sigma_{y_1}(Y_{13}\backslash Gy_1 \cdot \sigma(Y_{13})^{-1}) = 1 = - \sigma_{y_1}(Y_{12}\backslash Gy_1 \cdot (-\sigma(Y_{12})^{-1}))$. 
We apply Lemma~\ref{mod.elim} with choosing $Y_1 = Y_{13}\backslash Gy_1 \cdot \sigma(Y_{13})^{-1}$, $Y_2 = Y_{12}\backslash Gy_1 \cdot (-\sigma(Y_{12})^{-1})$, $e_1 = y_3$, $e= y_4$ and $Y_e = Y_{14}\backslash Gy_1$, and get
\begin{align*}
\Mod(Y_1, Y_2, p)(y_4) & = - (Y_{13}(y_4)\cdot \sigma(Y_{13})^{-1}) \cdot 1 \cdot \sigma(Y_{14})\cdot Y_{14}(y_3)^{-1} \cdot (-Y_{12}(y_3) \cdot \sigma(Y_{12})^{-1})\\
& = - Y_{13}(y_4)\cdot \sigma(Y_{13})^{-1} \cdot \sigma(Y_{14})\cdot \alpha^{-1} \cdot \sigma(Y_{12})^{-1}.
\end{align*}
Similarly, we apply Lemma~\ref{mod.elim} with choosing $Y_1' = Y_{12}\backslash Gy_1 \cdot (-\sigma(Y_{12})^{-1})$, $Y_2' = Y_{13}\backslash Gy_1 \cdot \sigma(Y_{13})^{-1}$, $e_1 = y_2$, $e= y_4$ and $Y_e = Y_{14}\backslash Gy_1$, and get
\begin{align*}
\Mod(Y_1', Y_2', p)(y_4) & = - (-Y_{12}(y_4)\cdot \sigma(Y_{12})^{-1}) \cdot (-1) \cdot \sigma(Y_{14})\cdot Y_{14}(y_2)^{-1} \cdot (Y_{13}(y_2) \cdot \sigma(Y_{13})^{-1})\\
& = - Y_{12}(y_4)\cdot \sigma(Y_{12})^{-1} \cdot \sigma(Y_{14})\cdot \beta^{-1} \cdot \sigma(Y_{13})^{-1}.
\end{align*}
As $\Mod(Y_1, Y_2, p)(y_4) = \Mod(Y_1', Y_2', p)(y_4)$, then
$$ - Y_{13}(y_4)\cdot \sigma(Y_{13})^{-1} \cdot \sigma(Y_{14})\cdot \alpha^{-1} \cdot \sigma(Y_{12})^{-1} = - Y_{12}(y_4)\cdot \sigma(Y_{12})^{-1} \cdot \sigma(Y_{14})\cdot \beta^{-1} \cdot \sigma(Y_{13})^{-1}.$$
Multiplying on the left by $-Y_{23}(y_4)^{-1}$ and we get
\begin{align*}
\beta^{-1}\cdot \sigma(Y_{13})^{-1} \cdot \sigma(Y_{14})\cdot \alpha^{-1} \cdot \sigma(Y_{12})^{-1} & = \alpha^{-1}\cdot \sigma(Y_{12})^{-1} \cdot \sigma(Y_{14})\cdot \beta^{-1} \cdot \sigma(Y_{13})^{-1}. \\
(\sigma(Y_{13}) \cdot \beta)^{-1} \cdot \sigma(Y_{14}) \cdot (\sigma(Y_{12}) \cdot \alpha)^{-1} & = (\sigma(Y_{12}) \cdot \alpha)^{-1} \cdot \sigma(Y_{14})\cdot (\sigma(Y_{13}) \cdot \beta)^{-1}. \numberthis\label{eq10}
\end{align*}

Next, we consider the rank 2 contraction $\cM/Gy_4$. Let $\sigma_{y_4}: \cC^*(\cM/Gy_4) \rightarrow T$ be the function induced from $\sigma$. As $\sigma$ is a localization on the rank 2 contraction $\cM/Gy_4$, then $\sigma_{y_4}$ is a localization of $\cM/Gy_4$. Now we consider the $T$-cocircuits $Y_{24}\backslash Gy_4, Y_{34}\backslash Gy_4$ and $Y_{14}\backslash Gy_4$. We list their values on elements $\{y_1, y_2, y_3\}$ and the corresponding $\sigma_{y_4}$ values in Table~\ref{table5}.

\begin{table}[htb]
\centering
\begin{tabular}{ |c | c |c|c|c| }\hline
& $y_1$ & $y_2$ & $y_3$ & $\sigma_{y_4}$   \\ \hline 
 $Y_{24}\backslash Gy_4$ & 1 & 0 & $Y_{24}(y_3)$   & $\sigma(Y_{24})$ \\ \hline
 $Y_{34}\backslash Gy_4$ & -1 & $Y_{34}(y_2)$ &0  & $\sigma(Y_{34})$ \\ \hline
 $Y_{14}\backslash Gy_4$ & 0 & $\beta$ & $Y_{14}(y_3)$ & $\sigma(Y_{14})$ \\ \hline
\end{tabular}
\smallskip
\caption{}
\label{table5}
\end{table}

As $Y_{34}(y_2) = \beta$, then $Y_{14}\backslash Gy_4$ eliminates $y_1$ between $Y_{24}\backslash Gy_4$ and $Y_{34}\backslash Gy_4$. So
$$Y_{14}(y_3) = Y_{24}(y_3) = -\alpha.$$

As $\sigma_{y_4}$ is a localization of $\cM/Gy_4$, then by Lemma~\ref{loc.add} 
\begin{equation}\label{eq7}
 \sigma(Y_{24}) + \sigma(Y_{34}) -\sigma(Y_{14}) \in N_G.
\end{equation}

We still consider $\cM/Gy_4$. $(Y_{24}\backslash Gy_4 \cdot \sigma(Y_{24})^{-1}, Y_{34}\backslash Gy_4 \cdot (-\sigma(Y_{34})^{-1}))$ forms a modular pair in $\cC^*(\cM/Gy_4)$ and $\sigma_{y_4}(Y_{24}\backslash Gy_4 \cdot \sigma(Y_{24})^{-1}) = 1 = - \sigma_{y_4}(Y_{34}\backslash Gy_4 \cdot (-\sigma(Y_{34})^{-1}))$. 
We apply Lemma~\ref{mod.elim} with choosing $Y_1 = Y_{24}\backslash Gy_4 \cdot \sigma(Y_{24})^{-1}$, $Y_2 = Y_{34}\backslash Gy_4 \cdot (-\sigma(Y_{34})^{-1})$, $e_1 = y_2$, $e= y_1$ and $Y_e = Y_{14}\backslash Gy_4$, and get
\begin{align*}
\Mod(Y_1, Y_2, p)(y_1) & = - (Y_{24}(y_1)\cdot \sigma(Y_{24})^{-1}) \cdot 1 \cdot \sigma(Y_{14})\cdot Y_{14}(y_2)^{-1} \cdot (-Y_{34}(y_2) \cdot \sigma(Y_{34})^{-1})\\
& = - \sigma(Y_{24})^{-1} \cdot \sigma(Y_{14})\cdot \beta^{-1} \cdot (-\beta \cdot \sigma(Y_{34})^{-1}) \\
& = \sigma(Y_{24})^{-1} \cdot \sigma(Y_{14})\cdot \sigma(Y_{34})^{-1}.
\end{align*}
Similarly, we apply Lemma~\ref{mod.elim} with choosing $Y_1' = Y_{34}\backslash Gy_4 \cdot (-\sigma(Y_{34})^{-1})$, $Y_2' = Y_{24}\backslash Gy_4 \cdot \sigma(Y_{24})^{-1}$, $e_1 = y_3$, $e= y_1$ and $Y_e = Y_{14}\backslash Gy_4$, and get
\begin{align*}
\Mod(Y_1', Y_2', p)(y_1) & = - (-Y_{34}(y_1)\cdot \sigma(Y_{34})^{-1}) \cdot (-1) \cdot \sigma(Y_{14})\cdot Y_{14}(y_3)^{-1} \cdot (Y_{24}(y_3) \cdot \sigma(Y_{24})^{-1})\\
& = \sigma(Y_{34})^{-1} \cdot \sigma(Y_{14})\cdot (-\alpha)^{-1} \cdot (-\alpha \cdot \sigma(Y_{24})^{-1}) \\
& = \sigma(Y_{34})^{-1} \cdot \sigma(Y_{14})\cdot \sigma(Y_{24})^{-1}.
\end{align*}
As $\Mod(Y_1, Y_2, p)(y_1) = \Mod(Y_1', Y_2', p)(y_1)$, then we get
\begin{equation} \label{eq11}
\sigma(Y_{24})^{-1} \cdot \sigma(Y_{14})\cdot \sigma(Y_{34})^{-1} = \sigma(Y_{34})^{-1} \cdot \sigma(Y_{14})\cdot \sigma(Y_{24})^{-1}.
\end{equation}

So now we can revise our previous Table~\ref{table1} to Table~\ref{table.main}.
\begin{table}[htb]
\centering
\begin{tabular}{ | l | c |c|c|c| c|} \hline 
& $y_1$ & $y_2$ & $y_3$ & $y_4$ & $\sigma$  \\ \hline 
 $Y_{23}$ & 1& 0&0 & $Y_{23}(y_4)$& 1 \\ \hline 
 $Y_{13}$ & 0 & 1 & 0 & $Y_{13}(y_4)$ & $\sigma(Y_{13})$ \\ \hline 
 $Y_{12}$ & 0& 0& 1& $Y_{12}(y_4)$ &$\sigma(Y_{12})$ \\ \hline 
 $Y_{24}$ &1& 0& $-\alpha$& 0 & $\sigma(Y_{24})$ \\ \hline 
 $Y_{34}$ &-1 &$\beta$ &0 &0 &$\sigma(Y_{34})$ \\ \hline 
 $Y_{14}$ &0 & $\beta$ & $-\alpha$ &0 & $\sigma(Y_{14})$ \\ \hline 
\end{tabular}
\smallskip
\caption{}
\label{table.main}
\end{table}

Then from Formula~\eqref{ext.phi.form.thm} and by values in Table~\ref{table.main}, we have
\begin{align*}
[Fy_2y_4, Fy_1y_2]_p & = -\overline{Y_{12}(y_4) \sigma(Y_{12})^{-1}\sigma(Y_{24}) Y_{24}(y_1)^{-1}} = - \overline{Y_{12}(y_4) \sigma(Y_{12})^{-1}\sigma(Y_{24})}, \\
[Fy_1y_3, Fy_3y_4]_p & = -\overline{Y_{34}(y_1) \sigma(Y_{34})^{-1}\sigma(Y_{13}) Y_{13}(y_4)^{-1}} = \overline{\sigma(Y_{34})^{-1}\sigma(Y_{13}) Y_{13}(y_4)^{-1}}, \\
[Fy_1y_4, Fy_1y_2]_p & = -\overline{Y_{12}(y_4) \sigma(Y_{12})^{-1}\sigma(Y_{14}) Y_{14}(y_2)^{-1}} = -\overline{Y_{12}(y_4) \sigma(Y_{12})^{-1}\sigma(Y_{14}) \beta^{-1}}, \\
[Fy_2y_3, Fy_3y_4]_p & = -\overline{Y_{34}(y_2) \sigma(Y_{34})^{-1}\sigma(Y_{23}) Y_{23}(y_4)^{-1}} = -\overline{\beta \sigma(Y_{34})^{-1}Y_{23}(y_4)^{-1}}. 
\end{align*}

Thus the inclusion~\eqref{maineq} can be written as
\begin{align*}
& N_G \\
\ni & -1 +[Fy_2y_4, Fy_1y_2]_p \cdot [Fy_1y_3, Fy_3y_4]_p + [Fy_1y_4, Fy_1y_2]_p \cdot [Fy_2y_3, Fy_3y_4]_p\\
= & -1 + (- \overline{Y_{12}(y_4) \sigma(Y_{12})^{-1}\sigma(Y_{24})}) \cdot \overline{\sigma(Y_{34})^{-1}\sigma(Y_{13}) Y_{13}(y_4)^{-1}} + \\
& \hspace{2in} (-\overline{Y_{12}(y_4) \sigma(Y_{12})^{-1}\sigma(Y_{14}) \beta^{-1}}) \cdot (-\overline{\beta \sigma(Y_{34})^{-1}Y_{23}(y_4)^{-1}} ) \\
= & \overline{ -1 - Y_{12}(y_4) \sigma(Y_{12})^{-1}\sigma(Y_{24})\sigma(Y_{34})^{-1}\sigma(Y_{13}) Y_{13}(y_4)^{-1} + Y_{12}(y_4) \sigma(Y_{12})^{-1}\sigma(Y_{14})\sigma(Y_{34})^{-1}Y_{23}(y_4)^{-1}} \\
= &\overline {-Y_{12}(y_4) \sigma(Y_{12})^{-1} \cdot (\sigma(Y_{12}) Y_{12}(y_4)^{-1}Y_{23}(y_4)\sigma(Y_{34}) + \sigma(Y_{24})\sigma(Y_{34})^{-1}\sigma(Y_{13}) Y_{13}(y_4)^{-1}Y_{23}(y_4) \sigma(Y_{34}) - \sigma(Y_{14}))} \\
& \,\,\,\,\,\,\,\, \overline{\cdot \sigma(Y_{34})^{-1}Y_{23}(y_4)^{-1}}\\
= &\overline {-Y_{12}(y_4) \sigma(Y_{12})^{-1}} \cdot \overline{(\sigma(Y_{12}) \alpha \sigma(Y_{34})+ \sigma(Y_{24})\sigma(Y_{34})^{-1}\sigma(Y_{13}) \beta \sigma(Y_{34})- \sigma(Y_{14}))} \cdot \overline{\sigma(Y_{34})^{-1}Y_{23}(y_4)^{-1}}.
\end{align*}

Now we will use Pathetic Cancellation to see that
$$\sigma(Y_{12}) \alpha \sigma(Y_{34})+ \sigma(Y_{24})\sigma(Y_{34})^{-1}\sigma(Y_{13}) \beta \sigma(Y_{34}) - \sigma(Y_{14}) \in N_G,$$
from which the inclusion~\eqref{maineq} will hold.

Let $$a = - \sigma(Y_{12})\alpha, \, b = \sigma(Y_{13})\beta, \, x = \sigma(Y_{24}), \, y = \sigma(Y_{34}), \, z = \sigma(Y_{14}).$$
Then we can restate the inclusion~\eqref{eq4}, \eqref{eq8}, \eqref{eq5}, \eqref{eq9}, \eqref{eq6}, \eqref{eq10}, \eqref{eq7} and \eqref{eq11} as
\begin{align*}
 1 + a - x &\in N_G, \\
 xa^{-1} & = a^{-1} x \hspace{0.3in}\text{ implying }  \hspace{0.2in} ax = xa, \\
-1 + b - y &\in N_G,\\
yb^{-1} &  = b^{-1} y \hspace{0.3in} \text{ implying } \hspace{0.2in} by = yb,\\
b + a - z &\in N_G,\\
b^{-1}za^{-1} & = a^{-1}zb^{-1},\\ 
x + y - z &\in N_G,\\
x^{-1}zy^{-1} & = y^{-1}zx^{-1}.
\end{align*}
All of these inclusions are exactly the hypothesis of Pathetic Cancellation property. So by Pathetic Cancellation,
$$xb - ay -z = x y^{-1} b y - ay - z = \sigma(Y_{24}) \sigma(Y_{34})^{-1} \sigma(Y_{13}) \beta \sigma(Y_{34}) + \sigma(Y_{12})\alpha \sigma(Y_{34}) -\sigma(Y_{14})\in N_G.$$

So the inclusion~\eqref{maineq} holds.
\end{enumerate}
\end{proof}

We can get the following Corollary to Lemma~\ref{loc.add}, Lemma~\ref{rank2.elim.loc} and Theorem~\ref{weakext.equi}.
 
\begin{cor} Let $\cM$ be a weak left (resp. right) $T$-matroid and let $\sigma: \mathcal{C^*(M)} \rightarrow T$ be a right (resp. left) $T^\times$-equivariant function. Then $\sigma$ is a weak localization if and only if for any modular pair $(Y_1, Y_2)$ of $T$-cocircuits with $Y_1(e) = -Y_2(e) \neq 0$ for some $e\in E$ and $Y_3$ eliminating $e$ between $Y_1$ and $Y_2$, we have $\sigma(Y_1) + \sigma(Y_2) - \sigma(Y) \in N_G$.
\end{cor}

This concludes the proof of our main theorem (Theorem~\ref{maintheorem}) in this paper.

\section{Stringent Skew Hyperfields}\label{sect.DD}

Stringent skew hyperfields behave in many ways like skew fields and the notions of weak and strong $F$-matroid coincide if $F$ is stringent (\cite{BP19}). This leads to a strengthening of Theorem~\ref{weakext.equi}, given as follows. 

\begin{thm}\label{DDext.equi}
Let $F$ be a stringent skew hyperfield, let $\cM$ be a strong left (resp. right) $F$-matroid on $E$ with $F$-cocircuit set $\mathcal{C^*(M)}$, and let
$$\sigma: \mathcal{C^*(M)} \rightarrow F$$ 
be a right (resp. left) $F^\times$-equivariant function. 

Then the following statements are equivalent.
\begin{enumerate}[(1)]
\item $\sigma$ defines a strong single-element extension of $\cM$.

\item $\sigma$ defines a strong single-element extension of every rank 2 contraction of $\cM$.

\item $\sigma$ defines a strong single-element extension of every rank 2 minor of $\cM$ on three elements.
\end{enumerate}
\end{thm}

\subsection{Characterization of strong Pathetic Cancellation for skew hyperfields}

In this section, we will first show a stronger property that implies Pathetic Cancellation for skew hyperfields. 

\begin{defn} Let $F$ be a skew hyperfield. We say $F$ satisfies the {\bf strong Pathetic Cancellation Property} if for every $a, b, x, y \in F^{\times}$ with $x\in 1\boxplus a$ and $y\in -1 \boxplus b$ we have
$$(x\boxplus y)\cap (a \boxplus b) \subseteq xb \boxplus -ay,$$
\end{defn}

It is easy to see that the strong Pathetic Cancellation Property and the Pathetic Cancellation Property will be the same when $F$ is a hyperfield.

\begin{thm}\label{strongPCimplyPC}
Let $F$ be a skew hyperfield. If $F$ satisfies strong Pathetic Cancellation,
then $F$ satisfies Pathetic Cancellation. 
\end{thm}

\begin{proof} Let $a, b, x, y, z\in F^\times$ with $x\in 1\boxplus a$, $ax=xa$, $y\in -1 \boxplus b$, $by = yb$, $z\in a\boxplus b$, $a^{-1}zb^{-1} = b^{-1}za^{-1}$, $z\in x\boxplus y$ and $x^{-1}zy^{-1} = y^{-1}zx^{-1}$. As $F$ satisfies strong Pathetic Cancellation, then
$$z\in (x\boxplus y)\cap (a \boxplus b) \subseteq xb \boxplus -ay.$$
So $F$ satisfies Pathetic Cancellation. 
\end{proof}

\begin{prop}\label{intersect.equal}
Let $F$ be a skew hyperfield. The following are equivalent.
\begin{enumerate}[(1)]
\item $F$ satisfies strong Pathetic Cancellation. 
\item For any $a, b \in F^{\times}$ and $x$, $y$ in $F$ with $x\in 1\boxplus a$, $y\in -1 \boxplus b$, we have
$$(x\boxplus y)\cap (a \boxplus b) \subseteq xb \boxplus -ay,$$
$$(a\boxplus b)\cap (xb \boxplus -ay) \subseteq x\boxplus y,$$
$$(x\boxplus y)\cap (xb \boxplus -ay) \subseteq a\boxplus b.$$
\item For any $a, b \in F^{\times}$ and $x$, $y$ in $F$ with $x\in 1\boxplus a$, $y\in -1 \boxplus b$, we have
$$(x\boxplus y)\cap (a \boxplus b)=(x\boxplus y)\cap (xb \boxplus -ay)=(a\boxplus b)\cap (xb \boxplus -ay).$$
\end{enumerate}
\end{prop}

\begin{proof} (1)$\Rightarrow$(3): By definition, for any $x\in 1\boxplus a$, $y\in -1 \boxplus b$, we have
$$(x\boxplus y)\cap (a \boxplus b) \subseteq xb \boxplus -ay.$$

Then 
\begin{equation}\label{intersect.eq1}
(x\boxplus y)\cap (a \boxplus b) \subseteq (xb \boxplus -ay)\cap (x\boxplus y)
\end{equation}

and
\begin{equation}\label{intersect.eq2}
(x\boxplus y)\cap (a \boxplus b) \subseteq (xb \boxplus -ay)\cap (a\boxplus b).
\end{equation}

As $x\in 1\boxplus a$ and $y\in -1 \boxplus b$, then $a^{-1}x \in 1 \boxplus a^{-1}$ and $-yb^{-1} \in -1 \boxplus b^{-1}$. By strong Pathetic Cancellation, we have
$$(a^{-1} \boxplus b^{-1})\cap (a^{-1}x \boxplus -yb^{-1}) \subseteq a^{-1}xb^{-1} \boxplus a^{-1}yb^{-1} = a^{-1} (x\boxplus y) b^{-1},$$
thus
$$(a\boxplus b)\cap (xb \boxplus -ay) \subseteq x\boxplus y.$$

Then
\begin{equation}\label{intersect.eq3}
(a\boxplus b)\cap (xb \boxplus -ay) \subseteq (x\boxplus y)\cap (a \boxplus b).
\end{equation}

and
\begin{equation}\label{intersect.eq4}
(a\boxplus b)\cap (xb \boxplus -ay) \subseteq (x\boxplus y)\cap (xb \boxplus -ay).
\end{equation}

Similarly, as $x\in 1\boxplus a$ and $y\in -1 \boxplus b$, then $x^{-1}a \in 1 \boxplus -x^{-1}$ and $-by^{-1} \in -1 \boxplus -y^{-1}$. By strong Pathetic Cancellation, we have
$$(-x^{-1} \boxplus -y^{-1})\cap ( x^{-1}a\boxplus -by^{-1}) \subseteq -x^{-1}ay^{-1} \boxplus -x^{-1}by^{-1} = -x^{-1}(a\boxplus b) y^{-1}, $$
thus
$$(x\boxplus y)\cap (xb \boxplus -ay) \subseteq a\boxplus b.$$

Then
\begin{equation}\label{intersect.eq5}
(x\boxplus y)\cap (xb \boxplus -ay) \subseteq (a\boxplus b)\cap (x \boxplus y).
\end{equation}
and
\begin{equation}\label{intersect.eq6}
(x\boxplus y)\cap (xb \boxplus -ay) \subseteq (a\boxplus b)\cap (xb \boxplus -ay).
\end{equation}

So by Inclusion~\eqref{intersect.eq1}, \eqref{intersect.eq2}, \eqref{intersect.eq3}, \eqref{intersect.eq4}, \eqref{intersect.eq5} and \eqref{intersect.eq6}, we have
$$(x\boxplus y)\cap (a \boxplus b)=(x\boxplus y)\cap (xb \boxplus -ay)=(a\boxplus b)\cap (xb \boxplus -ay).$$

(3)$\Rightarrow$(2): If for any $x$, $y$ in $F$ with $x\in 1\boxplus a$, $y\in -1 \boxplus b$, we have
$$(x\boxplus y)\cap (a \boxplus b)=(x\boxplus y)\cap (xb \boxplus -ay)=(a\boxplus b)\cap (xb \boxplus -ay).$$

Then $$(x\boxplus y)\cap (a \boxplus b)=(x\boxplus y)\cap (xb \boxplus -ay)\subseteq xb \boxplus -ay,$$
$$(a\boxplus b)\cap (xb \boxplus -ay) = (x\boxplus y)\cap (xb \boxplus -ay) \subseteq x\boxplus y,$$
$$(x\boxplus y)\cap (xb \boxplus -ay)=(a\boxplus b)\cap (xb \boxplus -ay)\subseteq a\boxplus b.$$

(2)$\Rightarrow$(1): This is immediate by definition of stong Pathetic Cancellation.
\end{proof}

\begin{cor} \label{PP.multi} For every skew hyperfield satisfying strong Pathetic Cancellation, we have
\begin{equation}\label{eq12} 
1\boxplus -1 \boxplus 1 \boxplus -1 = 1\boxplus -1.
\end{equation}
\end{cor}

\begin{proof} We first prove the direction $(\subseteq)$. Let $f\in 1\boxplus -1 \boxplus 1 \boxplus -1$. Then there exist $x, y$ such that $x, y \in 1\boxplus -1$ and $f\in x\boxplus y$. Let $a=-1$ and $b=1$. Then $x\in 1\boxplus a$ and $ y \in -1 \boxplus b$. Then by By Proposition~\ref{intersect.equal},
$$x\boxplus y = (x\boxplus y)\cap (xb\boxplus -ay) \subseteq a\boxplus b = 1\boxplus -1.$$
So $f\in x\boxplus y = 1\boxplus -1$, and so 
$$1\boxplus -1 \boxplus 1 \boxplus -1 \subseteq 1\boxplus -1.$$

Now we prove the direction $(\supseteq)$. It is clear that
$$1\boxplus -1 \subseteq 1\boxplus -1 \boxplus 0 \subseteq 1\boxplus -1 \boxplus 1 \boxplus -1.$$

So
$$1\boxplus -1 \boxplus 1 \boxplus -1 = 1\boxplus -1.$$
\end{proof}


\subsection{Pathetic Cancellation for stringent skew hyperfields}\label{subsect.DD.PC}

In this section, we will show that the stringency implies the strong Pathetic Cancellation property and so implies Pathetic Cancellation. This helps to generalize the characterization of localization of strong matroids over stringent skew hyperfields.

\begin{defn}\label{stringent} (\cite{BS20})
A skew hyperfield $F$ is said to be \textbf{stringent} if for any $a,b\in F$, $a\boxplus b$ is a singleton whenever $a\neq -b$.
\end{defn}

Fields, $\mathbb{K}$, and $\mathbb{S}$ are all stringent, but $\mathbb{P}$ does not satisfy stringency. 

\begin{thm}\label{DDPC} Any stringent skew hyperfield satisfies strong Pathetic Cancellation and Pathetic Cancellation.
\end{thm}

Before showing Theorem~\ref{DDPC}, we will first introduce a useful lemma for any skew hyperfields.

\begin{lem}\label{intersect} Let $F$ be a skew hyperfield and let $a, b\in F$. For any $x$, $y$ in $F$ with $x\in 1\boxplus a$, $y\in -1 \boxplus b$, we have 
$$(x\boxplus y)\cap (a \boxplus b) \neq \emptyset,$$
$$(x\boxplus y)\cap (xb \boxplus -ay) \neq \emptyset,$$
$$(a\boxplus b)\cap (xb \boxplus -ay) \neq \emptyset.$$
\end{lem}

\begin{proof} As $x\in 1\boxplus a$ and $y\in -1 \boxplus b$, we have $1\in x\boxplus -a$ and $-1\in y\boxplus -b$. Thus
$$0\in 1\boxplus -1 \subseteq x\boxplus -a \boxplus y\boxplus -b = (x\boxplus y)\boxplus -(a\boxplus b).$$
So $$(x\boxplus y)\cap (a \boxplus b) \neq \emptyset.$$

As $x\in 1\boxplus a$ and $y\in -1 \boxplus b$, we have $xy\in y\boxplus ay$ and $-xy\in x\boxplus -xb$. Thus
$$0\in xy\boxplus -xy \subseteq y\boxplus ay \boxplus x\boxplus -xb = (x\boxplus y)\boxplus -(xb\boxplus -ay).$$
So $$(x\boxplus y)\cap (xb \boxplus -ay) \neq \emptyset.$$

As $x\in 1\boxplus a$ and $y\in -1 \boxplus b$, we have $a\in -1\boxplus x$ and $b\in 1\boxplus y$. Then $ab\in -b\boxplus xb$ and $-ab\in -a\boxplus -ay$. Thus
$$0\in ab\boxplus -ab \subseteq -b\boxplus xb \boxplus -a\boxplus -ay = -(a\boxplus b)\boxplus (xb\boxplus -ay).$$
So $$(a\boxplus b)\cap (xb \boxplus -ay) \neq \emptyset.$$
\end{proof}

Now we prove Theorem~\ref{DDPC} by using the above lemma.

\begin{proof}[Proof of Theorem~\ref{DDPC}] Let $F = R\rtimes_{H, \psi} G$ be a stringent skew hyperfield. By Theorem~\ref{strongPCimplyPC}, we only need to show that $F$ satisfies strong Pathetic Cancellation. Let $a, b, x, y \in F^\times$ with $x\in 1\boxplus a $ and $ y\in -1\boxplus b,$. We need to show that
$$(a\boxplus b)\cap (x\boxplus y)\subseteq xb\boxplus -ay.$$

We will divide the proof into three cases.

\textbf{Case 1:} If $a=-1$ and $b=1$, then $xb\boxplus -ay = x\boxplus y$. So
$$(a\boxplus b)\cap (x\boxplus y) \subseteq x\boxplus y = xb\boxplus -ay.$$

\textbf{Case 2:} If $a=-1$ and $b\neq 1$ or $a\neq -1$ and $b = 1$, then without loss of generality we may assume that $a=-1$ and $b\neq 1$. So by stringency, $-1\boxplus b = \{y\}$. Then $a\boxplus b = -1\boxplus b = \{y\}$. By Lemma~\ref{intersect},
$$a\boxplus b = \{y\} \subseteq x\boxplus y \text{ and } a\boxplus b = \{y\} \subseteq xb\boxplus -ay.$$
So $$(a\boxplus b)\cap (x\boxplus y)=\{y\} \subseteq xb\boxplus -ay.$$

\textbf{Case 3:} If $a\neq -1$ and $b \neq 1$, then $1\boxplus a = \{x\}$ and $-1\boxplus b = \{y\}$. So $a^{-1} \boxplus 1 = \{a^{-1} \cdot x\} = \{x\cdot a^{-1}\}$, and so $xa = ax$.

If $a \neq -b$, then $a\boxplus b = \{z\}$ for some $z\in F^\times$. By Lemma~\ref{intersect},
$$a\boxplus b = \{z\} \subseteq x\boxplus y \text{ and } a\boxplus b = \{z\} \subseteq xb\boxplus -ay.$$
So $$(a\boxplus b)\cap (x\boxplus y)=\{z\} \subseteq xb\boxplus -ay.$$

If $a=-b$, then $x=-y$. As $1\boxplus a = \{x\}$, then $\psi(x) \geq 1_G$. So $\psi(ax) \geq \psi(a)$. Then
$$a\boxplus b= a(1\boxplus -1),$$
$$x\boxplus y= x(1\boxplus -1),$$
and $$xb\boxplus -ay = xa\boxplus -ax = ax \boxplus -ax = ax(1\boxplus -1).$$
So $$(a\boxplus b)\cap (x\boxplus y)\subseteq a(1\boxplus -1) \subseteq ax(1\boxplus -1) = xb\boxplus -ay.$$

So $F$ satisfies strong Pathetic Cancellation and also Pathetic Cancellation.
\end{proof}

From \cite{BP19}, we know the following theorem for matroids over stringent skew hyperfields.

\begin{thm}\label{ddweakstrong} (\cite{BP19}) Any weak matroid over a stringent skew hyperfield is strong.
\end{thm}

Then we can get a similar characterization of extension for a strong matroid over a stringent skew hyperfield as for a weak matroid, and so Theorem~\ref{DDext.equi} holds.

Next we would like to talk about the following important property for skew hyperfields.
\begin{defn} A skew hyperfield $F$ is {\bf doubly distributive} if for any $a, b, c, d \in F$, $(a\boxplus b) (c \boxplus d) = ac \boxplus ad \boxplus bc \boxplus bd$.
\end{defn}

Double distributivity is very important in matroid theory. Fields, $\mathbb{K}$, and $\mathbb{S}$ are all doubly distributive. 
In \cite{BS20}, Bowler and Su showed that any doubly distributive skew hyperfield is stringent. So the same statement in Theorem~\ref{DDext.equi} with `stringent' replaced by `doubly distributive' also holds.

\bibliographystyle{alpha}
\bibliography{biblioCM}

\newcommand{\etalchar}[1]{$^{#1}$}
\begin{thebibliography}{BLVS{\etalchar{+}}99}

\bibitem[AD12]{AD12}
Laura Anderson and Emanuele Delucchi.
\newblock Foundations for a theory of complex matroids.
\newblock {\em Discrete Comput. Geom.}, 48(4):807--846, 2012.

\bibitem[BB16]{Bak16}
Matthew Baker and Nathan Bowler.
\newblock Matroids over hyperfields.
\newblock {\em \tt arXiv:1601.01204}, 2016.

\bibitem[BB19]{Bak17}
Matthew Baker and Nathan Bowler.
\newblock Matroids over partial hyperstructures.
\newblock {\em Adv. Math.}, 343:821--863, 2019.

\bibitem[BLV78]{BLV78}
Robert~G. Bland and Michel Las~Vergnas.
\newblock Orientability of matroids.
\newblock {\em J. Combinatorial Theory Ser. B.}, 24(1):94--123, 1978.

\bibitem[BLVS{\etalchar{+}}99]{BLVS+99}
Anders Bj{\"o}rner, Michel Las~Vergnas, Bernd Sturmfels, Neil White, and
  G{\"u}nter~M. Ziegler.
\newblock {\em Oriented matroids}, volume~46 of {\em Encyclopedia of
  Mathematics and its Applications}.
\newblock Cambridge University Press, Cambridge, second edition, 1999.

\bibitem[BP19]{BP19}
Nathan Bowler and Rudi Pendavingh.
\newblock Perfect matroids over hyperfields.
\newblock {\em \tt arXiv:1908.03420}, 2019.

\bibitem[BS21]{BS20}
Nathan Bowler and Ting Su.
\newblock Classification of doubly distributive skew hyperfields and stringent
  hypergroups.
\newblock {\em J. Algebra}, 574:669--698, 2021.

\bibitem[CR70]{CR70}
Henry~H. Crapo and Gian-Carlo Rota.
\newblock {\em On the foundations of combinatorial theory: {C}ombinatorial
  geometries}.
\newblock The M.I.T. Press, Cambridge, Mass.-London, preliminary edition, 1970.

\bibitem[Cra65]{Crapo65}
Henry~H. Crapo.
\newblock Single-element extensions of matroids.
\newblock {\em J. Res. Nat. Bur. Standards Sect. B}, 69B:55--65, 1965.

\bibitem[Del11]{Del11}
Emanuele Delucchi.
\newblock Modular elimination in matroids and oriented matroids.
\newblock {\em European J. Combin.}, 32(3):339--343, 2011.

\bibitem[Dre86]{Dre86}
Andreas W.~M. Dress.
\newblock Duality theory for finite and infinite matroids with coefficients.
\newblock {\em Adv. in Math.}, 59(2):97--123, 1986.

\bibitem[FL78]{FL78}
Jon Folkman and Jim Lawrence.
\newblock Oriented matroids.
\newblock {\em Journal of Combinatorial Theory, Series B}, 25(2):199 -- 236,
  1978.

\bibitem[Kra57]{Kra57}
Marc Krasner.
\newblock Approximation des corps valu\'es complets de caract\'eristique
  {$p\not=0$} par ceux de caract\'eristique {$0$}.
\newblock In {\em Colloque d'alg\`ebre sup\'erieure, tenu \`a {B}ruxelles du 19
  au 22 d\'ecembre 1956}, Centre Belge de Recherches Math\'ematiques, pages
  129--206. \'Etablissements Ceuterick, Louvain; Librairie Gauthier-Villars,
  Paris, 1957.

\bibitem[Kra83]{Kra83}
Marc Krasner.
\newblock A class of hyperrings and hyperfields.
\newblock {\em Internat. J. Math. Math. Sci.}, 6(2):307--311, 1983.

\bibitem[LV75]{LV75}
Michel Las~Vergnas.
\newblock Matro\"\i des orientables.
\newblock {\em C. R. Acad. Sci. Paris S\'er. A-B}, 280:Ai, A61--A64, 1975.

\bibitem[LV78]{LV78}
Michel Las~Vergnas.
\newblock Extensions ponctuelles d'une g\'eom\'etrie combinatoire orient\'ee.
\newblock In {\em Probl\`emes combinatoires et th\'eorie des graphes ({C}olloq.
  {I}nternat. {CNRS}, {U}niv. {O}rsay, {O}rsay, 1976)}, volume 260 of {\em
  Colloq. Internat. CNRS}, pages 265--270. CNRS, Paris, 1978.

\bibitem[Oxl92]{Oxl92}
James~G. Oxley.
\newblock {\em Matroid theory}.
\newblock Oxford Science Publications. The Clarendon Press, Oxford University
  Press, New York, 1992.

\bibitem[Pen18]{Pen18}
Rudi Pendavingh.
\newblock Field extensions, derivations, and matroids over skew hyperfields.
\newblock {\em \tt arXiv:1802.02447}, 2018.

\bibitem[Sta12]{Sta12}
Richard~P. Stanley.
\newblock {\em Enumerative combinatorics. {V}olume 1}, volume~49 of {\em
  Cambridge Studies in Advanced Mathematics}.
\newblock Cambridge University Press, Cambridge, second edition, 2012.

\bibitem[Su18]{Su18}
Ting Su.
\newblock {\em Extensions of Matroids over Tracts and Doubly Distributive
  Partial Hyperfields}.
\newblock PhD thesis, (Order No. 10973187, State University of New York at
  Binghamton), 2018.
\newblock ProQuest Dissertations and Theses.

\bibitem[Su23]{Su20}
Ting Su.
\newblock Matroids over skew tracts.
\newblock {\em European Journal of Combinatorics}, 109:103643, 2023.

\bibitem[SW96]{SW96}
Charles Semple and Geoff Whittle.
\newblock Partial fields and matroid representation.
\newblock {\em Adv. in Appl. Math.}, 17(2):184--208, 1996.

\bibitem[Vir10]{Viro10}
Oleg Viro.
\newblock Hyperfields for tropical geometry i. hyperfields and dequantization.
\newblock {\em \tt arXiv:1006.3034v2}, 2010.

\bibitem[Wel76]{Wel76}
D.~J.~A. Welsh.
\newblock {\em Matroid theory}.
\newblock Academic Press [Harcourt Brace Jovanovich, Publishers], London-New
  York, 1976.
\newblock L. M. S. Monographs, No. 8.

\end{thebibliography}

\end{document}